\documentclass[10pt,a4paper]{article}
\usepackage[latin1]{inputenc}
\usepackage{amsmath}
\usepackage{amsfonts}
\usepackage{amssymb}
\usepackage[english]{babel}
\usepackage{makeidx}
\usepackage{enumerate}
\usepackage{fancyhdr}
\usepackage[all]{xy}
\usepackage{lscape}
\usepackage{pict2e}

\usepackage{graphicx}
\graphicspath{{./}{SigPicture/}}
\DeclareGraphicsExtensions{.ps}

\newtheorem{theo}{Theorem}[section]
\newtheorem{prop}[theo]{Proposition}
\newtheorem{pte}[theo]{Property}
\newtheorem{defi}[theo]{Definition}
\newtheorem{lem}[theo]{Lemma}
\newtheorem{cor}[theo]{Corollary}

\newcommand{\eq}[1][r]
   {\ar@<-3pt>@{-}[#1]
    \ar@<-1pt>@{}[#1]|<{}="gauche"
    \ar@<+0pt>@{}[#1]|-{}="milieu"
    \ar@<+1pt>@{}[#1]|>{}="droite"
    \ar@/^2pt/@{-}"gauche";"milieu"
    \ar@/_2pt/@{-}"milieu";"droite"}
\newcommand{\citeg}[1]{\textbf{\cite{#1}}}
\newcommand{\refg}[1]{\textbf{\ref{#1}}}
\newcommand{\reffig}[1]{\textbf{Figure \ref{#1}}}

\newcommand{\RR}{\mathbb R}
\newcommand{\CC}{\mathbb C}
\newcommand{\TT}{\mathcal T}
\newcommand{\Aa}{\mathcal A}
\newcommand{\EE}{\mathbb E}
\newcommand{\ZZ}{\mathbb Z}
\newcommand{\QQ}{\mathbb Q}
\newcommand{\NN}{\mathbb N}

\newcommand{\HH}{\mathcal{H}}

\newcommand{\E}[1]{\mathcal{E}(\Omega,#1;\CC \oplus \CC)}
\newcommand{\hti}{\hat{\otimes}}
\newcommand{\Qi}{\Xi_{00}}
\newcommand{\wo}{\omega_0}
\newcommand{\itemb}{\item[$\bullet$]}

\title{\textbf{\huge{An index theorem to solve the gap-labeling conjecture for the pinwheel tiling}}}
\author{Ha\"ija MOUSTAFA}
\date{}

\begin{document}
\maketitle

\vspace{-0.3cm}
\begin{abstract}
\noindent
In this paper, we study the $K_0$-group of the $C^*$-algebra associated to a pinwheel tiling. We prove that it is given by the sum of $\ZZ \oplus \ZZ^6$ with a cohomological group. The $C^*$-algebra is endowed with a trace that induces a linear map on its $K_0$-group. We then compute explicitly the image, under this map, of the summand $\ZZ \oplus \ZZ^6$, showing that the image of $\ZZ$ is zero and the image of $\ZZ^6$ is included in the module of patch frequencies of the pinwheel tiling (see \citeg{Hai-coho}). We finally prove that we can apply the measured index theorem due to A. Connes (\citeg{Con1}) to relate the image of the last summand of the $K_0$-group to a cohomological formula which is more computable. This is the first step in the proof of the gap-labeling conjecture for the pinwheel tiling, the second step is done in \citeg{Hai-coho} where we study the cohomological formula obtained by the index theorem.
\end{abstract}

\vspace{-0.1cm}
\tableofcontents

\newpage

\section{Introduction}

\begin{figure}[!ht]
\begin{center}
\includegraphics[scale=0.2]{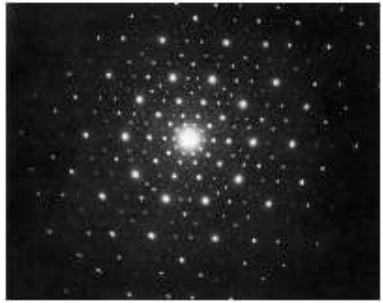}
\end{center}
\caption{Diffraction diagram.}
\label{diffraction}
\end{figure}

In 1982, D. Shechtman discovered, in a rapidly solidified aluminium alloy, a phase similar to the one obtained from crystals (\citeg{Sheetal}). To study the atomic distribution of such a solid, he realized its diffraction diagram that is shown in \reffig{diffraction}.\\

\noindent
This diagram is discrete and is similar to the one obtained from crystals. From its study, he could deduce that the atoms are distributed in such a way that the distance between two atoms is bounded below and such that there is only a finite number of local configurations upto translations.\\
However, the diffraction diagram of this solid has a 10-fold symmetry that is forbidden by crystallographic classification.\\
The structure of this solid is thus different from the one of crystals.\\
In particular, the atomic distribution isn't stable by any translation. Nevertheless, each local configuration is repeted uniformly in the alloy.\\
For these reasons, such a solid was called \textit{quasicrystal}.

\bigskip
This discovery gave rise to a big interest in solid physics and in mathematics.\\
In mathematics, this atomic distribution is naturally modeled by aperiodic tilings encoding its geometric properties into combinatoric properties.

\medskip
\noindent
A \textit{tiling} of $\RR^n$ is a countable family $P = \{t_0, t_1, \ldots \}$ of non empty compact sets $t_i$, called \textit{tiles}, each tile being homeomorphic to the unit ball, such that the $t_i$'s cover $\RR^n$, with a fixed origin, and the tiles meet each other only on their border.\\
In the sequel, we consider the particular case in which the $t_i$'s are $G$-copies of elements of a finite family $\{p_0, \ldots, p_n \}$, where $G$ is some subgroup of the group of rigid motions of $\RR^n$. The $p_k$'s are called \textit{prototiles}.\\
A \textit{patch} is a finite collection of tiles in a tiling.\\
A tiling is of \textbf{finite $\mathbf{G}$-type} if, for all $R>0$, there are only finitely many patches of diameter less than or equal to $R$ upto the $G$-action.\\
A tiling $\TT'$ is \textbf{$\mathbf{G}$-repetitive} if, for each patch $\mathcal{A}$ in $\TT'$, there exists $R(\mathcal{A})>0$ such that every ball of radius $R(\mathcal{A})$ in $\RR^n$ meets $\TT'$ on a patch containing a $G$-copy of $\mathcal{A}$.\\
In this paper, a \textbf{$\mathbf{G}$-aperiodic tiling} is a nonperiodic tiling (for translations of $\RR^n$) of finite $G$-type and $G$-repetitive.\\
The atoms of a quasicrystal are located in the tiles of the aperiodic tiling modeling it.

\bigskip
Motivated by electromagnetic and macroscopic properties of such solids, J. Bellissard studied the electonic motion in a quasicrystal (see \citeg{Bel82}, \citeg{Bel86} and \citeg{Bel1}).\\
This motion is closely related to the spectral gaps of the Shr\"odinger operator associated to the solid. This operator is defined as follows :
$$H = -\dfrac{\hslash}{2m} \Delta + V(x)$$
where $-\frac{\hslash}{2m}$ is some constant, $\Delta$ is the Laplacian on $\RR^n$ and $V$ is some potential depending only on the structure of the studied quasicrystal.\\
$V$ is given by $V(x)={\displaystyle \sum_{y \in L} v(x-y) }$ where $L$ is the point set of equilibrium atomic positions in the quasicrystal and $v$ the effective potential for valence electrons near an atom. In other words, $v$ is the function governing the interaction between an electron and an atom (see \citeg{BelHerZar}).\\
Bellissard looked for a mathematical way to label the gaps of the spectrum of this Schr\"odinger operator and this is the aim of its gap-labeling conjecture.

\bigskip
Tilings were studied before the discovery of quasicrystals but the interest generated by this new material greatly increased progresses made in their study, looking tiling properties from a new point of view, the noncommutative one.\\
In 2000, in \citeg{KelPut}, J. Kellendonk and I. Putnam associated a $C^*$-algebra to a $G$-aperiodic tiling to study its properties. It is the crossed product of the continuous functions on a topological space $\Omega$ by $G$. The space $\Omega$ encodes the combinatorial properties  of the tiling into  topological and dynamical properties.

\bigskip
We recall, in section $\mathbf{2}$, how this space is obtained together with other definitions and the construction of the $C^*$-algebra associated to tilings.\\
To construct $\Omega$, let's fix an orthonormal basis of $\RR^n$ and a $G$-aperiodic tiling $\TT$.\\
Letting $G$ act on the right on $\TT$, we have a family of tilings $\TT.G$ and one can take its closure under some distance (adapted to such tilings) to obtain the \textbf{continuous hull} $\Omega$ of $\TT$ endowed with a natural action of $G$ and containing tilings with the same properties as $\TT$ (this construction is due to Kellendonk \citeg{Kelgap}) .\\
The combinatorial properties of $\TT$ ensure that $\Omega$ is a compact topological space and that each orbit under the action of $G$ is dense in $\Omega$.

\bigskip
The \textbf{canonical transversal} $\Xi$ is the set of all the tilings of $\Omega$ satisfying some rigid conditions.\\
This is a Cantor space and allows us to see $\Omega$ as a foliated space, i.e $\Omega$ is covered by open sets given by $K \times U$, where $K$ is a clopen (closed and open) subset of $\Xi$ and $U$ is an open subset of $G$, with some smooth conditions on transition maps. The space $\Xi$ is then a transversal of $\Omega$ and the leaves of $\Omega$ are homeomorphic to $G$ or to a quotient of $G$ by a discrete subgroup.

\bigskip
We can associate to the dynamical system $(\Omega,G)$ the $C^*$-crossed product $C(\Omega) \rtimes G$ which is defined in the end of section $2$.\\
Since $G$ is amenable, $\Omega$ is endowed with a $G$-invariant, ergodic probability measure $\mu$. This measure induces a densely defined trace $\tau^\mu$ on $C(\Omega) \rtimes G$ and a linear map $\tau^\mu_* : K_0\big( C(\Omega) \rtimes G \big) \rightarrow \RR$ on the $K$-theory of $C(\Omega) \rtimes G$.

\bigskip
In section $\mathbf{3}$, we present the works made by J. Bellissard.\\
For $G=\RR^n$, J. Bellissard linked the spectral gaps of the Schr\"odinger operator $H$ to this $K$-theory, noting that each spectral gap gives rise to an element of $K_0 \big( C(\Omega) \rtimes \RR^n \big)$ (see \citeg{Bel1} and \citeg{BelHerZar}).\\
He also showed that the image under $\tau^\mu_*$ of $K_0 \big( C(\Omega) \rtimes \RR^n \big)$ is related to some physical function $\mathcal{N} : \RR \rightarrow \RR^{+}$, called the \textit{integrated density of states}, abreviated IDS (see p.\pageref{IDS}).\\
For $E \in \RR$, $\mathcal{N}(E)$ can be seen as the number of eigenvalues of $H$ per unit volume up to $E$.\\
Thus, this is a nonnegative, nondecreasing function on $\RR$ and thereby, it defines a Stieljes-Lebesgue measure $d\mathcal{N}$ on $\RR$ given by $d\mathcal{N} \big( [E ; E'] \big) := \mathcal{N}(E') - \mathcal{N}(E)$.\\
$d\mathcal{N}$ is absolutely continuous with respect to the Lebesgue measure $dE$ (i.e $\int_A dE= 0 \Rightarrow \int_A d\mathcal{N}=0$) and, by the Radon-Nikodym theorem, we can define the \textit{density of states} to be the derivative $d\mathcal{N}/dE$.\\
This density of states is a well known quantity in solid state physics which is in fact accessible by scattering experiments (see \citeg{Bel1}).\\
This quantity is important in physics to deduce conductivity properties of quasicrystals and it is crucial to predict its values.\\
For this, J. Bellissard proved that this function is related to a mathematical object thanks to \textbf{Shubin's formula}. Shubin's formula expresses that, if $E$ is in a spectral gap, there is a projection $p_E$ in $C(\Omega) \rtimes \RR^n$ such that $\mathcal{N}(E)$ is equal to $\tau^\mu(p_E)$ and that this does not depend on $E$ chosen in the spectral gap.\\
Thereby, the integrated density of states $\mathcal{N}$ takes values in the image, under $\tau^\mu_*$, of $K_0 \big( C(\Omega) \rtimes \RR^n \big)$.\\
To study conductivity properties of quasicrystals, it is thus important to know this image.

\bigskip
The gap-labeling conjecture established by Bellissard predicts this image : $\mu$ induces a measure $\mu^t$ on the canonical transversal $\Xi$ (given locally by the quotient of $\mu(K \times U)$ by the volume of $U$) and the conjecture expresses the link between  the image of $K_0 \big( C(\Omega) \rtimes \RR^n \big)$ under $\tau^\mu_*$ and the image, under $\mu^t$, of continuous functions on $\Xi$ with integer values:

\bigskip
\noindent
\textbf{Conjecture : }  (\citeg{Bel1}, \citeg{BelHerZar}) \textit{
$$\tau_*^\mu \Big( K_0 \big( C(\Omega) \rtimes \RR^n \big) \Big) = \mu^t \big( C(\Xi,\ZZ) \big)$$}

\vspace{-0.3cm}
\noindent
In this conjecture, $\mu^t \big( C(\Xi,\ZZ) \big)$ is the space ${\displaystyle \left \{ \int_\Xi f d\mu^t ; f \in C(\Xi,\ZZ) \right \}}$.

\bigskip
Since then, several works have been done to prove this conjecture.\\
The Pimsner-Voiculescu exact sequence was used by J. Bellissard in \citeg{Bel1} to prove the conjecture for $n=1$ and A. van Elst iterated this sequence twice in \citeg{vanElst} to prove the case $n=2$.\\
J. Bellissard, J. Kellendonk and A. Legrand have proved the conjecture for $n=3$ in \citeg{BelKelLeg} using a spectral sequence due to Pimsner.\\
Finally, in 2002, the conjecture was proved in full generality for $\RR^n$-aperiodic tilings in three independent papers : by J. Bellissard, R. Benedetti and J.-M. Gambaudo in \citeg{BelBenGam}, by M.-T. Benameur and H. Oyono-Oyono in \citeg{BenOyo} and by J. Kaminker and  I. Putnam in \citeg{KamPut}.

\bigskip
A natural question is to know if this conjecture remains true for more general tilings.\\
This paper gives a first step to solve the gap-labeling conjecture in the particular case of the \textit{$\mathit{(1,2)}$-pinwheel} tiling of the plane introduced by Conway and Radin (see \citeg{Rad} and \citeg{Rad1}) which is not a $\RR^2$-aperiodic tiling but a $\RR^2 \rtimes S^1$-aperiodic one.\\
It is a tiling built from two right triangles of side $1$, $2$ and $\sqrt{5}$, one being the mirror image of the other.\\
The \textit{substitution} method gives a pinwheel tiling as follows : we can cover the stretched image, by a $\sqrt{5}$ factor, of the two triangles by the union of copies by a rigid motion of these triangles, these copies meeting only on their border, as follows :
\begin{figure}[!ht]
\begin{center}
\includegraphics[scale=0.25]{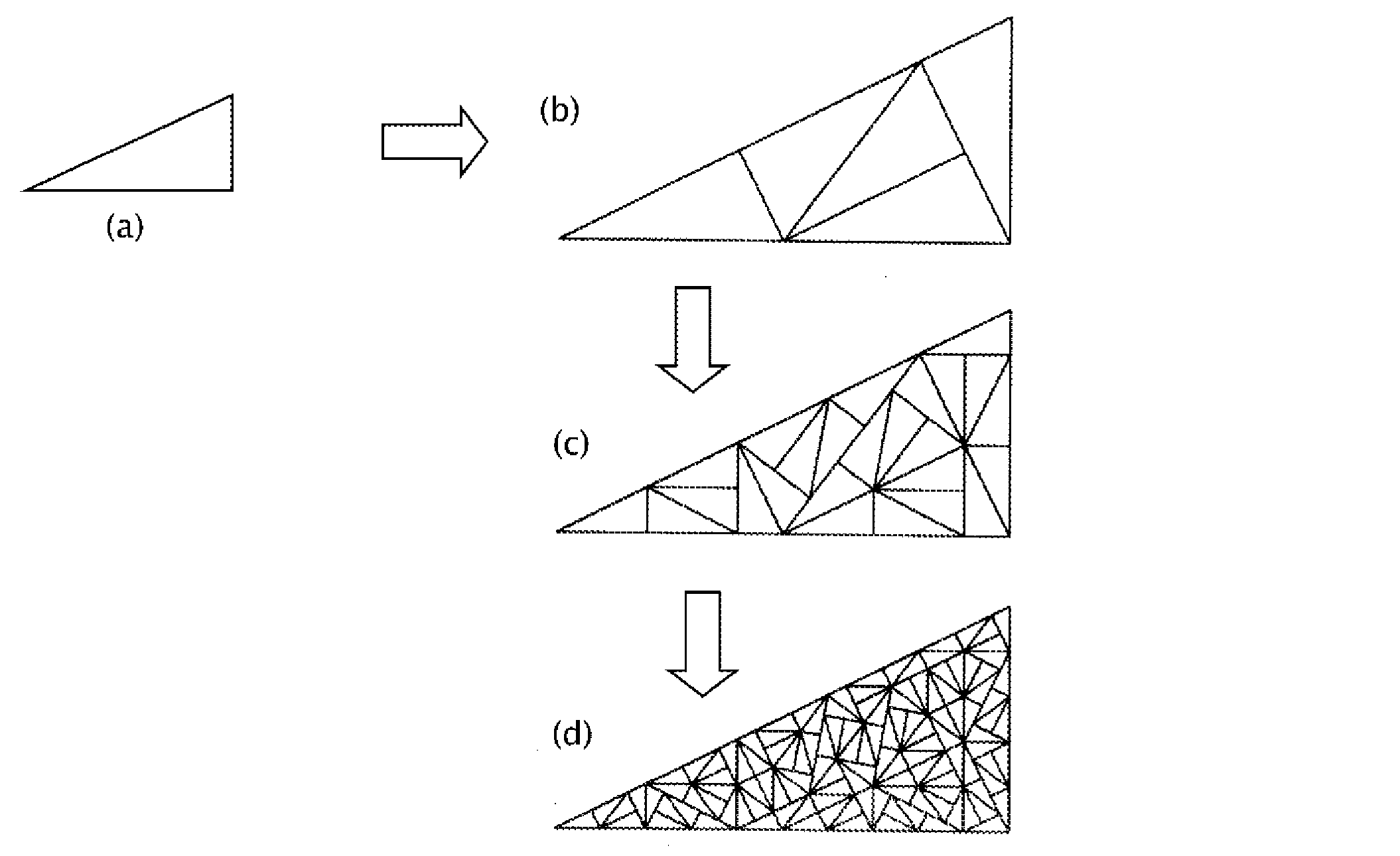}
\end{center}
\caption{Substitution of the pinwheel tiling.}
\label{substintro}
\end{figure}

\begin{figure}[!ht]
\begin{center}
\includegraphics[width=\textwidth]{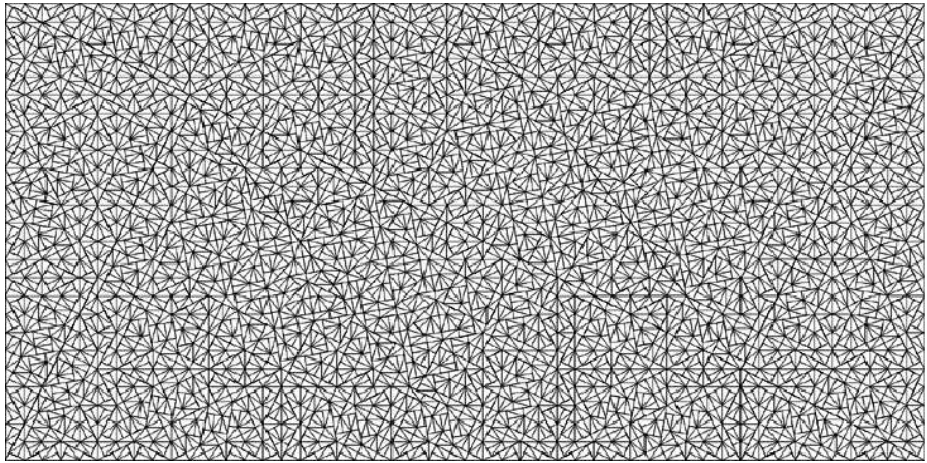}
\end{center}
\caption{Patch of a pinwheel tiling.}
\label{pinwhintro}
\end{figure}

\bigskip
Iterating this process, we build a union of triangles that covers bigger and bigger regions of the plane and that, in the limit, covers the whole plane and, thus, gives a tiling of the plane called a \textbf{pinwheel tiling}. A patch of this tiling is shown in \reffig{pinwhintro}.

\noindent
This tiling is nonperiodic for translations of the plane, of finite $\RR^2 \rtimes S^1$-type and $\RR^2 \rtimes S^1$-repetitive, where $\RR^2 \rtimes S^1$ is the group of rigid motion of the plane.\\
It is not of finite $\RR^2$-type since the triangles appear in infinitly many orientations.\\
The continuous hull $\Omega$ of this tiling contains six circles fixed by a rotation of angle $\pi$ around the origin. The union of these $6$ circles is denoted $F$ and these are the only fixed points for rotations.

\bigskip
The constrution explained above can be applied to obtain a dynamical system $(\Omega, \RR^2 \rtimes S^1)$, where $\Omega$ is a compact topological space and such that each $\RR^2 \rtimes S^1$-orbit is dense in $\Omega$.\\
The canonical transversal $\Xi$ can still be built and $\Omega$ is a foliated space with leaves homeomorphic to $\RR^2 \rtimes S^1$ except six that are homeomorphic to the quotient of $\RR^2 \rtimes S^1$ by the subgroup generated by $(0,R_\pi)$.

\bigskip
\noindent
In this paper, we give a first step of the proof of the following theorem, the second step is given in \citeg{Hai-coho} :

\medskip
\noindent
\textbf{Theorem :} \textit{If $\TT$ is a pinwheel tiling, $\Omega = \Omega(\TT)$ its hull provided with an invariant ergodic probability measure $\mu$ and $\Xi$ its canonical transversal provided with the induced measure $\mu^t$, then}
$$\tau^\mu_* \Big( K_0 \big( C(\Omega \times \RR^2 \rtimes S^1) \big) \Big) = \mu^t \big( C(\Xi,\ZZ) \big).$$

\noindent
This theorem aims to link the image of $ K_0 \big( C(\Omega \times \RR^2 \rtimes S^1) \big)$ under $\tau^\mu_*$ to the $\ZZ$-module of patch frequencies $\mu^t \big( C(\Xi,\ZZ) \big)$.

\bigskip
For this, we follow the guideline of the works made by M.-T. Benameur and H. Oyono-Oyono in \citeg{BenOyo}.\\
In their paper, they used the measured index theorem proved by A. Connes (see \citeg{Con1} or \citeg{MooSch}) to link the analytical part $\tau^\mu_* \Big( K_0 \big( C(\Omega \times \RR^n) \big) \Big)$ to a topological part, easier to compute,  $Ch_\tau \big( K_n(C(\Omega)) \big)$ which is in $H^*_\tau \big( \Omega \big)$, the tangential cohomology group of $\Omega$ ($Ch_\tau$ is the tangential Chern character, see \citeg{MooSch}).

\bigskip
\noindent
In section $\mathbf{4}$, we prove that this theorem can still be applied for pinwheel tilings.\\
For this, we want to prove that any element of $K_0 \big( C(\Omega) \rtimes \RR^2 \rtimes S^1 \big)$ can be seen as the analytical index of a Dirac operator "twisted by a unitary" of $K_1 \big( C(\Omega) \big)$, i.e for any $b \in K_0 \big( C(\Omega) \rtimes \RR^2 \rtimes S^1 \big)$, we want to find a $[u] \in K_1 \big( C(\Omega) \big)$ such that $b=[u] \otimes_{C(\Omega)} [D_3]$, where $[D_3] \in KK_1 \big( C(\Omega), C(\Omega) \rtimes \RR^2 \rtimes S^1 \big)$ is the class of the tangential Dirac operator $D_3$ along the leaves of $\Omega$ and $\otimes_{C(\Omega)}$ is the Kasparov product over $C(\Omega)$.\\
In fact, we prove that $K_0 \big( C(\Omega) \rtimes \RR^2 \rtimes S^1 \big)$ is isomorphic to the sum of $\ZZ^7$ and a subgroup $H$ (see below for more details on $H$).\\
The image under $\tau_*^\mu$ of $\ZZ^7$ is explicitly computable and for any $b \in H$, there is some $[u] \in K_1 \big( C(\Omega) \big)$ such that $\tau^\mu_*(b)=\tau^\mu_* \big( [u] \otimes_{C(\Omega)} [D_3] \big)$, which is enough to apply the index theorem to obtain :

\bigskip 
\noindent
\textbf{Theorem \refg{indice} :}  \textit{
$\forall b \in H$, $\exists [u] \in K_1 \big( C(\Omega) \big)$ such that
$$\tau^\mu_*(b)=\tau^\mu_*([u] \otimes_{C(\Omega)} [D_3]) = \big \langle Ch_\tau ( [u]) , [C_{\mu^t}] \big \rangle,$$
where $[C_{\mu^t}] \in H^\tau_3(\Omega)$ is the Ruelle-Sullivan current associated to $\mu^t$ (locally, it is given by integration on $\Xi \times \RR^2 \times S^1$, see \citeg{MooSch}) and $\langle \, , \, \rangle$ is the pairing between the tangential cohomology and homology groups.}

\bigskip
To obtain this theorem, we see $\Omega$ as a double foliated space.\\
First, it is foliated by the action of $S^1 \subset \RR^2 \rtimes S^1$ and we can build the map $\psi_1 : K_1\big( C(\Omega) \big) \rightarrow K_0 \big( C(\Omega) \rtimes S^1 \big)$ given by the Kasparov product, over $C(\Omega)$, with the class of the cycle induced by the tangential Dirac operator $d_1$ along the leaves $S^1$.\\
Then, $\Omega$ is also foliated by the action of $\RR^2 \rtimes S^1$ and thus, the Kasparov cycle defined by the tangential Dirac operator $D_2$ transverse to the inclusion of the foliation by $S^1$ into the foliation by $\RR^2 \rtimes S^1$ is used to construct the map $\psi_2 : K_0 \big( C(\Omega) \rtimes S^1 \big) \rightarrow K_0 \big( C(\Omega) \rtimes \RR^2 \rtimes S^1 \big)$ given as the Kasparov product, over $C(\Omega) \rtimes S^1$, with the class of this cycle in $KK_0 \big( C(\Omega) \rtimes S^1 , C(\Omega) \rtimes \RR^2 \rtimes S^1 \big)$.\\
Using the Dirac-dual Dirac construction (see \citeg{Kaspa2}, \citeg{SkanKK}), one can show that $\psi_2$ is an isomorphism.\\
We then prove that $\psi_2 \circ \psi_1$ is given by the Kasparov product with the class of the cycle induced by the tangential Dirac operator of dimension $3$ (see \citeg{HilSkan} where, in the case of foliations, they studied the product of Gysin maps associated to double foliations).

\bigskip
To prove these results, we first use a six term exact sequence in $K$-theory to show ($\check{H}^2(\cdot ; \ZZ)$ is the \v{C}ech cohomology group with integer coefficients and $\check{H}^2_c(\cdot ; \ZZ)$ is the one with compact support) :

\bigskip
\noindent
\textbf{Proposition \refg{somdirecte} :} 
$$K_0(C(\Omega) \rtimes S^1) \simeq \ZZ\oplus \ZZ^6 \oplus \, \check{H}^2_c \Big( (\Omega \setminus F)/S^1 \, ; \, \ZZ \Big)$$

\noindent
Then, we prove that $\tau^\mu_* \big( \psi_2(\ZZ) \big)=0$ and to study the image, under $\tau_*^\mu$, of $H:=\psi_2 \Big(  \check{H}^2_c \big( (\Omega \setminus F)/S^1 \, ; \, \ZZ \big) \Big )$, we prove that $\check{H}^2_c \Big( (\Omega \setminus F)/S^1 \, ; \, \ZZ \Big)$ is isomorphic to $\check{H}^3 \big( \Omega \, ; \, \ZZ \big)$ and that this isomorphism is $\psi_1$, modulo the Chern character.\\
Thus, $\psi_2 \circ \psi_1$ is surjective on the $H$ summand of $K_0 \big( C(\Omega) \rtimes \RR^2 \rtimes S^1 \big)$.\\
Using a result proved by Douglas, Hurder and Kaminker in \citeg{DouHurKam} on the odd index theorem for foliated spaces, we obtain the main theorem of this paper : 

\bigskip
\noindent
\textbf{Theorem \refg{indice} :}  \textit{
$\forall b \in H$, $\exists [u] \in K_1 \big( C(\Omega) \big)$ such that :
$$\tau^\mu_*(b)=\tau^\mu_*([u] \otimes_{C(\Omega)} [D_3]) = \big \langle Ch_\tau ( [u]) , [C_{\mu^t}] \big \rangle. $$
}

\bigskip
\noindent
Finally, an explicit computation, using the index theorem for foliated spaces on $\Omega$ seen as a foliated space for the $\RR^2$-action, gives the inclusion
$$\tau^\mu_* \big( \psi_2(\ZZ^6) \big) \subset \mu^t(C(\Xi,\ZZ)).$$

In this paper, we thus obtain the following result :

\bigskip
\noindent
\textbf{Theorem :}  \textit{If $\TT$ is a pinwheel tiling, $\Omega = \Omega(\TT)$ its hull provided with an invariant ergodic probability measure $\mu$ and $\Xi$ its canonical transversal provided with the induced measure $\mu^t$, then
$$K_0(C(\Omega) \rtimes \RR^2 \rtimes S^1) \overset{\psi_2}{\simeq} \ZZ\oplus \ZZ^6 \oplus \, \check{H}^2_c \Big( (\Omega \setminus F)/S^1 \, ; \, \ZZ \Big).$$
And
\begin{itemize}
	\itemb $\tau^\mu_* \big( \psi_2(\ZZ) \big)=0$.
	\itemb $\tau^\mu_* \big( \psi_2(\ZZ^6) \big) \subset \mu^t(C(\Xi,\ZZ))$.
	\itemb $\forall b \in H$, $\exists [u] \in K_1 \big( C(\Omega) \big)$ such that :
$$\tau^\mu_*(b)=\tau^\mu_*([u] \otimes_{C(\Omega)} [D_3]) = \big \langle Ch_\tau ( [u]) , [C_{\mu^t}] \big \rangle. $$
\end{itemize}
}

\medskip
\noindent
To prove the gap-labeling conjecture for pinwheel tilings, it thus remains to study the topological part of the last point of the theorem to prove that $$\tau_*^\mu \left( K_0 \big( C(\Omega) \rtimes \RR^2 \rtimes S^1 \big) \right) \subset \mu^t \big( C(\Xi,\ZZ) \big).$$
This is done in \citeg{Hai-coho} and the inclusion in the other direction is easily obtained in this paper too.

\bigskip
\noindent
\textbf{Aknowledgements. } It is a pleasure for me to thank my advisor Herv\'e Oyono-Oyono who always supported and advised me during this work.\\
I also want to thank Jean Bellissard for useful discussions on the gap-labeling conjecture.

\newpage

\section{Reminders}

\subsection{Pinwheel tiling and continuous hull}

\vspace{0.5cm}
A \textbf{tiling of the plane} is a countable family $P=\{t_1, t_2 \ldots \}$ of non empty compact subsets $t_i$ of $\RR^2$, called \textbf{tiles} (each tile being homeomorphic to the unit ball), such that:
\begin{itemize}
\item[$\bullet$] ${\displaystyle \bigcup_{i \in \mathbb{N}} t_i = E_2}$ where $E_2$ is the euclidean plane with a fixed origin $O$;
\item[$\bullet$] Tiles meet each other only on their border ;
\item[$\bullet$] Tiles's interiors are pairwise disjoint.
\end{itemize} 

\bigskip
\noindent
We are interested in the special case where there exists a finite family of tiles $\{p_1, \ldots, p_n \}$, called \textbf{prototiles}, such that each tile $t_i$ is the image of one of these prototiles under a \textit{rigid motion} (i.e. a direct isometry of the plane).\\
In fact this paper will focus on the particular tiling called \textbf{pinwheel tiling} or \textbf{(1,2)-pinwheel tiling} which is obtained by a substitution explained below.

\begin{figure}[ht]
\begin{center}
\includegraphics[scale=0.25]{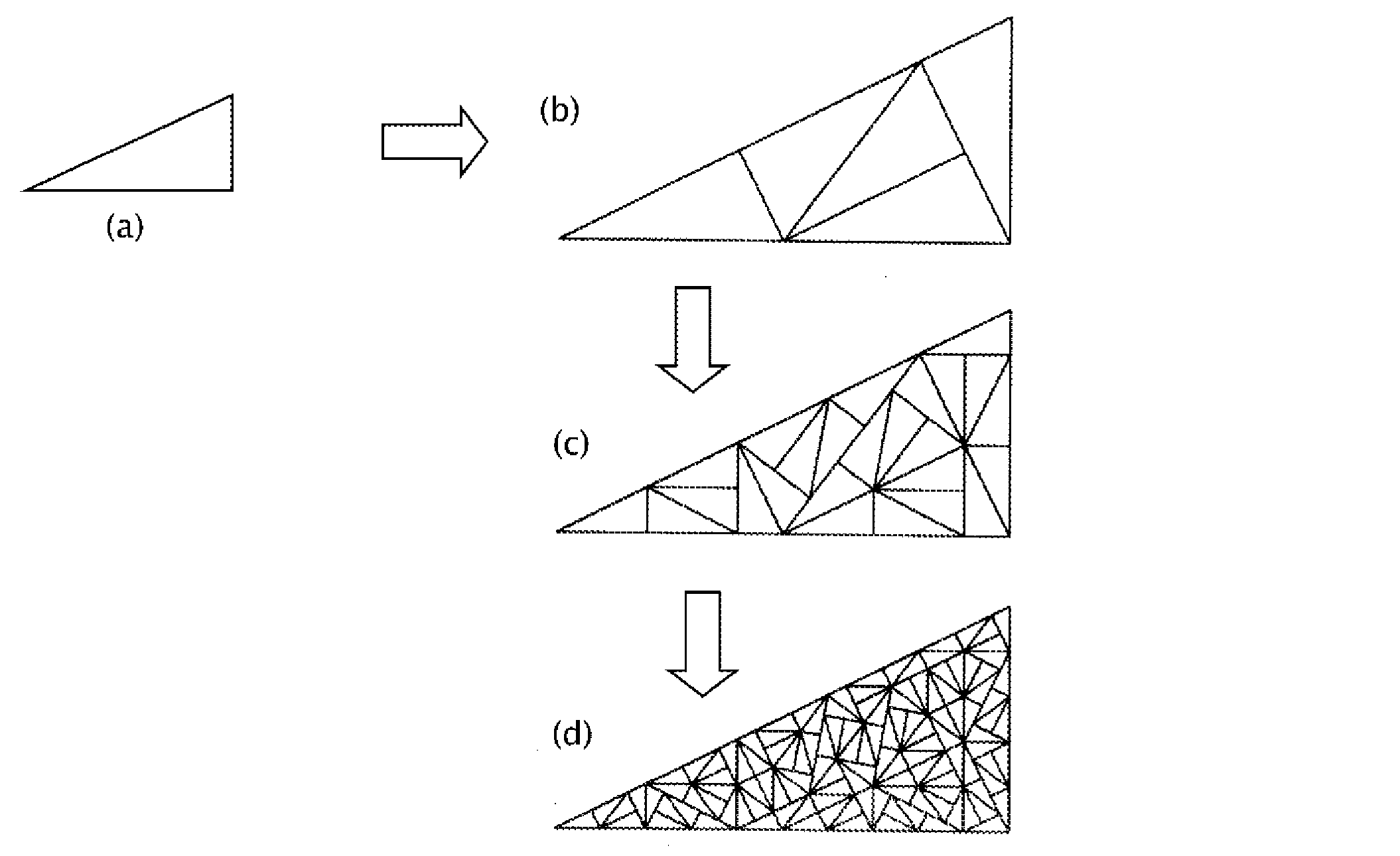}
\end{center}
\caption{Substitution of the pinwheel tiling.}
\label{Pin}
\end{figure}

\bigskip
\noindent
Our construction of a pinwheel tiling is based on the construction made by Charles Radin in \citeg{Rad}. It's a tiling of the plane obtained by the substitution described in \reffig{Pin}.

\bigskip
\noindent
This tiling is constructed from two prototiles, the right triangle in \reffig{Pin}\textbf{.(a)} with legs $1$, $2$ and $\sqrt{5}$ and its mirror image.

\begin{figure}[ht]
\begin{center}
\includegraphics[width=\textwidth]{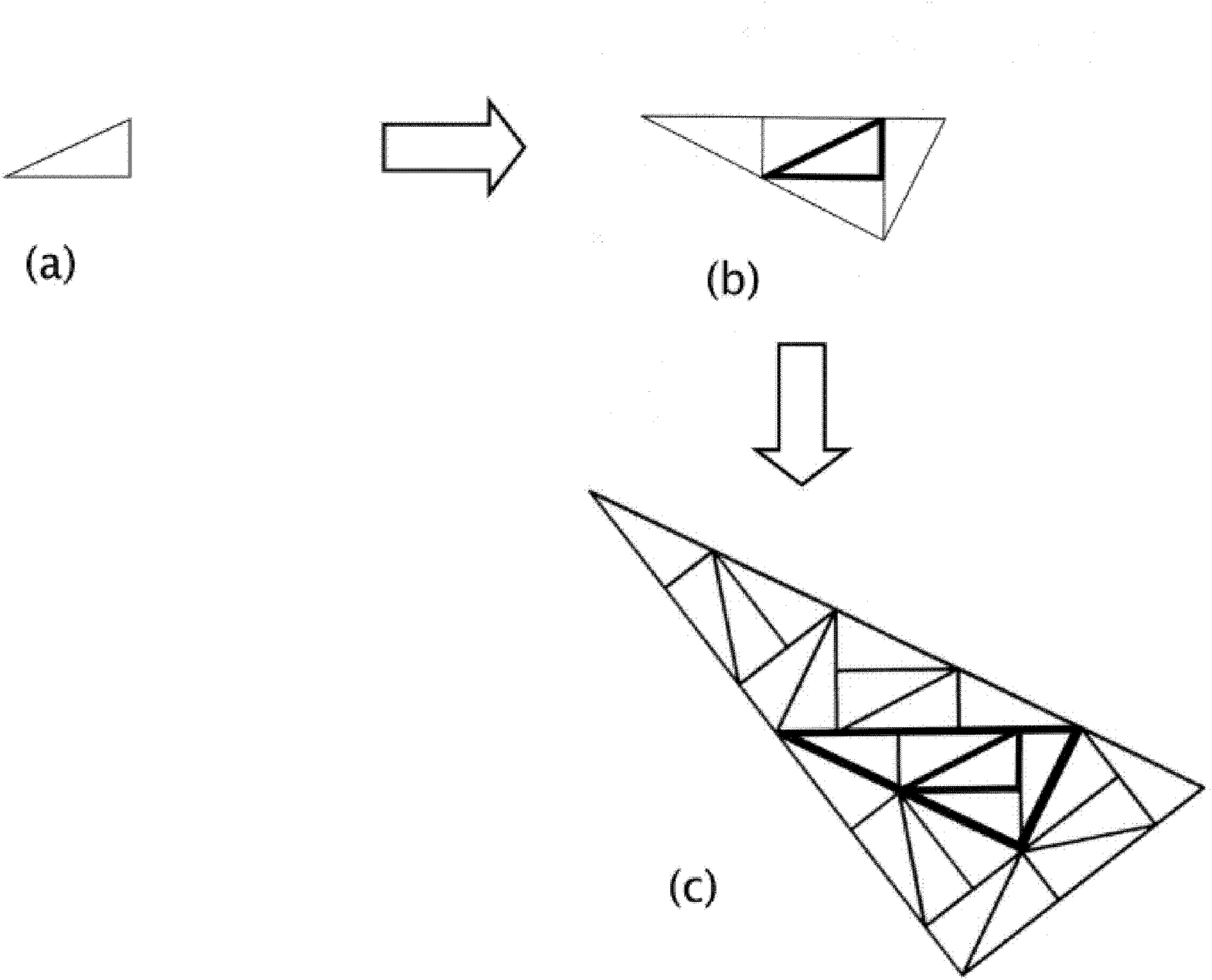}
\end{center}
\caption{Construction of a pinwheel tiling.}
\label{pinwrot}
\end{figure}

\bigskip
\noindent
To obtain this tiling, we begin from the right triangle with the following vertices in the plane : $(0,0)$ , $(2,0)$ and $(2,1)$.\\
This tile and its reflection are called \textbf{supertiles of level 0} or \textbf{0-supertiles}.\\
We will next define \textbf{1-supertiles} as follows : take the right triangle with vertices $(-2,1)$, $(2,-1)$ and $(3,1)$ and take the decomposition of \reffig{Pin}\textbf{.(b)}. This $1$-supertile is thus decomposed in five $0$-supertiles, which are isometric copies of the first tile, with the beginning tile in its center (see \reffig{pinwrot}\textbf{.(b)}).

\bigskip
\noindent
We next repeat this process by mapping this $1$-supertile in a $2$-supertile with vertices $(-5,5)$, $(1,-3)$ and $(5,0)$ (see \reffig{pinwrot}\textbf{.(c)}).\\
Including this $2$-supertile in a $3$-supertile with correct orientation and so on, this process leads to the desired \textbf{pinwheel tiling} $\TT$.

\bigskip
\noindent
We will now attach to this tiling a topological space reflecting the combinatorial properties of the tiling into topological and dynamical properties of this space.

\bigskip
\noindent
For this, we observe that the direct isometries of the plane are acting on the euclidean plane $E_2$ where we have fixed the origin $O$.\\
Direct isometries $\EE^2=\RR^2 \rtimes S^1$ thus act naturally on our tiling $\TT$ on the right .\\
If $R_\theta$ denotes the rotation about the origin with angle $\theta$ and $s \in \RR^2$, $\TT.(s,R_\theta):=R_{-\theta}(\TT-s)$. 
We will also denote $(s,R_\theta)$ by $(s,\theta)$.
\begin{defi}
A \textbf{patch} is a finite union of tiles of a tiling.\\
A tiling $\TT'$ is of \textbf{finite $\mathbf{\EE^2}$-type} or of \textbf{Finite Local Complexity} (\textbf{FLC}) if for any $R > 0$, there is only a finite number of patches in $\TT'$ of diameter less than $R$ up to direct isometries.\\
A tiling $\TT'$ of finite $\mathbf{\EE^2}$-type is \textbf{$\mathbf{\EE^2}$-repetitive} if, for any patch $\Aa$ in $\TT'$, there is $R(\Aa)>0$ such that any ball of radius $R(\Aa)$ intersects $\TT'$ in a patch containing a $\EE^2$-copy of $\Aa$.
\end{defi}

\noindent
The tiling $\TT$ is of finite $\EE^2$-type, $\mathbf{\EE^2}$-repetitive and nonperiodic for translations (see \citeg{Pet}).

\bigskip
\noindent
To attach a topological space to $\TT$, we define a metric on $\TT.\EE^2$ :\\
If $\TT_1$ and $\TT_2$ are two tilings in $\TT.\EE^2$, we define 
$$\hspace*{-2cm} A=\left \{\varepsilon \in \Big [ 0,\tfrac{1}{\sqrt{2}} \Big ] \; / \; \exists s,s' \in B^2_\varepsilon(0) \, , \, \theta,\theta' \in B^1_\varepsilon(0) \text{ s.t. } \right.$$ 
$$\hspace*{3cm} \left. \TT_1.(s,\theta) \cap B_{\frac{1}{\varepsilon}}(O) =  \TT_2.(s',\theta') \cap B_{\frac{1}{\varepsilon}}(O) \right \}$$
where $B_{\frac{1}{\epsilon}}(O)$ is the euclidean ball centered in $O$ with radius $\frac{1}{\epsilon}$ and $B_\epsilon^i(0)$ are the euclidean balls in $\RR^i$ centered in $0$ and with radius $\epsilon$ (i.e. we consider direct isometries near $Id$).\\
Then, define :
$$d(\TT_1,\TT_2) = \left\{
\begin{array}{ll}
Inf A & \text{ if } A \neq \emptyset\\
\frac{1}{\sqrt{2}} & \text{ else } \quad \\
\end{array} \right..$$
$d$ is a bounded metric on $\TT.\EE^2$. For this topology, a base of neighborhoods is defined by: two tilings $\TT_1$ and $\TT_2$ are close if, up to a small direct isometry, they coincide on a large ball around the origin.\\
The topology defined here works because the tilings considered are of finite $\EE^2$-type. There exist other topologies (equivalent) that one can put on this space.\\
The topology thus obtained is metrizable but none of the metrics that can be defined to produce the topology is canonical.\\
A more canonical way to define the topology was given in \citeg{BelHerZar}.\\
The problem of non-uniqueness of the metric has been investigated in \citeg{PeaBel}.

\begin{defi} \text{ }\\
The \textbf{continuous hull} of $\TT$ is then the completion of $(\TT.\EE^2,d)$ and will be denoted $\Omega(\TT)$.
\end{defi}

\noindent
Let's enumerate some well known properties of this continuous hull:

\begin{pte}  (\citeg{BelBenGam}, \citeg{BelHerZar}, \citeg{BenGam}, \citeg{KelPut}, \citeg{LagPle}, \citeg{Rad}) \label{decompo}
\begin{itemize}
	\item $\Omega(\TT)$ is formed by finite $\EE^2$-type, $\mathbf{\EE^2}$-repetitive and nonperiodic (for translations) tilings and each tiling of $\Omega(\TT)$ has the same patches as $\TT$.
	\item $\Omega(\TT)$ is a compact space since $\TT$ is of finite $\EE^2$-type.
	\item Each tiling in $\Omega(\TT)$ are uniquely tiled by $n$-supertiles, for all $n \in \NN$.
	\item The dynamical system $(\Omega(\TT), \EE^2)$ is minimal since $\TT$ is repetitive, i.e each orbit under direct isometries is dense in $\Omega(\TT)$.
\end{itemize}
\end{pte}

\noindent
The last property of $\Omega(\TT)$ allows us to write $\Omega$ without mentioning the tiling $\TT$ (in fact, if $\TT' \in \Omega(\TT)$, $\Omega(\TT') = \Omega(\TT)$).

\begin{defi}
Any tiling in $\Omega$ is called a \textbf{pinwheel tiling}.
\end{defi}

\noindent
\textbf{Remark :} we can easily see that our continuous hull is the compact space $X_\phi$ defined by Radin and Sadun in \citeg{RadSad}.

\subsection{The canonical transversal}

\vspace{0.5cm}
In this section, we will construct a Cantor transversal for the action of $\EE^2$ and we show that this transversal gives the local structure of a foliated space.

\noindent
For this, we fix a point in the interior of the two prototiles of the pinwheel tiling. This, in fact, gives for any tiling $\TT_1$ in $\Omega$ (i.e constructed by these two prototiles), a punctuation of the plane denoted $\TT_1^{punct}$.\\
Define then $\Omega_0$ to be the set of every tilings $\TT_1$ of $\Omega$ such that $O \in \TT_1^{punct}$.\\
The \textbf{canonical transversal} is the space $\Omega_0 / S^1$.

\bigskip
\noindent
We can identify this space with a subspace of $\Omega$ by constructing a continuous section $s:\Omega_0/S^1 \longrightarrow \Omega$.\\
To obtain such a section, we fix an orientation of the two prototiles of our tilings once for all. Hence when we consider a patch of a tiling in the transversal $\Omega_0$, there is only one orientation of this patch where the tile containing the origin have the orientation chosen for the prototiles.\\
Let's then $[\omega]$ be in $\Omega_0 / S^1$, there is only one $\theta \in [0;2\pi[$ such that the tile in $R_\theta(\omega)$ containing the origin has the good orientation.\\
We define $s([\omega]) := R_\theta(\omega)$.\\
$s$ is well defined because $\theta$ depends on the representative $\omega$ chosen but not $R_\theta(\omega)$.\\
$s:\Omega_0 / S^1 \longrightarrow s(\Omega_0 / S^1)$ is then a bijection. We easily see that $s$ is continuous and thus it is a homeomorphism from the canonical transversal onto a compact subspace $\Xi$ of $\Omega$.\\
We also call this space the \textbf{canonical transversal}.\\
We can see $\Xi$ as the set of all the tilings $\TT_1$ in $\Omega$ with the origin on the punctuation of $\TT_1$ and with the tile containing the origin in the orientation chosen for the prototiles.

\bigskip
\noindent
We then have :

\begin{prop} (\citeg{BenGam})\\
The canonical transversal is a Cantor space.
\end{prop}
 
\noindent
A base of neighborhoods is obtained as follows : consider $\TT' \in \Xi$ and $\Aa$ a patch around the origin in $\TT'$ then 
$$U(\TT',\mathcal{A}) = \{\TT_1 \in \Xi \; \vert \; \TT_1 = \TT'  \text{ on } \Aa \}$$ \label{clopen}
is a closed and open set in $\Xi$, called a \textbf{clopen} set.

\bigskip
\noindent
Before defining the foliated stucture on $\Omega$, we must study the rotations which can fix tilings in $\Omega$.\\
In pinwheel tilings, we can sometimes find regions tiled by supertiles of any level and so we introduce the following definition:

\begin{defi}
A region of a tiling which is tiled by $n$-supertiles for all $n \in \NN$ is called an \textbf{infinite supertile} or \textbf{supertile of infinite level}.
\end{defi}

\noindent
If a ball in a tiling $\TT_1$ fails to lie in any supertile of any level n, then $\TT_1$ is tiled by two or more supertiles of infinite level, with the offending ball straddling a boundary.\\
We can, in fact, construct a pinwheel tiling with two half-planes as infinite supertiles as follows: Consider the rectangle consisting of two $(n-1)$-supertiles  in the middle of a $n$-supertile. For each $n \geqslant 1$, orient this rectangle with its center at the origin and its diagonal on the $x$-axis, and fill out the rest of a (non-pinwheel) tiling $\TT_n$ by periodic extension. By compactness this sequence has a convergent subsequence, which will be a pinwheel tiling and which will consist of two infinite supertiles (this example comes from \citeg{RadSad}).

\bigskip
\noindent
Note that the boundary of an infinite supertile must be either a line, or have a single vertex, since it is tiled by supertiles of all levels.\\
We call such a line a \textbf{fault line}.

\begin{lem}
If $(s,\theta)$ fixes a pinwheel tiling $\TT'$ then $\theta \in \{ 0,\pi \} mod(2\pi)$.\\
Moreover, if $\theta = 0$ then $s=0$. In other terms, translations can't fix a pinwheel tiling.
\end{lem}

\begin{landscape}
\begin{figure}[ht] 
\begin{center}
\includegraphics[scale=0.24]{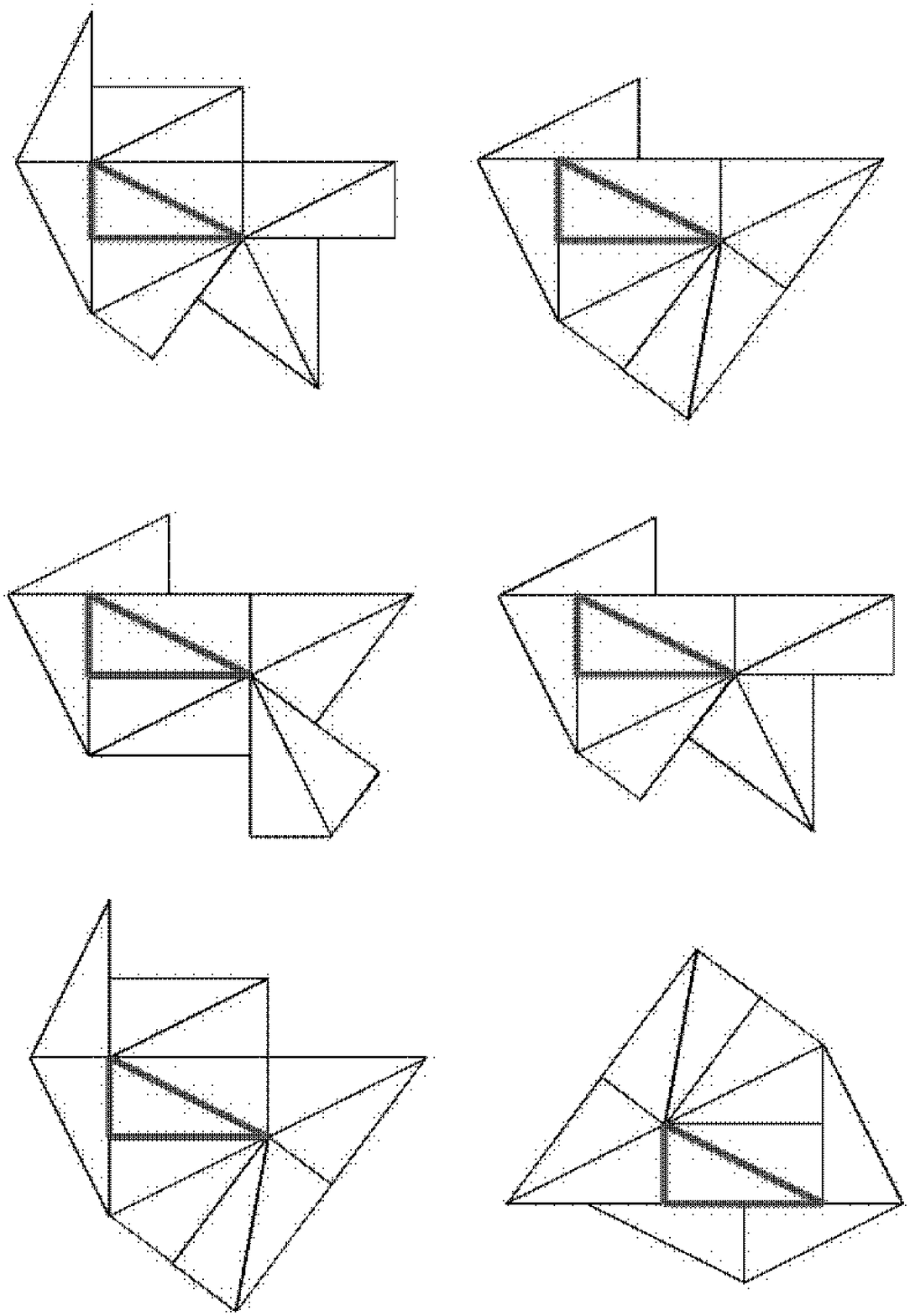}
\hspace{0.2cm}
\includegraphics[scale=0.24]{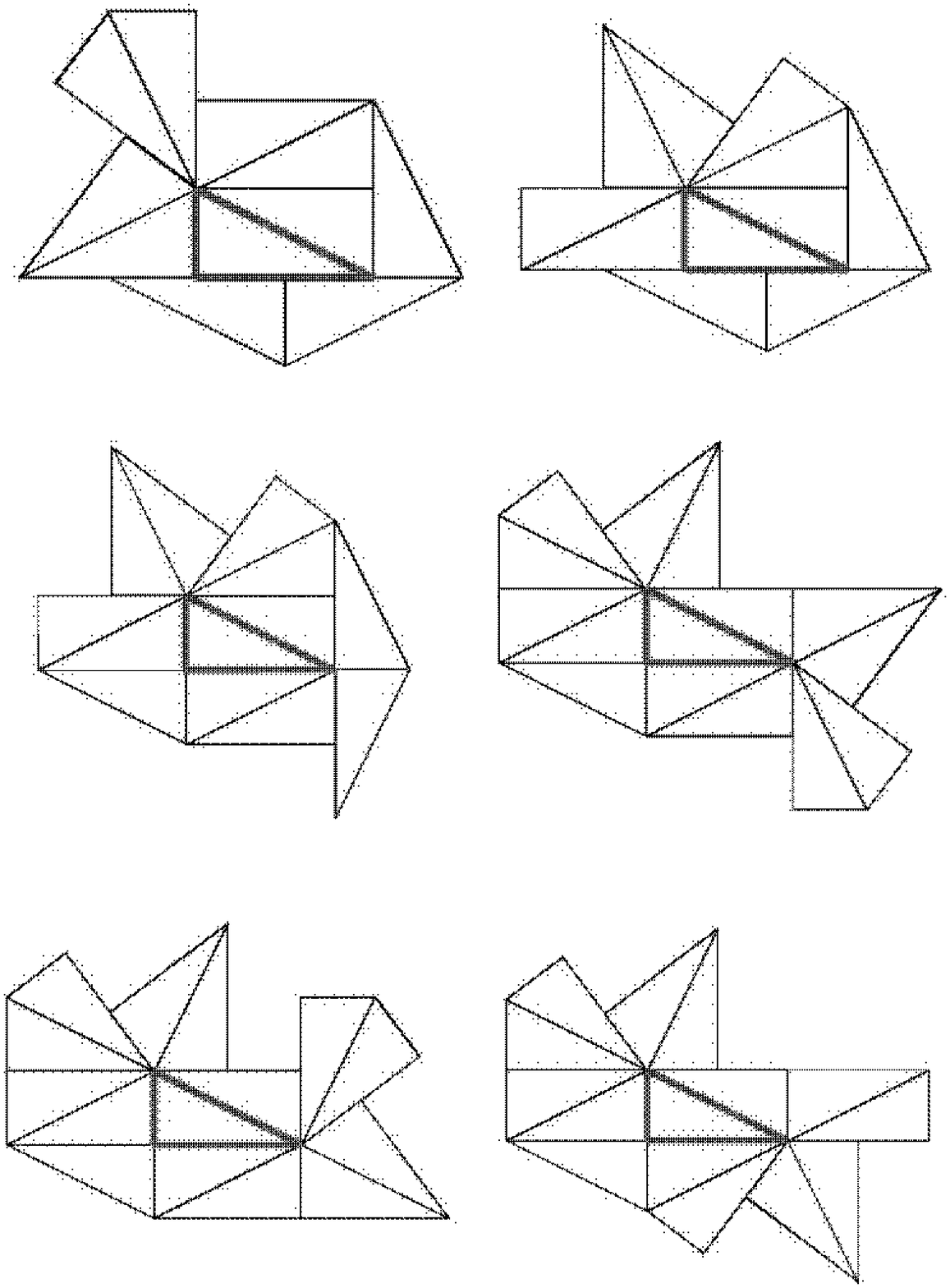}
\hspace{0.2cm}
\includegraphics[scale=0.24]{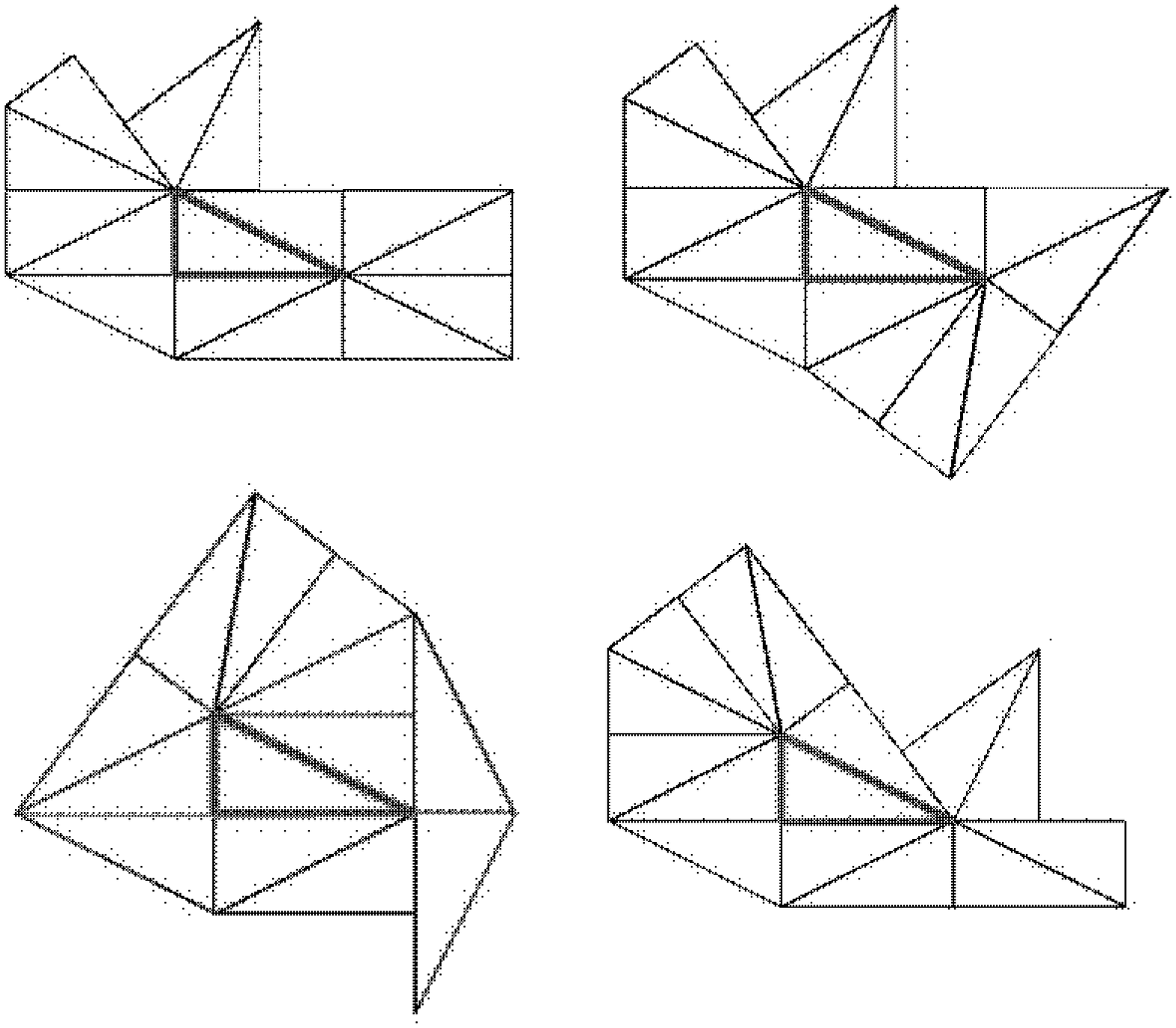}
\hspace{0.2cm}
\includegraphics[scale=0.24]{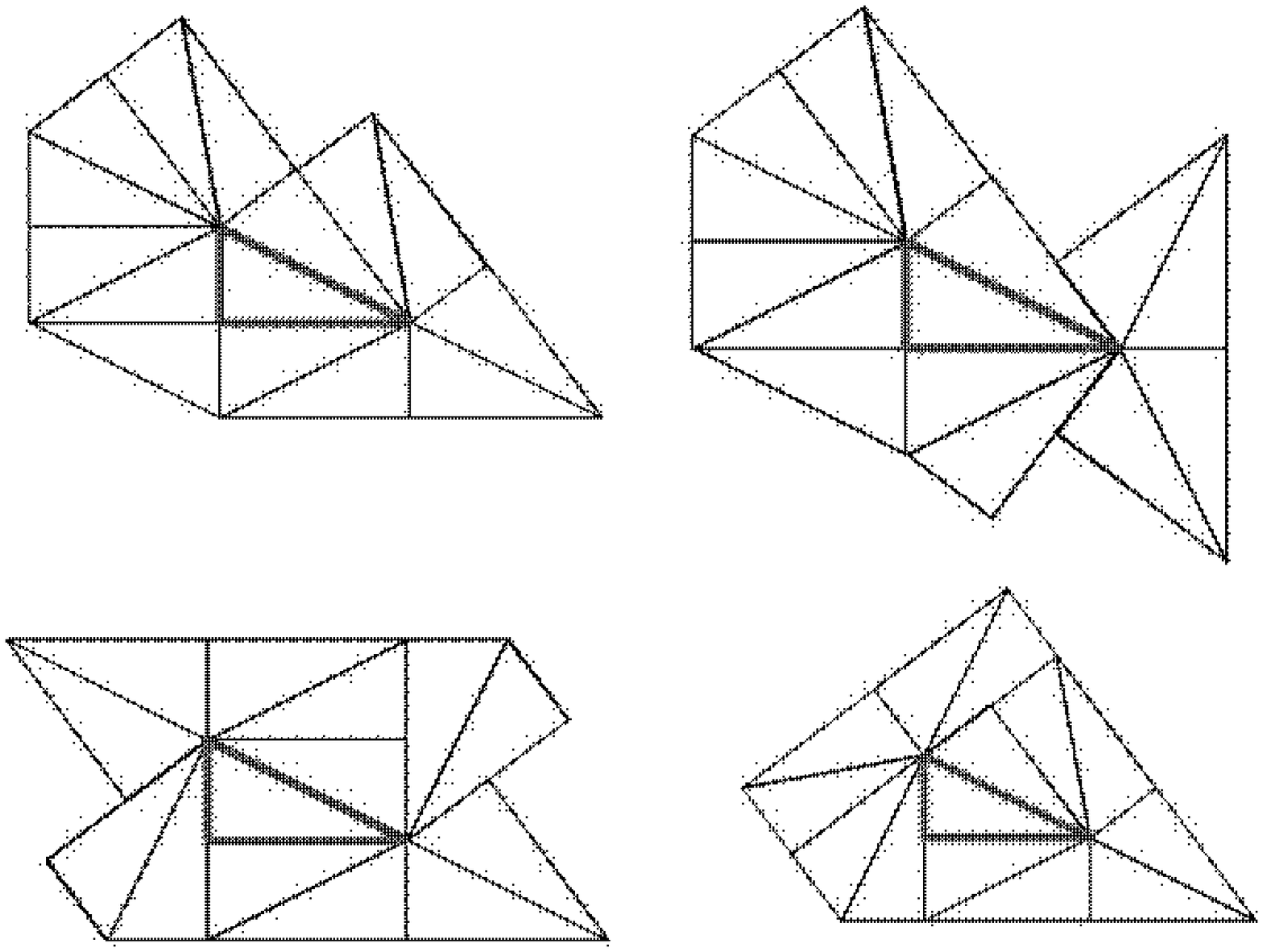}
\hspace{0.5cm}
\includegraphics[scale=0.24]{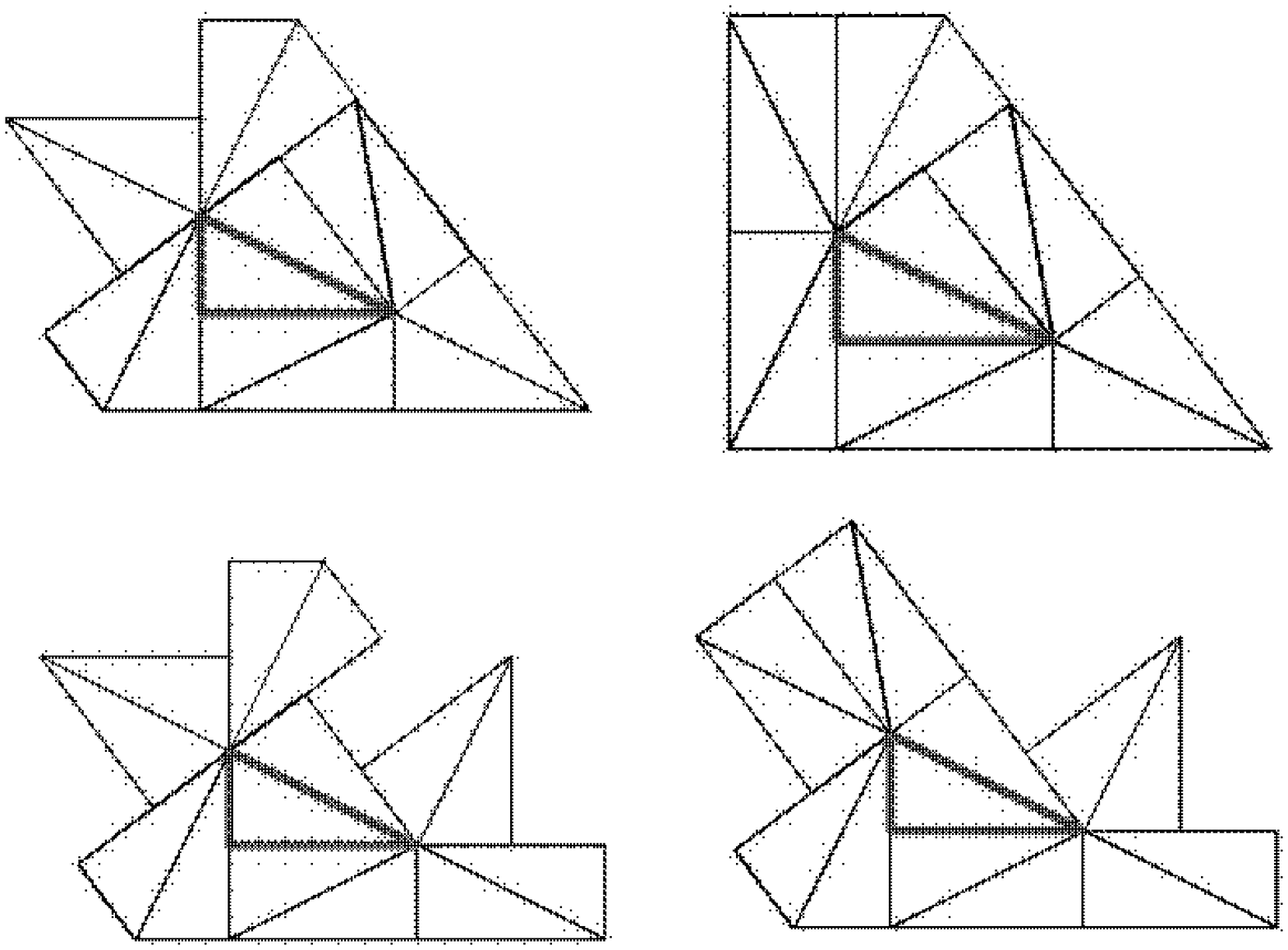}
\hspace{0.5cm}
\includegraphics[scale=0.24]{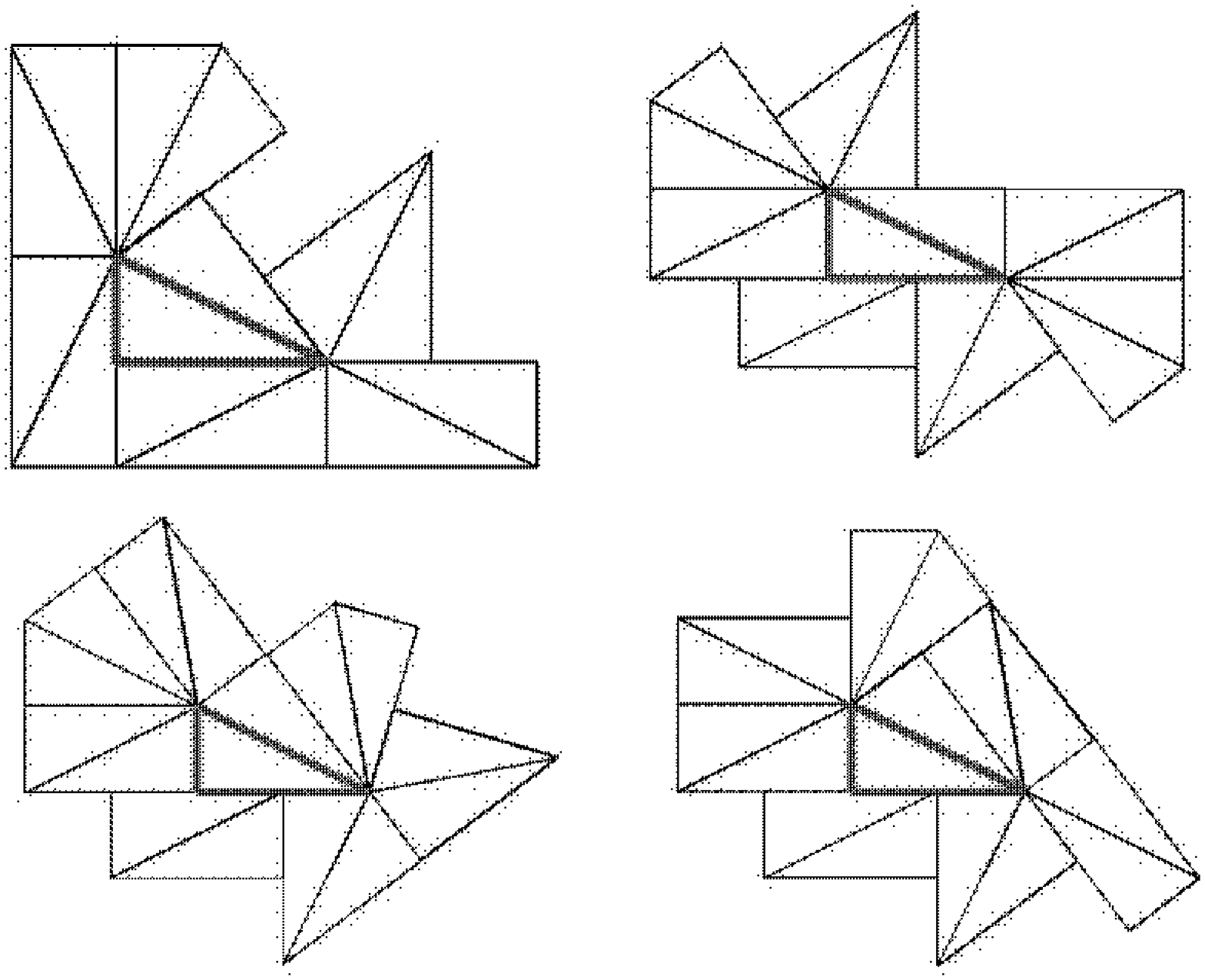}
\hspace{0.5cm}
\end{center}
\vspace{-3.5cm}\caption{28 collared prototiles. \label{prototuilescouronnees}} 
\end{figure}
\end{landscape}

\begin{landscape}
\begin{figure}[ht] 
\vspace{-1cm}
\begin{center}
\includegraphics[scale=0.24]{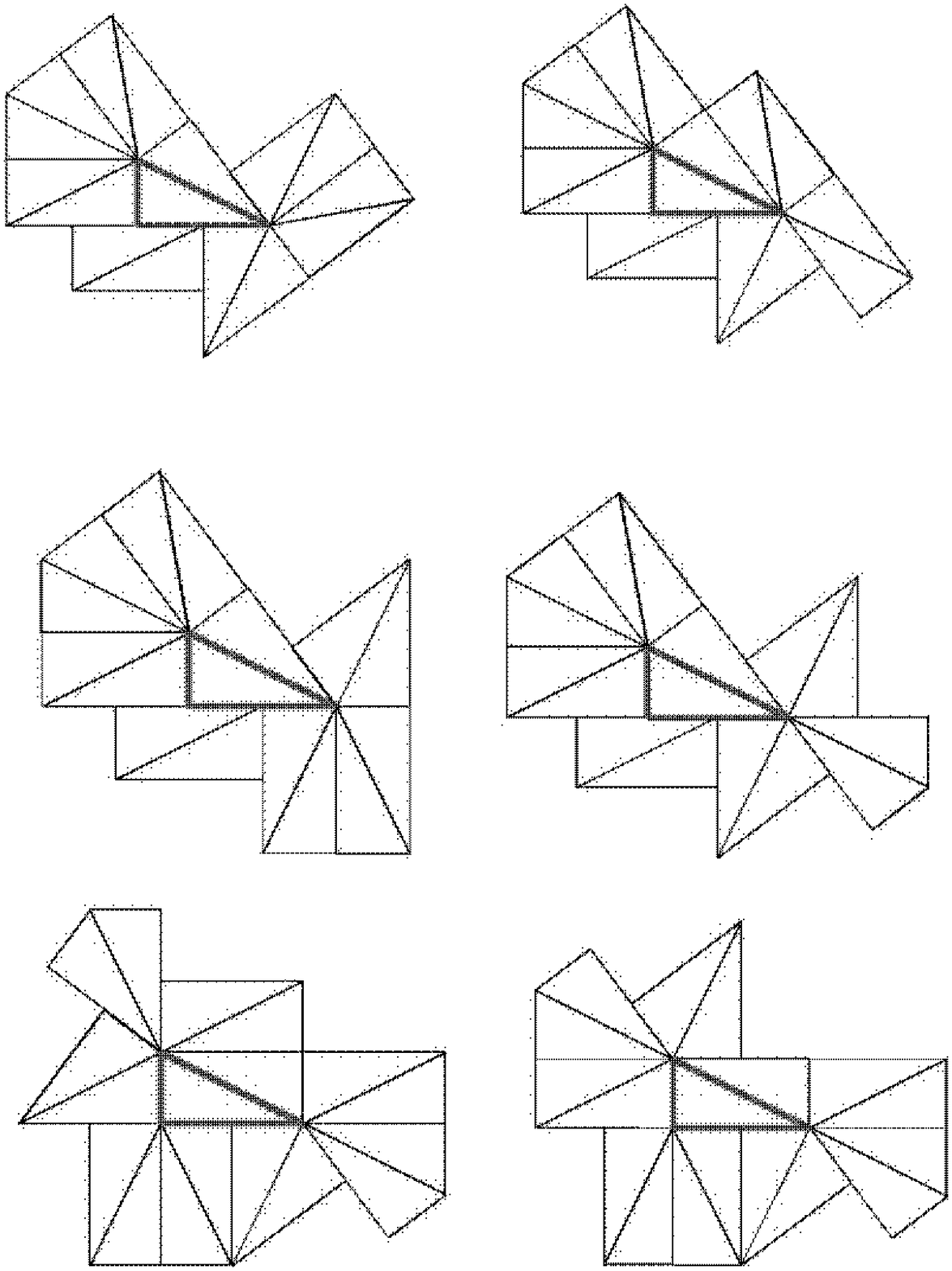}
\hspace{0.5cm}
\includegraphics[scale=0.24]{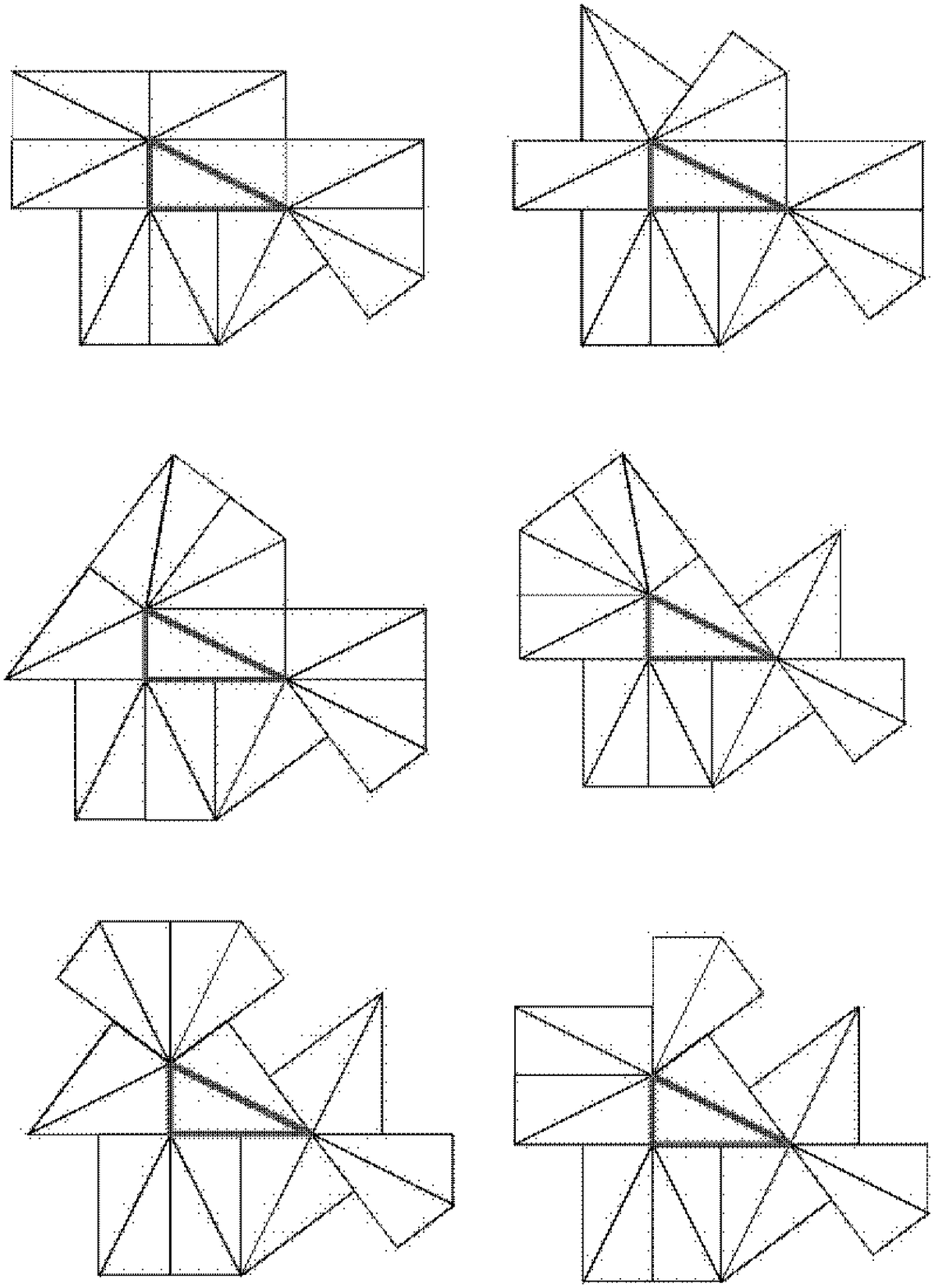}
\hspace{0.5cm}
\includegraphics[scale=0.24]{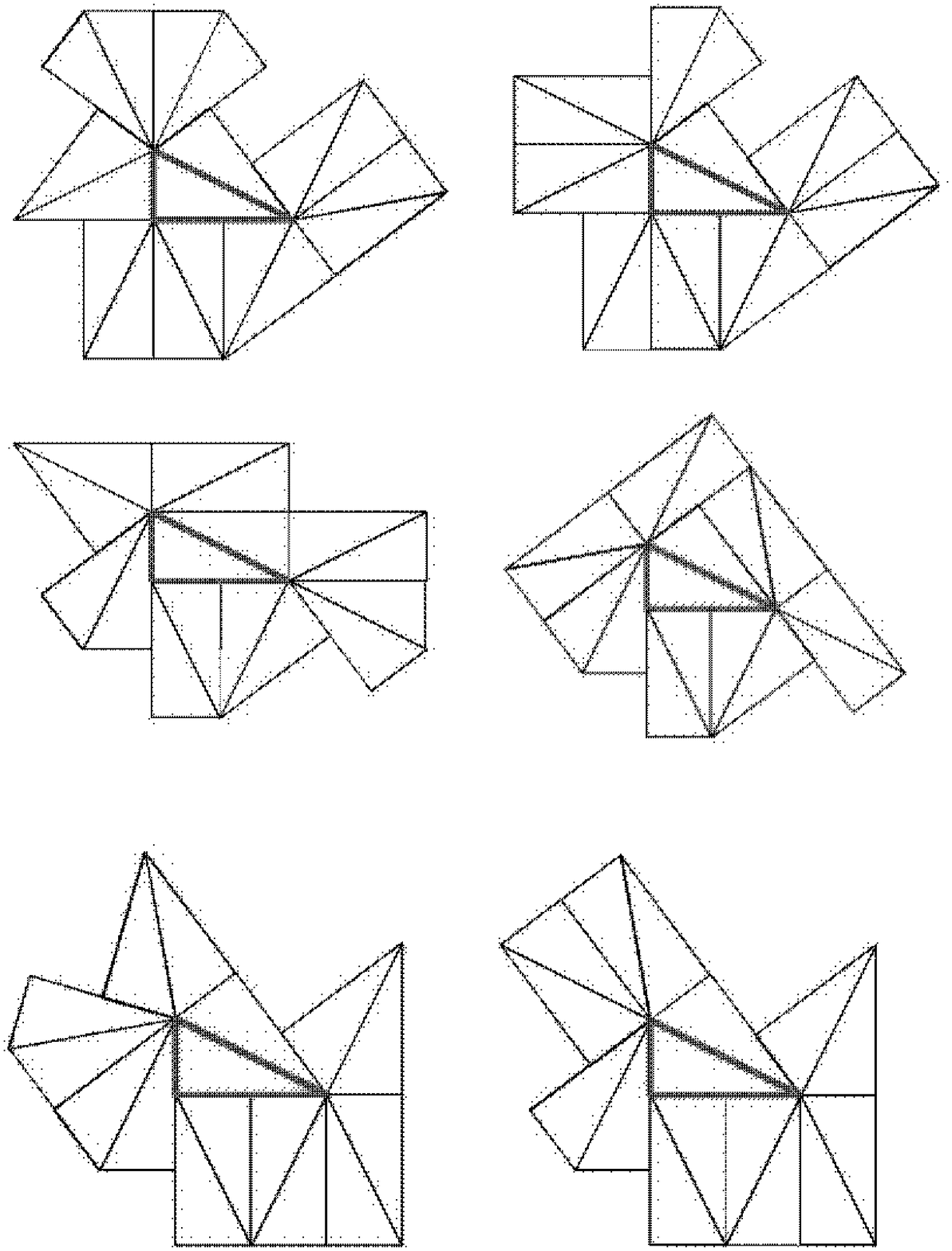}
\hspace*{0.5cm}
\includegraphics[scale=0.24]{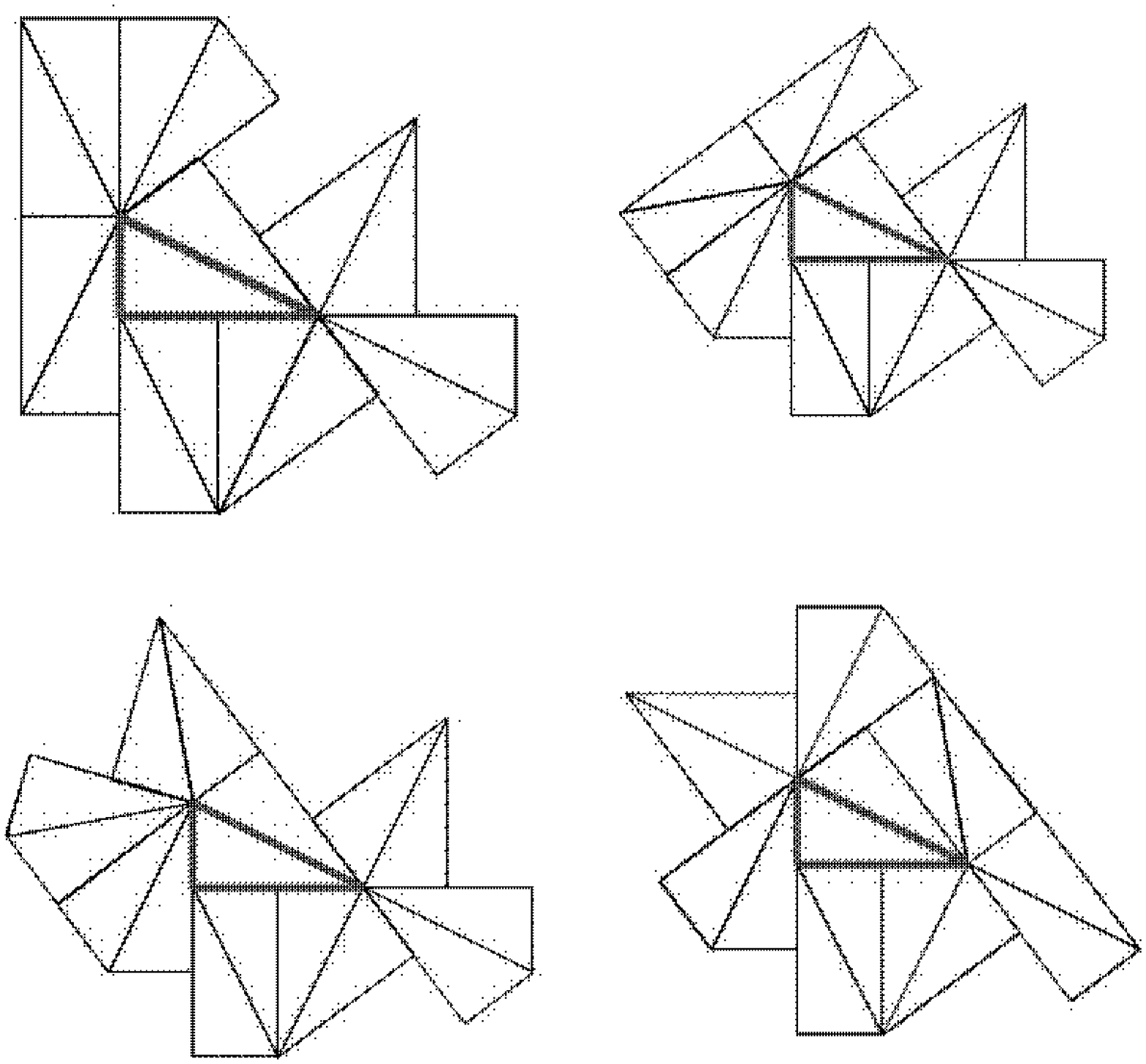}
\hspace*{0.5cm}
\includegraphics[scale=0.34]{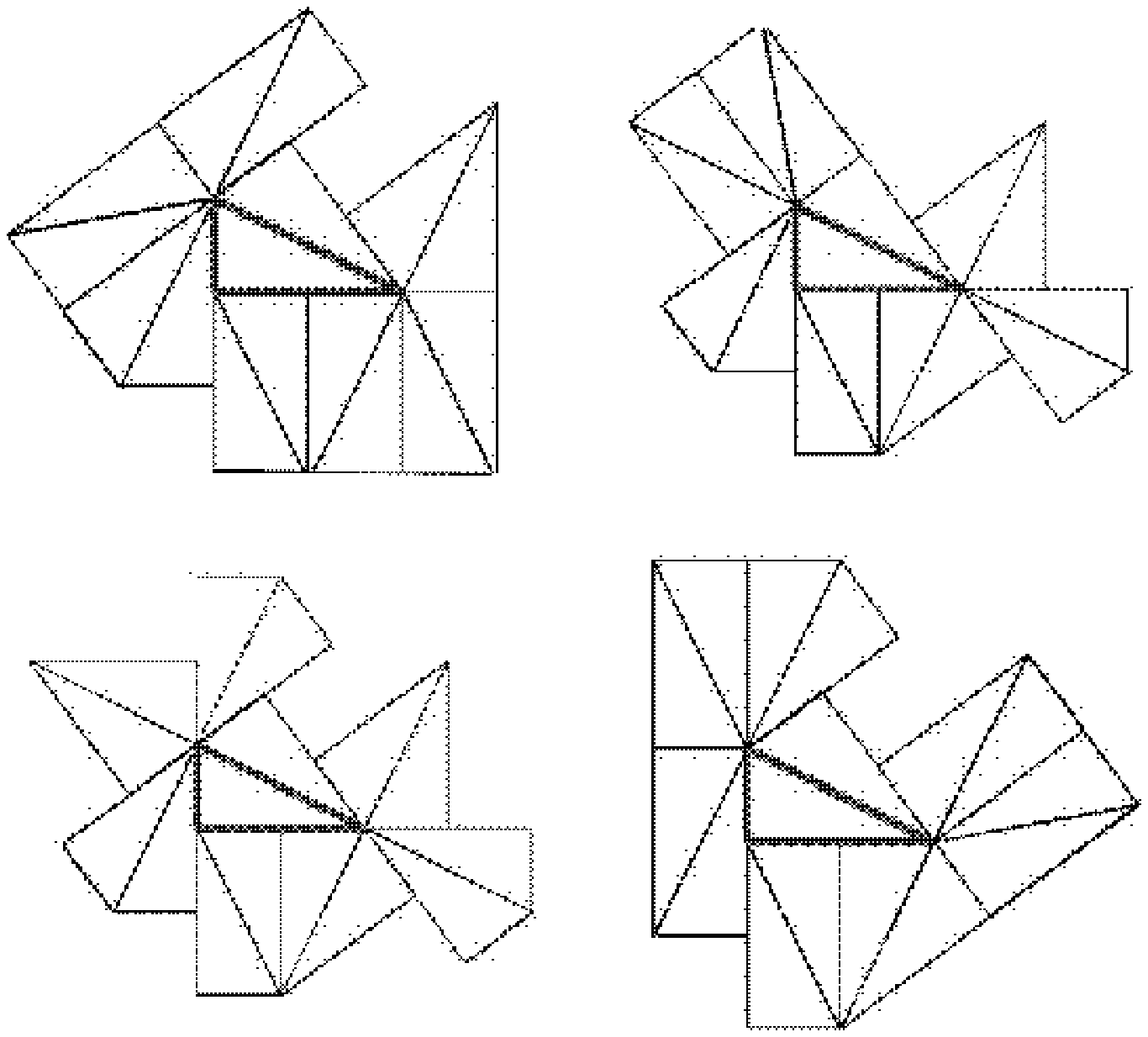}
\end{center}
\vspace{-2.5cm}\caption{26 other collared prototiles. \label{prototuilescouronnees2}}
\end{figure}
\end{landscape}

\noindent
\textbf{Proof :} Let's consider the different cases:
\begin{enumerate}
\item First, if the tiling $\TT'$ which is fixed by $(s,\theta)$ have no fault line (i.e have no infinite supertile), then $s=0$ and $\theta=0 (mod 2 \pi)$.\\
Indeed, let's $x$ be in $E_2$ be such that $\overrightarrow{Ox} = s$ then $O$ and $x$ is in the interior of a $m$-supertile since there isn't infinite supertiles (see p.29 in \citeg{RadSad}).\\
Since no direct isometry fixes our prototiles, $s$ and $\theta$ must be zero.

\item Let's see the case in which $\TT'$ have some infinite supertiles.\\
By \citeg{RadSad} p.30, the number of infinite supertiles in $\TT'$ is bounded by a constant $K$ (in fact for pinwheel tilings, we can take $K=\frac{2\pi}{\alpha}$ where $\alpha$ is the smallest angle in the prototiles).\\
Thus, $\TT'$ doesn't contain more than $K$ infinite supertiles and in fact, it has only a finite number of fault lines. This will give us the result.\\
Indeed, since $(s,\theta)$ fixes $\TT'$, if $F$ is a fault line, $(s,\theta)$ sends it on another fault line $F_1$ in $\TT'$ and thus, $R_\theta$ sends $F$ on a line  parallel to $F_1$.\\
As there is only a finite number of fault lines in $\TT'$, there is $M \in \NN^*$, $m \in \ZZ^*$ such that $M \theta = 2 \pi m$.\\
If we now use results obtained in \citeg{RadSad} p.32, $\theta$ must be in the group of relative orientations $G_{RO}(Pin)$ of pinwheel tilings which is the subgroup of $S^1$ generated by $\frac{\pi}{2}$ and $2\alpha$.\\
Hence, if $\theta=2k\alpha + l \frac{\pi}{2}$ with $k \in \ZZ^*$ and $l \in \ZZ$, it would mean that $\alpha$ is rationnal with respect to $\pi$, which is impossible (see \citeg{Rad} p.664) hence $k=0$ and $\theta \in \left \{ 0, \, \dfrac{\pi}{2}, \, \pi, \, \dfrac{3\pi}{2} \right \}$ $mod(2 \pi)$.\\
Now, if we study the first vertex coronas (i.e the minimal patches around a vertex or around the middle point of the hypothenuses), there is only patches with a $2$-fold symmetry (\citeg{SadPriv} or see \reffig{prototuilescouronnees} and \reffig{prototuilescouronnees2} p.\pageref{prototuilescouronnees} and p.\pageref{prototuilescouronnees2}). Thus, $\theta \in \left \{ 0, \, \pi \right \}$ $mod(2 \pi)$ and if $\theta = 0$, $s=0$ since pinwheel tilings are not fixed by translations.
\begin{flushright}
$\square$
\end{flushright}
\end{enumerate}

\bigskip
\noindent
We note that, in fact, there exists only $6$ pinwheel tilings with a $2$-fold symmetry up to rotations (\citeg{SadPriv}).\\
Hence, there is only $6$ orbits with fixed points for the $\RR^2 \rtimes S^1$ action on $\Omega$.\\
Moreover, there is exactly $6$ circles $F_1, \ldots, F_6$ containing fixed points for the $S^1$ action on $\Omega$ (of course, therefore, the $6$ orbits of these circles contain all the fixed points of the $\RR^2 \rtimes S^1$-action).

\bigskip
\noindent
We thus obtain the following important result on the dynamic of our tiling space:

\begin{theo} (\citeg{BenGam})
The continuous hull is a minimal foliated space.
\end{theo}

\noindent
\textbf{Proof :}
\begin{itemize}
\item[] The proof follows the one in \citeg{BenGam} except that, locally, $\Omega$ looks like an open subset of $S^1$ $\times$ an open subset of $\RR^2$ $\times$ a Cantor set instead of $S^1$ $\times$ an open subset of $\RR^2$ $\times$ a Cantor set, like in \citeg{BenGam}.

\bigskip
\noindent
$\Omega$ is covered by a finite number of open sets $U_i = \phi_i (V_i \times T_i)$ where :
\begin{itemize}
\item[$\bullet$] $T_i$ is a clopen set in $\Xi$;
\item[$\bullet$] $V_i$ is an open subset of $\RR^2 \rtimes S^1$ which read $V_i = \Gamma_i \times W_i$ where $W_i$ is an open subset of $\RR^2$ and $\Gamma_i$ an open subset of $S^1$ of the form $]l\pi/2-\pi/3 ; l\pi/2+\pi/3[$, $l \in \{0,1,2,3\}$;
\item[$\bullet$] $\phi_i:V_i \times T_i \longrightarrow \Omega$ is defined by $\phi_i(v,\omega_0) = \omega_0.v$.
\end{itemize}

\noindent
As we can find finite partitions of $\Xi$ in clopen sets with arbitrarily small diameter, it is possible to choose this
diameter small enough so that:

\begin{itemize}
\item[$\bullet$] the maps $\phi_i$ are homeomorphisms on their images;
\item[$\bullet$] whenever $\TT_1 \in U_i \cap U_j$, $\TT_1 = \phi_i(v,\omega_0) = \phi_j(v',\omega_0')$, the element $v'.v^{-1}$ is independent  of the choice of $\TT_1$ in $U_i \cap U_j$, we denote it by $g_{ij}$.
\end{itemize}
The transition maps read : $(v',\omega_0') = (g_{ij}.v,\omega_0.g_{ij}^{-1})$.\\
It follows that the boxes $U_i$ and charts $h_i = \phi_i^{-1} : U_i \longrightarrow V_i \times T_i$ define a foliated structure on $\Omega$.\\
By construction, the leaves of $\Omega$ are the orbits of $\Omega$ under the $\EE^2$-action.
\end{itemize}

\begin{flushright}
$\square$
\end{flushright}

\bigskip
\noindent
We must do several remarks on the actions.\\
$\EE^2$ isn't acting freely on $\Omega$, even if the translations are, but we could adapt results of Benedetti and Gambaudo obtained in their paper \citeg{BenGam} studying the possible symmetries in our pinwheel tilings.\\
The $\EE^2$-action is not free on $\Omega_0$ too but the $S^1$-action is.

\bigskip
\noindent
Using the group of relative orientations $G_{RO}(Pin)$, we can see that each $\RR^2$-orbit of $\Omega$ is in fact a dense subset of $\Omega$ (see \citeg{HolRadSad}).

\bigskip
\noindent
\subsection{$C^*$-algebra associated to a tiling}
\medskip
\noindent
As in the paper of Kellendonk and Putnam, one can naturally associate a $C^*$-algebra to pinwheel tilings.\\
It is the crossed product $C^*$-algebra $\mathcal{A}=C(\Omega) \rtimes (\RR^2 \rtimes S^1)$ obtained from the dynamical system $(\Omega, \RR^2 \rtimes S^1)$.

\bigskip
\noindent
For the sake of the next session, let's consider a topological space $\Omega$ endowed with a right action of a locally compact group $G$ (for pinwheel tilings, $G=\RR^2 \rtimes S^1$). For simplicity, let's suppose that $G$ is unimodular and let's fix a Haar measure $\lambda$ on $G$.\\
In order to construct $C(\Omega) \rtimes G$, let's first consider the vector space $C_c \left( \Omega \times G \right)$ of continuous functions with compact support on $\Omega \times G$. One endows this space with a convolution product and an involution :
$$ f*g(\omega,g) := \int_{G} f(\omega,h) g(\omega.h,h^{-1}g) d\lambda(h)$$
$$f^*(\omega,g) := \overline{f(\omega.g,g^{-1})}$$
with $f,g\in C_c(\Omega \times G)$ and $\omega \in \Omega$, $g \in G$.\\
One then defines a norm on it :
$$\parallel f \parallel_{\infty,1} = Max \left \{  \sup_{\omega \in \Omega} \int_{G} \mid \! f(\omega,g) \! \mid d\lambda(g) , \sup_{\omega \in \Omega} \int_{G} \mid \! f^*(\omega,g) \! \mid d\lambda(g) \right \}$$
Define the involutive algebra $L^{\infty,1}(\Omega \times G)$ as the completion of $C_c(\Omega \times G)$ under this norm.\\
This algebra is represented on $L^2(G)$ by the family $\{ \pi_\omega , \omega \in \Omega \}$ of representations given by :
$$\pi_\omega (f)\psi(g) := \int_{G} f(\omega.g,g^{-1}.h)\psi(h) d\lambda(h), \quad \quad \psi \in L^2(G).$$
Namely, $\pi_\omega$ is linear, $\pi_\omega(fg)=\pi_\omega(f)\pi_\omega(g)$, $\pi_\omega(f)^*=\pi_\omega(f^*)$ and $\pi_\omega$ is boun\-ded, $\parallel \! \pi_\omega(f) \! \parallel \leqslant \parallel \! f \! \parallel_{\infty,1}$.\\
Set $\parallel \! f \! \parallel = \sup_{\omega \in \Omega} \parallel \! \pi_\omega(f) \! \parallel$ which defines a $C^*$-norm on $L^{\infty,1}(\Omega \times G)$ and permits to define $C(\Omega) \rtimes G$ as the completion of $C_c(\Omega \times G)$ or of $L^{\infty,1}(\Omega \times G)$ under this norm.

\bigskip
\noindent
For pinwheel tilings, $G=\RR^2 \rtimes S^1$ and this $C^*$-algebra $C(\Omega) \rtimes (\RR^2 \rtimes S^1)$ is isomorphic to $(C(\Omega) \rtimes \RR^2) \rtimes S^1$ (\citeg{Chab} p.15).\\
$S^1$ acts on $C_c(\Omega \times \RR^2)$ by : $\theta . f(w,s) := f(R_{-\theta}(w),R_{-\theta}(s))$ for $f \in C_c(\Omega \times \RR^2)$.\\
This action extends to the crossed product $C(\Omega) \rtimes \RR^2$.

\bigskip
\noindent
This algebra is interesting for two reasons. First, it contains dynamical informations related to the combinatorial properties of the tiling. Secondly, the gap-labeling conjecture made by J. Bellissard in \citeg{Bel1} links the electronic motion in a quasicrystal to the $K_0$-group of the $C^*$-algebra of the tiling. This is the subject of the next session.\\
A natural question is to know if, in the case of pinwheel tilings, we can relate the $K$-theory of $C(\Omega) \rtimes \RR^2 \rtimes S^1$ to the $\ZZ$-module of patch frequencies of the pinwheel tiling (see the end of next session for motivations).

\bigskip
\noindent
\section{Gap-labeling}
\subsection{Quasicrystals}
\medskip
\noindent
Before beginning the presentation of the works of Bellissard on the gap-labeling conjecture, we introduce the definition of aperiodic tilings which model quasicrystals. It is a tiling of the space $\RR^n$ with the same combinatorial properties as pinwheel tilings.

\medskip
\noindent
\begin{defi}
A \textbf{tiling} of $\RR^n$ is a countable family $P = \{t_0, t_1, \ldots \}$ of non empty compact sets $t_i$, called \textbf{tiles} (each tile is suppoed to be homeomorphic to the unit ball), such that :
\begin{itemize}
\item[$\bullet$] ${\displaystyle \bigcup_{i \in \mathbb{N}} t_i = E_n}$ where $E_n$ is the Euclidean space $\RR^n$ endowed with a fixed origin $O$;
\item[$\bullet$] Tiles meet only on their border;
\item[$\bullet$] The interior of two tiles are disjoint.
\end{itemize}

\bigskip
\noindent
Let $G$ be a subgroup of the group of rigid motions of $E_n$.\\
We consider in the sequel the particular case where there exists a finite family $\{p_0, \ldots, p_n \}$ such that each compact set $t_i$ is a $G$-copy of some $p_k$. The $p_k$'s are then called \textbf{prototiles}.\\
A tiling is of \textbf{finite $\mathbf{G}$-type} if, for any $R>0$, there are only finitely many patches of diameter less than or equal to $R$ upto the $G$-action.\\
A tiling $\TT'$ is \textbf{$\mathbf{G}$-repetitive} if, for each patch $\mathcal{A}$ in $\TT'$, there exists $R(\mathcal{A})>0$ such that every ball of radius $R(\mathcal{A})$ in $\RR^n$ meets $\TT'$ on a patch containing a $G$-copy of $\mathcal{A}$.\\
A \textbf{$\mathbf{G}$-aperiodic tiling} is a tiling of $\RR^n$ of finite $G$-type, $G$-repetitive and nonperiodic with respect to translations of $\RR^n$.
\end{defi}

\bigskip
\noindent
For example, pinwheel tilings are aperiodic tilings with $n=2$ and $G=\RR^2 \rtimes S^1$, the group of direct isometries of the plane.\\
Many more examples are known, the most famous being the one constructed by Penrose ($n=2$, $G=\RR^2$) (see \citeg{Pet}).\\
In this section on the gap-labeling, we fix a $\RR^n$-aperiodic tiling $\TT$.

\bigskip
\noindent
One can then associate to such a tiling a topological space $\Omega$ as we have done for pinwheel tilings. Let $\TT.\RR^n$ be the space of all the translations of $\TT$. $\Omega$ is then the completion of this space for the metric $d$ defined as follows :
$$\hspace*{-2cm} A=\left \{\varepsilon \in \Big [ 0,\tfrac{1}{\sqrt{2}} \Big ] \; / \; \exists s,s' \in B^2_\varepsilon(0) \, , \, \theta,\theta' \in B^1_\varepsilon(0) \text{ s.t. } \right.$$
\vspace{-0.6cm} 
$$\hspace*{3cm} \left. \TT_1.(s,\theta) \cap B_{\frac{1}{\varepsilon}}(O) =  \TT_2.(s',\theta') \cap B_{\frac{1}{\varepsilon}}(O) \right \}$$
where $B_{\frac{1}{\varepsilon}}(O)$ is the ball of $E_n$ centered at $O$ with radius $\frac{1}{\varepsilon}$ and $B^n_\varepsilon(0)$ is the ball of $\RR^n$ centered at $0$ with radius $\varepsilon$ (i.e one only considers translations near $Id$) .

\bigskip
\noindent
$d$ is then defined by : 
$$d(\TT_1,\TT_2) = \left\{
 \begin{array}{cc}
Inf A & \text{ if } A \neq \emptyset\\
\frac{1}{\sqrt{2}} & \text{ else } \quad
\end{array} \right..$$

\noindent
Then $\Omega(\TT):=\overline{\TT.\RR^n}$ is the completion of $\TT.\RR^n$ under $d$.\\
One can prove that $\Omega(\TT)$ is a compact metric space, that every $\RR^n$-orbits are dense in this space and thus, all the tilings of $\Omega(\TT)$ have the same completion.\\
One can then define, as in the first section, the canonical transversal by fixing a point in each prototile and letting the \textbf{canonical transversal} $\Xi$ be the set of all $\TT' \in \Omega$ with one point on the origin $O$.\\
This space is still a Cantor set (see \citeg{BenGam}) and $\Omega$ is a foliated space with leaves homeomorphic to $\RR^n$.

\bigskip
\subsection{Integrated density of states - IDS} \label{IDS}

\bigskip
\noindent
This section summarizes the works of Bellissard \citeg{Bel1} (see also \citeg{Ypma}).\\
The gap-labeling conjecture describes qualitatively the spectrum of the operator associated to the electronic motion in a quasicrystal (which is a solid with a particular atomic distribution, modeled by aperiodic tilings).\\
This motion is described by the Schr\"odinger operator on the Hilbert space $L^2(\RR^n)$ :
$$H = -\dfrac{\hslash}{2m} \Delta + V(x)$$
where $\hslash$ is the Dirac constant, equals to $\dfrac{h}{2\pi}$ where $h$ is the Planck constant, $m$ is the electron mass, $\Delta$ the Laplacian on $\RR^n$ and $V$ is a potential in $L^\infty(\RR^n)$ depending on the atomic distribution of the quasicrystal.\\
The domain of $H$ is $\mathcal{D}(H) = \{ \psi \in L^2(\RR^n) \mid \Delta \psi \in L^2(\RR^n) \}$.

\bigskip
\noindent
Let fix a $\RR^n$-aperiodic tiling $\TT$ which is a model for the quasicrystal. One can then consider the two spaces, introduced in the last section, $\Omega$ and $\Xi$ together with the $C^*$-algebra, associated to the dynamical system, $C(\Omega) \rtimes \RR^n$ called the \textbf{$\mathbf{C^*}$-algebra of observables}.

\bigskip
\begin{defi}
A \textbf{covariant} family of selfadjoint operators $\{H_\omega ; \omega \in \Omega \}$ satisfies $U_a.H_\omega.U^*_a = H_{\omega-a}$ for all $x \in \RR^n$ and $U_a:L^2(\RR^n) \rightarrow L^2(\RR^n)$ defined by $U_a(f)(x):=f(x+a)$.\\
Such a family is \textbf{affiliated} to the $C^*$-algebra $A=C(\Omega) \rtimes \RR^n$ if, for every $f \in C_0(\RR)$, the bounded operator $f(H_\omega)$ can be represented as $\pi_\omega(h_f)$ for some $h_f \in A$ such that the map $h:C_0(\RR) \rightarrow A : f \mapsto h_f$ is a bounded *-homomorphism.
\end{defi}
\noindent
In the sequel, we fix such a covariant family of selfadjoint operators $\{H_\omega ; \omega \in \Omega \}$ affiliated to $C(\Omega) \rtimes \RR^n$ (see Bellissard's papers \citeg{Bel1} and \citeg{BelHerZar} for the existence of such family).\\
We are then interested in the spectrum of $H_\omega$.\\
For this, Bellissard linked a physical function defined on $\RR$ (the Integrated Density of States IDS)  to a map obtained from $C(\Omega) \rtimes \RR^n$.

\bigskip
\noindent
Let's define the IDS $\mathcal{N}(E)$ which can be seen as the number of eigenvalues per unit volume up to $E$.

\begin{defi}
Let $G$ be a locally compact group. A \textbf{F$\mathbf{\scriptstyle{\varnothing}}$lner sequence} is a sequence $(\Lambda_n)$ of open subsets of $G$, each with finite Haar measure $\mid \Lambda_n \mid$, such that $G = \cup \Lambda_n$ and such that for all $x \in G$,
$$ \lim_{n \rightarrow \infty} \frac{\mid \Lambda_n \Delta x.\Lambda_n \mid}{\mid \Lambda_n \mid} = 0$$ 
where $V \Delta W = (V \cup W) \setminus (V \cap W)$.
\end{defi}

\bigskip
\noindent
It can be shown that there exists a F$\scriptstyle{\varnothing}$lner sequence in $G$ if and only if $G$ is amenable (see \citeg{Gre}).\\
One can now define
$$\mathcal{N}_{\omega,\Lambda_n}(E) := \# \{E' \in Sp(H_{\omega,\Lambda_n}) \mid E' \leqslant E \},$$
where $H_{\omega,\Lambda_n}$ is the Hamiltonian $H_\omega$ restricted to $\Lambda_n$, acting on the Hilbert space $L^2(\Lambda_n)$, subject to certain boundary conditions (see \citeg{Bel1}).\\
The IDS $\mathcal{N}_\omega : \RR \longrightarrow \RR^+$ is then defined by
$$\mathcal{N}_\omega(E) := \lim_{n \rightarrow \infty} \frac{1}{\mid \Lambda_n \mid} \mathcal{N}_{\omega,\Lambda_n}(E).$$

\noindent
Jean Bellissard proved that this limit exists and is independent of the chosen boundary conditions (\citeg{Bel1}).

\bigskip

\subsection{Shubin's formula}

\bigskip
\noindent
The link between the IDS and the $C^*$-algebra $C(\Omega) \rtimes \RR^n$ is given by the Shubin's formula.\\
On a Hilbert space $\mathcal{H}$ with an orthonormal basis $\{e_i\}$, one can define the operator trace of a bounded operator $A$ by : 
$$ Tr(A) := \sum_{i=1}^{\infty} \langle e_i,Ae_i\rangle .$$
This trace is independent of the choice of a basis in $\mathcal{H}$ if $A$ is \textit{trace-class}, i.e if $Tr(\mid A \mid) < \infty$ where $\mid A \mid  := \sqrt{A^*A}$.\\
One then has :
$$\mathcal{N}_\omega(E) = \lim_{n \rightarrow \infty} \frac{1}{\mid \Lambda_n \mid} Tr_{\Lambda_n} \left ( \chi_{]-\infty;E]}(H_{\omega,\Lambda_n}) \right )$$
where $Tr_{\Lambda_n}$ is the restriction to $L^2(\Lambda_n)$ of the operator trace on $L^2(\RR^n)$ and $\chi_{]-\infty;E]}$ is the characteristic function of $]-\infty;E]$.\\
To link this formula to the $C^*$-algebra of observables, one defines an ergodic invariant measure on $\Omega$. This is given by the next proposition which is a consequence of the amenability of $\RR^n$ and of the Krein-Milman theorem :

\begin{prop}
There exists a  translation invariant, ergodic probability measure $\mu$ on $\Omega$. 
\end{prop}

\noindent
A trace is then densely defined on $C(\Omega) \rtimes \RR^n$ by : $\forall \, f \in C_c(\Omega \times \RR^n)$,
$$\tau^\mu(f):= \int f(\omega,0) d\mu(\omega).$$
Since $\mu$ is translation invariant, $\tau^\mu$ has the properties of a positive trace i.e $\tau^\mu(f*f^*) \geqslant 0$ and $\tau^\mu(f*g)=\tau^\mu(g*f)$.\\
This functional can be extended in a faithful trace (which follows from the ergodicity of $\mu$ and the minimality of the action of $\RR^n$ on $\Omega$) and semi-finite (see \citeg{MooSch} p150 and following) on the Von Neumann algebra of $C(\Omega) \rtimes \RR^n$. Moreover, this trace is finite on projections of the Von Neumann algebra (see \citeg{MooSch} p154 and \citeg{Ped} 5.6.7).

\bigskip
\noindent
The link between the IDS and this trace is highlighted by the following Birkhoff theorem :

\begin{theo}
Let $\Omega$ be a compact space with a probability measure $\mu$ that is ergodic and invariant under the action $T$ of an amenable group $G$ on $\Omega$.\\
Then the left Haar measure $\lambda$ on $G$ can be normalized in such a way that for all $f \in C(\Omega)$, for $\mu$-almost all $\omega \in \Omega$, we have :
$$ \int_\Omega f(\omega') d\mu(\omega') = \lim_{n \rightarrow \infty} \frac{1}{\mid \Lambda_n \mid} \int_{\Lambda_n} f(T_g\omega)d\lambda(g),$$
where $(\Lambda_n)$ is a F$\scriptstyle{\varnothing}$lner sequence in $G$.
\end{theo}

\bigskip
\noindent
Applying this theorem, we get for all $f \in C_c(\Omega \times \RR^n)$ and $\mu$-almost all $\omega \in \Omega$ :
$$\tau^\mu(f) =  \displaystyle{\lim_{n \rightarrow \infty} \frac{1}{\mid \Lambda_n \mid} \int_{\Lambda_n} f(\omega.g,0)d\lambda(g)}.$$
By definition of $\pi_\omega$ and $Tr$, the formula of the theorem becomes :
\begin{eqnarray}
\tau^\mu(f) & =  & \displaystyle{\lim_{n \rightarrow \infty} \frac{1}{\mid \Lambda_n \mid} Tr_{\Lambda_n}(\pi_\omega(f))}, \text{ for } \mu-\text{almost all } \omega.
\end{eqnarray}

\bigskip
\noindent
Note that, if $E \in \RR$, then $\chi_{]-\infty;E]}(H_\omega) = \chi_{H \leqslant E}$ for some $\chi_{H \leqslant E}$ in the Von Neumann algebra of the $C^*$-algebra of observables.\\
Moreover, if $E \in \mathbf{\mathfrak{g}}$ for some gap $\mathbf{\mathfrak{g}}$ in the spectrum of $H_\omega$ (i.e a connected component of the complement of the spectrum) then $\chi_{]-\infty;E]}$ is a bounded continuous map on $Sp(H_\omega)$ and thus $\chi_{]-\infty;E]}(H_\omega)$ can be represented as $\pi_\omega(\chi_{H \leqslant E})$ for some element $\chi_{H \leqslant E} \in C(\Omega) \rtimes \RR^n$, since $H_\omega$ is affiliated to this $C^*$-algebra. 

\bigskip
\noindent
There is then a link between a physically measurable function by experiments and our $C^*$-algebra :

\begin{defi}
A covariant family $(H_\omega)$ affiliated to $C(\Omega) \rtimes \RR^n$ is said to satisfy the \textbf{Shubin's formula} if for $\mu$-almost all $\omega \in \Omega$, we have:
$$\mathcal{N}_\omega (E) = \tau^\mu (\chi_{H \leqslant E}).$$
\end{defi}
The common value is denoted by $\mathcal{N}(E)$.

\subsection{Gap-labeling}

\bigskip
\noindent
In this section, we fix a covariant family of selfadjoint operators $(H_\omega)$ affiliated to $C(\Omega) \rtimes \RR^n$ and satisfying the Shubin's formula (see \citeg{Bel1} for the existence of such a family).\\
As we said, if $E \notin Sp(H_\omega)$, $\chi_{]-\infty;E]}(H_\omega)$ is a projection of the $C^*$-algebra.\\
Moreover, if $\mathbf{\mathfrak{g}}$ is a gap in $Sp(H_\omega)$ then for any two values $E, E' \in \mathbf{\mathfrak{g}}$, we have $\chi_{]-\infty;E]}(H_\omega) = \chi_{]-\infty;E']}(H_\omega)$. Thus, there is a labeling of the gaps $\mathbf{\mathfrak{g}}$ by projections $P(\mathbf{\mathfrak{g}})$ of  $C(\Omega) \rtimes \RR^n$ : $P(\mathbf{\mathfrak{g}}) := \chi_{]-\infty;E]}(H_\omega)$ for any $E \in \mathbf{\mathfrak{g}}$.\\
Now, an important property of the trace $\tau^\mu$ is its invariance under unitary transformations which implies that it induces a trace on the set of equivalence classes of projections under unitary transformations. Thus, $\tau^\mu$ induces a linear map $\tau^\mu_*$ on $K_0(C(\Omega) \rtimes \RR^n)$.\\
We then have another form for the Shubin's formula on gaps of $Sp(H_\omega)$ : 
$$\mathcal{N}(\mathbf{\mathfrak{g}}) = \tau^\mu_* [P(\mathbf{\mathfrak{g}})].$$
Using the fact that $C(\Omega) \rtimes \RR^n$ is separable, we can state :

\begin{prop}
On gaps of $Sp(H_\omega)$, the Integrated Density of States takes its values in $\tau^\mu_* \left( K_0(C(\Omega) \rtimes \RR^n) \right)$, which is a countable subset of $\RR$.
\end{prop}

\bigskip
\noindent
The goal is then to compute the image under $\tau^\mu_*$ of the $K$-theory of $C(\Omega) \rtimes \RR^n$.

\bigskip
\noindent
For $\RR^n$-aperiodic tilings, this image was conjectured in \citeg{BelHerZar} in 2000 and was proven independently by Bellissard, Benedetti and Gambaudo \citeg{BelBenGam} on one hand, Benameur and Oyono-Oyono \citeg{BenOyo} on the other and finally by Kaminker and Putnam \citeg{KamPut}:

\begin{theo}
Let $\Omega$ be the continuous hull of a $\RR^n$-aperiodic tiling with a totally disconnected canonical transversal $\Xi$. Let $\mu$ be a translation invariant, ergodic probability measure and $\mu^t$ the induced measure on $\Xi$. We then have : 
$$\tau^\mu_* \Big(K_0 \big( C(\Omega) \rtimes \RR^n \big) \Big)  = \mu^t \big( C(\Xi,\ZZ) \big)$$
where $C(\Xi,\ZZ)$ is the set of continuous functions on $\Xi$ with integer values and 
$$\mu^t \big( C(\Xi,\ZZ) \big) := \left \{ \int_\Xi f d\mu^t ,\, f \in C(\Xi,\ZZ) \right \}.$$
\end{theo}

\noindent
A natural question is whether this theorem remains true for pinwheel tilings i.e :
\begin{theo}
Let $\Omega$ be the continuous hull of a pinwheel tiling, $\mu$ a $\RR^2 \rtimes S^1$-invariant, ergodic probability measure and $\mu^t$ the induced measure on $\Xi$. We have : 
$$\tau^\mu_* \Big( K_0 \big( C(\Omega) \rtimes \RR^2 \rtimes S^1 \big) \Big) = \mu^t \big( C(\Xi,\ZZ) \big).$$
\end{theo}

\noindent
In other words, is there still this link between the image of the $K$-theory of the $C^*$-algebra of pinwheel tilings under the linear map induced by the trace and the $\ZZ$-module of patch frequencies?\\
In this paper, we make a first step in the direction of this theorem. More precisely, we link the image of the $K$-theory to the tangential cohomology group  of the continuous hull (see \citeg{MooSch}).\\
This cohomological part is more computable as could be seen in \citeg{Hai-coho} and the gap-labeling is proved in this paper, showing moreover that the module of patch frequencies is given explicitly by $\frac{1}{264}\ZZ\left[\frac{1}{5} \right]$.

\bigskip
\noindent
To make this link between $K$-theory and cohomology, we follow the ideas of \citeg{BenOyo}.\\
Let's remind the important steps of their proof. In their article, Benameur and Oyono-Oyono used the index theorem for foliations established by Alain Connes in \citeg{Con1} and more precisely, the version for foliated spaces proved in \citeg{MooSch}.\\
They then proceeded in several steps.\\
First, the inclusion $\mu^t(C(\Xi,\ZZ)) \subset \tau^{\mu^t}_*(K_0(C(\Xi) \rtimes \ZZ^n))$ is easy and thus they only have to prove the inclusion in the other side. Then, they notice that it is enough to prove it for even integers $n$.\\
Next, they prove that the $K$-theory of $C(\Omega)$ is isomorphic to the one of $C(\Xi) \rtimes \ZZ^n$ and that this isomorphism is given by the map $e \mapsto Ind_\Omega(\partial^e_{\Xi,\RR^n})$ where $\partial^e_{\Xi,\RR^n}$ is a Dirac operator with coefficients in the fiber bundle associated to $e$ and $Ind_\Omega$ is the analytical index.\\ 
The second crucial point of the proof is the fact that the top dimensional tangential cohomology group is isomorphic, by a map $\Psi_{\ZZ^n}$, to the integer group of coinvariants $C(\Xi,\RR)_{\ZZ^n}$ of $C(\Xi,\RR)$ under the action of $\ZZ^n$.
The last step is to link $\tau^{\mu^t}_* \Big( K_0 \big( C(\Xi) \rtimes \ZZ^n \big) \Big)$ to this cohomology group using  :

\begin{theo}
$$\tau^\mu_*(Ind_\Omega(\partial^e_\Xi)) = \langle ch^n_l([e]),[C_{\ZZ^n,\mu}] \rangle$$
where $C_{\ZZ^n,\mu}$ is some current on $\Omega$ and $ch^n_l([e])$ is an element of the top cohomology group (it is the image of $[e]$ under the component of degree $n$ of the tangential Chern character).
\end{theo}

\noindent
Finally, they prove that the image of $ch^n_l([e])$ under $\Psi_{\ZZ^n}$ is integer valued i.e $\Psi_{\ZZ^n}(ch^n_l([e])) \subset C(\Xi,\ZZ)_{\ZZ^n}$.

\bigskip
\noindent
\section{Index theorem for the gap-labeling of the pinwheel tiling}
\medskip
\noindent
In a first step, we study the $K$-theory of the $C^*$-algebra associated to the dynamical system $C(\Omega) \rtimes \RR^2 \rtimes S^1$ and then we compute its image under the linear map $\tau^\mu_*$.\\
In the sequel, one identifies the class of an unbounded triple defined in \citeg{BaajJulg} with the class that it defines in $KK$-theory.\\
The reader can refer to \citeg{SkanKK} for definitions and properties of equivariant $KK$-theory groups of Kasparov that will be used in this section.

\subsection{Study of $K_0\big( C(\Omega) \rtimes \RR^2 \rtimes S^1)$} \label{studyK}

\text{ }\\
\noindent
To compute this $K$-theory, we will proceed in two steps. The first one consists in using the Dirac-Dual Dirac construction.\\
For this, let's consider the Dirac operator $\partial_2$ on $\RR^2$. Then, $F=\partial_2 (1+\partial_2^2)^{-\frac{1}{2}}$ is an elliptic pseudodifferential operator of order 0 on $H=L^2(\RR^2,\CC \oplus \, \CC)$.\\
Letting $C_0(\RR^2)$ act on $H$ by multiplication $f \mapsto M(f)$, one has a Kasparov cycle $(H,M,F)$.\\
Moreover, since $G=\RR^2 \rtimes S^1$ acts naturally on the left of $C_0(\RR^2)$ and of $H$, and since $F$ is then $G$-invariant, the class of $(H,M,F)$ defines an element $\alpha_G$ of $KK^G_0(C_0(\RR^2),\CC)$, called the \textit{fundamental element}.

\bigskip
\noindent
There exists an element $\sigma_G \in KK^G_0(\CC,C_0(\RR^2))$ such that the Kasparov product of $\alpha_G$ with $\sigma_G$ over $\CC$ is $\alpha_G \otimes_\CC \sigma_G = 1_{C_0(\RR^2)}$ i.e $\sigma_G$ is a right inverse for $\alpha_G$ (see \citeg{Kaspa1}, \citeg{Kaspa2}, \citeg{SkanKK}).\\
Furthermore, since $G=\RR^2 \rtimes S^1$ is amenable, $\sigma_G \otimes_{C_0(\RR^2)} \alpha_G = 1_{\CC}$.

\bigskip
\noindent
$\alpha_G$ is thus an invertible element in $KK^G_0(C_0(\RR^2),\CC)$ and we can prove :

\begin{prop} \label{KK1}
{\small$K_0 \left( C(\Omega) \rtimes \RR^2 \rtimes S^1 \right)$} is isomorphic to  {\small $K_0 \left( C(\Omega) \rtimes S^1 \right)$}.
\end{prop}

\bigskip
\noindent
\textbf{Proof :}
\begin{itemize}
\item[] Denote $[\partial_2]$ the class of the above cycle $(H,M,F)$.\\
We have $\tau_\Omega [\partial_2] \in KK^G_0 \left( C_0(\Omega \times \RR^2),C(\Omega) \right)$ where $\tau_\Omega [\partial_2] :=1_\Omega \otimes_{\CC} [\partial_2]$.\\
Using the descent homomorphism
$$J_G : KK^G_0 \left( C_0(\Omega \times \RR^2),C(\Omega) \right) \longrightarrow KK_0 \left( C_0(\Omega \times \RR^2) \rtimes G,C(\Omega)\rtimes G \right),$$ on this element, we obtain the element
$$J_G \left( \tau_\Omega \left( [\partial_2] \right) \right) \in KK_0 \left( C_0(\Omega \times \RR^2) \rtimes G,C(\Omega) \rtimes G \right).$$
Since the homomorphisms $\tau_\Omega$ and $J_G$ are compatible with Kasparov product and unit elements, $J_G \left( \tau_\Omega \left([\partial_2] \right) \right)$ is an invertible element and thus defines an isomorphism $\beta:=\otimes_{C_0(\Omega \times \RR^2) \rtimes G} \; J_G \left( \tau_\Omega \left ( [\partial_2] \right) \right)$ :
$$\beta : KK_0 \left( \CC,C_0(\Omega \times \RR^2) \rtimes G \right) \longrightarrow KK_0 \left( \CC,C(\Omega) \rtimes G \right).$$
Moreover, since $C(\Omega) \rtimes S^1$ is Morita equivalent to $C_0(\Omega \times \RR^2) \rtimes \RR^2 \rtimes S^1$ (see lemma \refg{morita}), we have an isomorphism 
$$\delta :  KK_0 \left( \CC,C(\Omega) \rtimes S^1 \right) \longrightarrow  KK_0 \left( \CC,C_0(\Omega \times \RR^2) \rtimes G \right).$$
The desired isomorphism of the proposition is then obtained as the composition $\beta o \delta$. 

\begin{flushright}
$\square$\\
\end{flushright}
\end{itemize}

\noindent
We recall a result that we used in the proof of the proposition and that we will often use in the sequel :

\begin{lem} \label{morita}
Let $G$ be a locally compact group and $X$ a locally compact space with a right $G$-action.\\
Let $\mathcal{K}\big( L^2(G) \big)$ be the compact operators on $L^2(G)$.\\
Then the $C^*$-algebras $\mathcal{K}\big( L^2(G) \big) \otimes C_0(X)$ and $C_0(X \times G) \rtimes G$ are isomorphic.\\
Moreover, for $G=\RR^2$, this isomorphism is $S^1$-equivariant and $C_0(X) \rtimes S^1$ is Morita equivalent to $C_0(X \times \RR^2) \rtimes \RR^2 \rtimes S^1$.
\end{lem}

\bigskip
\noindent
Thus, we proved that, to study $K_0 \big( C(\Omega) \rtimes \RR^2 \rtimes S^1 \big)$, it suffices to study $K_0(C(\Omega) \rtimes S^1)$.\\
To investigate this group, we use the 6 circles $F_1$, ... , $F_6$ stable under a rotation of angle $\pi$ around the origin and $F:= \bigcup F_i$, to obtain the following six term exact sequence :

$$\xymatrix{   K^{S^1}_0(C_0(\Omega \setminus F)) \ar[r]   & K^{S^1}_0(C(\Omega)) \ar[r] & K^{S^1}_0(C(F)) \ar[d]^{\partial} \\
                        K^{S^1}_1(C(F)) \ar[u] & K^{S^1}_1(C(\Omega)) \ar[l]  & K^{S^1}_1(C_0(\Omega \setminus F)) \ar[l]}$$

\noindent
But, $K^{S^1}_i(C(F)) \simeq K_i \big( C(F) \rtimes S^1 \big) \simeq K_i \Big( (\CC \oplus \CC)^6 \Big)$, the first isomorphism is proved in \citeg{Julg} and the second comes from the isomorphism of $C(F_i) \rtimes S^1$ with $\text{Ind}^{S^1}_{\{-1,1\}}(\CC) \rtimes S^1$ and from the fact that this $C^*$-algebra is Morita equivalent to $C^*(\ZZ/2\ZZ) \simeq \CC \oplus \CC$, see \citeg{Gra} (the $S^1$-action on $F_i$ is obtained as the $S^1$-action on $S^1$ defined by $e^{i\theta}.z:=e^{2i\theta}z$  and $\text{Ind}^{S^1}_{\{-1,1\}}(\CC) $ is the induced $C^*$-algebra defined as the space of continuous functions $f$ on $S^1$ with values in $\CC$ such that $f(-z)=-f(z)$).\\
Thus, $K^{S^1}_0(C(F)) \simeq \ZZ^{12}$ ($= \ZZ^6 \oplus \ZZ^6)$ and $K^{S^1}_1(C(F)) \simeq 0$.\\
Note that $K^{S^1}_0(C(F_i)) \simeq \ZZ^2$ is generated by two $S^1$-equivariant vector bundles. The first one is the trivial fiber bundle $F_i \times \CC$ with the diagonal action of $S^1$ where $S^1$ acts trivially on $\CC$ and the second one is the fiber bundle $F_i \times \CC$ with the diagonal action of $S^1$ where $S^1$ acts by multiplication on $\CC$.\\
Furthermore, $S^1$ acts freely and properly on $\Omega \setminus F$ since all the fixed points have been removed. Thus, 
$$K^{S^1}_i(C_0(\Omega \setminus F)) \simeq K_i(C_0(\Omega \setminus F) \rtimes S^1) \simeq K_i \left( C_0\big( (\Omega \setminus F)/S^1 \big) \right),$$ 
where the last isomorphism is obtained by Morita equivalence (see \citeg{Rief}).

\bigskip
\noindent
The six term exact sequence becomes :

{\small $$\xymatrix{   K_0( C_0(\Omega \setminus F)/S^1 )) \ar[rr]  & & K^{S^1}_0(C(\Omega)) \ar[rr] && \hspace{0.5cm} \ZZ^{12} \hspace{0.4cm} \ar[d]^{\partial} \\
                      \hspace*{0.5cm} 0 \hspace{0.5cm} \ar[u] && K^{S^1}_1(C(\Omega)) \ar[ll] && K_1( C_0( (\Omega \setminus F)/S^1)) \ar[ll]}$$}

\noindent To study the connecting map $\partial$, we need the following lemma.

\bigskip
\noindent
It links the $K$-theory group of a topological space of low dimension with its cohomology group (see \citeg{Mat}, section 3.4) : 

\begin{lem}  \label{chernentierMathey}
Let $X$ be a connected finite $CW$-complex of  dimension $\leqslant$ 3.
Then, there exist canonical  isomorphisms
$$ ch_{ev}^{\ZZ} := ch_0^{\ZZ} \oplus ch_2^{\ZZ} : K_0\big( C(X)  \big) \longrightarrow \check{H}^0(X;\ZZ) \oplus \check{H}^2(X;\ZZ)$$
$$ ch_{odd}^{\ZZ} := ch_1^{\ZZ} \oplus ch_3^{\ZZ} : K_1 \big( C(X) \big) \longrightarrow \check{H}^1(X;\ZZ) \oplus \check{H}^3(X;\ZZ)$$
that are natural for such complexes and compatible with the usual Chern character i.e such that the following diagram commutes
$$\xymatrix{ K_j \big( C(X) \big) \ar[d]_{ch_n^{\ZZ}} \ar[drr]^{ch_n} \\
\check{H}^n(X;\ZZ) \ar[rr] & & \check{H}^n(X;\QQ)}$$
where $j = n \text{ mod } 2$, $\check{H}^n(X;\ZZ) \rightarrow \check{H}^n(X;\QQ)$ is the canonical homomorphism induced by the inclusion $\ZZ \hookrightarrow \QQ$ and $ch_n$ is the component of degree $n$ of the Chern character.
\end{lem}

\bigskip
\noindent
Since $\Omega/S^1$ is the inverse limit of connected finite $CW$-complexes of dimension $2$ (see \citeg{OrmRadSad} or \citeg{Hai-coho}), we obtain the following lemma:
\bigskip
\noindent
\begin{lem} \label{coho}
The \v{C}ech cohomology groups with compact support and with integer coefficients $H_c^k((\Omega \setminus F)/S^1;\ZZ)$ vanish for $k\geqslant 3$.\\
Moreover, $$K_0 \left ( C_0 \Big((\Omega \setminus F)/S^1 \Big) \right) \simeq \check{H}^2_c \Big( (\Omega \setminus F)/S^1 \, ; \, \ZZ \Big)$$ and 
$$K_1 \left ( C_0 \Big((\Omega \setminus F)/S^1 \Big) \right) \simeq \check{H}^1_c \Big( (\Omega \setminus F)/S^1 \, ; \, \ZZ\Big).$$
\end{lem}

\noindent
\textbf{Proof :}

\begin{itemize}
\item[] In fact, according to \citeg{OrmRadSad} or \citeg{Hai-coho}, $\Omega/S^1$ is the inverse limit of $CW$-complexes of dimension $2$. Thus,  $\check{H}^k(\Omega/S^1 \, ; \, \ZZ)=0$ for $k \geqslant 3$.\\
The long exact sequence of relative cohomology groups associated to the pair $(\Omega/S^1, F/S^1)$ then gives :
{\small $$\hspace*{-0.5cm} \xymatrix{ \ldots \ar[r] & \check{H}^2(\Omega/S^1,F/S^1;\ZZ) \ar[r] & \check{H}^2(\Omega/S^1;\ZZ) \ar[r] & \check{H}^2(F/S^1;\ZZ) \\ 
\ar[r]  & \check{H}^3(\Omega/S^1,F/S^1;\ZZ) \ar[r] & \check{H}^3(\Omega/S^1;\ZZ) \ar[r] & \check{H}^3(F/S^1;\ZZ) \ar[r] & \ldots\\}$$}

Since $F/S^1$ is composed by $6$ points, $\check{H}^k(F/S^1;\ZZ) = 0$ for $k \geqslant 1$ and thus, $\check{H}^{k+1}(\Omega/S^1,F/S^1;\ZZ) \simeq \check{H}^{k+1}(\Omega/S^1;\ZZ)$ for $k \geqslant 1$.\\
From the result reminded above on cohomology groups of $\Omega/S^1$, we proved that the relative cohomology groups associated to $(\Omega/S^1, F/S^1)$ with integer coefficients vanish for degrees greater than 3.\\
To conclude, we use lemma 11 p.321 from \citeg{Span} to state that, for any $k$, $\check{H}^k( \Omega/S^1,F/S^1;\ZZ) \simeq \check{H}^k_c \Big( \big( \Omega/S^1 \big) \setminus \big( F/S^1\big) ; \ZZ  \Big)$.\\
Moreover, since $F$ is stable for the $S^1$-action, $(\Omega/S^1 \setminus F/S^1) = (\Omega \setminus F)/S^1$, which completes the proof of the first point of the lemma.

\bigskip
\noindent
Furthermore, we have :
$$K_0 \left ( C_0 \Big((\Omega \setminus F)/S^1 \Big) \right) = \tilde{K}_0 \left ( C \Big( \big( (\Omega \setminus F)/S^1 \big)^+ \Big) \right)$$
and
$$K_1 \left ( C_0 \Big((\Omega \setminus F)/S^1 \Big) \right) = K_1 \left ( C \Big( \big( (\Omega \setminus F)/S^1 \big)^+ \Big) \right)$$
where $\big( (\Omega \setminus F)/S^1 \big)^+$ is the Alexandroff compactification of $(\Omega \setminus F)/S^1$ and $\tilde{K}_0$ is the reduced $K$-theory.\\
Using results from \citeg{Hai-coho}, we can easily prove that $\big( (\Omega \setminus F)/S^1 \big)^+$ is the inverse limit of $CW$-complexes of dimension 2.\\
Thus, applying the results of lemma \refg{chernentierMathey}, we have proved the second point of the lemma.
\end{itemize}

\begin{flushright}
$\square$\\
\end{flushright}
We can then compute the kernel of $\partial$ :
\begin{lem}  \label{noyaubord}
$ \text{ Ker }\partial= \bigoplus \limits_{i=1}^7 \ZZ.q_i$ where
$$q_1=(1,1,1,1,1,1,0,0,0,0,0,0)$$
$$q_2=(0,1,1,1,1,1,1,0,0,0,0,0)$$
$$q_3=(1,0,1,1,1,1,0,1,0,0,0,0)$$
$$q_4=(1,1,0,1,1,1,0,0,1,0,0,0)$$
$$q_5=(1,1,1,0,1,1,0,0,0,1,0,0)$$
$$q_6=(1,1,1,1,0,1,0,0,0,0,1,0)$$
$$q_7=(1,1,1,1,1,0,0,0,0,0,0,1).$$
\end{lem}

\noindent
The proof is given in Appendix \refg{appnoyaubord}.

\bigskip
\noindent
Since $q_1$ lifts on the constant projection equal to $1$ on $\Omega$, we thus have proved :

\begin{prop} \label{somdirecte} We have
$$K_0(C(\Omega) \rtimes \RR^2 \rtimes S^1) \simeq K_0(C(\Omega) \rtimes S^1)  \simeq \ZZ \oplus \left ( \bigoplus \limits_{i=2}^7 \ZZ.q_i \right) \oplus \, \check{H}^2_c \Big( (\Omega \setminus F)/S^1 \, ; \, \ZZ \Big),$$
where $\ZZ$ is generated by the constant projection equal to $1$.
\end{prop}

\bigskip
\noindent
\subsection{Computation of the image under $\tau^\mu_*$ of the summand $\mathbb{Z}$} \label{copieZ}
\medskip
\noindent
In this section, we compute the image under $\tau^\mu_*$ of $\beta \circ \delta(\ZZ)$ in $K_0(C(\Omega)\rtimes \RR^2 \rtimes S^1)$ where $\beta \circ \delta$ is the isomorphism contructed in the proof of proposition \refg{KK1}.

\bigskip
\noindent
To compute it, we will consider the maps $\phi:\CC \rightarrow C(\Omega)$ and $\phi_A : A \rightarrow C(\Omega) \otimes A$ ($A$ a $C^*$-algebra with a $\RR^2 \rtimes S^1$-action) given by $\phi(z)(\omega):=z$ for any $\omega \in \Omega$ and $\phi_A:= \phi \otimes Id$. We will also use the induced map $\phi^{\RR^2}_A:A \rtimes \RR^2 \rightarrow \big( C(\Omega) \otimes A \big) \rtimes \RR^2$. The maps $\phi$ and $\phi_A$ are $\RR^2 \rtimes S^1$-equivariant and thus define, by functoriality of $KK$-theory, the following maps :
$$\phi_*:KK^{S^1}_0(\CC,\CC) \longrightarrow KK^{S^1}_0(\CC,C(\Omega)),$$
$$(\phi_A)_*:KK^{S^1}_0(\CC,A) \longrightarrow KK^{S^1}_0(\CC,C(\Omega) \otimes A),$$
$$(\phi^{\RR^2}_A)_*:KK^{S^1}_0(\CC,A \rtimes \RR^2) \longrightarrow KK^{S^1}_0 \Big ( \CC,\big( C(\Omega) \otimes A \big) \rtimes \RR^2 \Big),$$
$$\phi_A^*:KK^{S^1}_0(C(\Omega) \otimes A,C(\Omega)) \longrightarrow KK^{S^1}_0(A,C(\Omega)).$$ 
We will denote $\phi_*$ and $\phi^*$ these homomorphisms when no confusion would be possible or if $A=\CC$.

\bigskip
\begin{prop} \label{diagcomm}
We have the following commutative diagram :
$$\xymatrix{  K^{S^1}_0(\CC) \ar[rr]^{\phi_*} \ar[d]^{\eta_1} \ar@/_5pc/[dd]_{\delta_1} && K^{S^1}_0(C(\Omega)) \ar[d]^\eta \ar@/^6pc/[dd]^\delta\\
  K^{S^1}_0 \big( \mathcal{K}(L^2(\RR^2)) \big) \ar[rr]^{ \big( \phi_\mathcal{K} \big)_*} \ar[d]^{\Psi_*} && K^{S^1}_0 \big ( C(\Omega) \otimes \mathcal{K}(L^2(\RR^2)) \big) \ar[d]^{\Psi_*} \\
 K^{S^1}_0 \big( C_0(\RR^2) \rtimes \RR^2 \big) \ar[rr]^{\hspace{-0.3cm} \big( \phi^{\RR^2}_{C_0(\RR^2)} \big)_*} \ar[d]^{\beta_1} && K^{S^1}_0 \big( C_0(\Omega \times \RR^2) \rtimes \RR^2 \big) \ar[d]^{\beta}\\
 K^{S^1}_0 \big( C^*(\RR^2) \big) \ar[rr]^{\big( \phi^{\RR^2}_{\CC} \big)_*} && K^{S^1}_0 \big( C(\Omega) \rtimes \RR^2 \big) }$$
where the isomorphisms $\delta_1$ and $\beta_1$ are constructed in a similar way as the maps used in the proof of proposition \refg{KK1}, taking $\CC$ instead $C(\Omega)$, and by using the partial descent homomorphism $J^{S^1}_{\RR^2}$ defined in \citeg{Chab}. Moreover, $\Psi_*$ is the isomorphism induced by the $S^1$-equivariant isomorphism of lemma \ref{morita}, between $C(X) \otimes \mathcal{K} \big( L^2(\RR^2) \big)$ and $C_0(X \times \RR^2) \rtimes \RR^2$, with $X=\Omega$ and $X=$ one point. $\eta$ and $\eta_1$ are the isomorphisms given by Morita equivalence.
\end{prop}

\noindent
\textbf{Proof :}

\begin{itemize}
\item[] The proof is essentially based on naturality properties of the Kasparov product and the descent homomorphism.
\item The middle diagram commutes trivially.
\medskip
\item The isomorphism $\eta_1^{-1}:K^{S^1}_0 \Big( \mathcal{K} \big ( L^2(\RR^2) \big) \Big) \rightarrow K^{S^1}_0(\CC)$ is the Kasparov product by the class of the $S^1$-equivariant bimodule giving the Morita equivalence between $\mathcal{K} \big( L^2(\RR^2) \big)$ and $\CC$ i.e 
$$\eta_1^{-1}=\otimes_{ \mathcal{K} ( L^2(\RR^2) )} \Big [ \big( L^2(\RR^2),i,0 \big) \Big ].$$ 
We denote $y_1$ the class of $\big( L^2(\RR^2), \, i, \, 0 \big)$ in $KK^{S^1}_0\Big (\mathcal{K} \big( L^2(\RR^2) \big),\CC \Big)$.\\
Similarly, $\eta^{-1}:K^{S^1}_0 \Big ( \mathcal{K} \big( L^2(\RR^2) \big) \otimes C(\Omega) \Big) \rightarrow K^{S^1}_0(C(\Omega))$ is defined by $\otimes_{\mathcal{K} ( L^2(\RR^2) ) \otimes C(\Omega)} y_2$ where
$$y_2:=\Big[ \big(L^2(\RR^2) \otimes C(\Omega), \, \theta, \, 0 \big) \Big] \in KK^{S^1}_0 \Big( \mathcal{K} \big( L^2(\RR^2) \big) \otimes C(\Omega),C(\Omega) \Big)$$
where $\theta$ is the action of $\mathcal{K} \big( L^2(\RR^2) \big) \otimes C(\Omega)$ on the $S^1$-\'equivariant bimodule $L^2(\RR^2) \otimes C(\Omega)$ (this is the natural action of compact operators on $L^2(\RR^2)$ and multiplication operator on $C(\Omega)$).\\
We then have 2 commutative diagrams  (see \citeg{SkanKK}) :
$$\xymatrix{ K^{S^1}_0(\CC) \ar[rrr]^{\phi_*} &&& K_0^{S^1}(C(\Omega)) \\
&&&\\
  K^{S^1}_0 \big( \mathcal{K} \big( L^2(\RR^2) \big) \big) \ar[uu]^{\bigotimes_{\mathcal{K} ( L^2(\RR^2) )} \; y_1} \ar[rrruu]_{\hspace{0.5cm}\bigotimes_{\mathcal{K} ( L^2(\RR^2) )} \; \phi_*(y_1)}  }$$

$$\xymatrix{  &&& K_0^{S^1}(C(\Omega)) \\
&&&\\
  K^{S^1}_0 \Big( \mathcal{K} \big(L^2(\RR^2) \big) \Big)  \ar[rrruu]^{\hspace{-1cm} \bigotimes_{\mathcal{K} ( L^2(\RR^2)  )}  \phi_{\mathcal{K}}^*(y_2)}  \; \ar[rrr]_{\hspace{-0.3cm} \big( \phi_\mathcal{K} \big)_*}  &&&  K^{S^1}_0 \Big( \mathcal{K} \big( L^2(\RR^2) \big) \otimes C(\Omega) \Big)  \ar[uu]_{\bigotimes_{\mathcal{K} ( L^2(\RR^2) ) \otimes C(\Omega)} \; y_2} }$$
But, $$\phi_\mathcal{K}^*(y_2) = \big[ (L^2(\RR^2) \otimes C(\Omega), \, \theta \circ \phi_\mathcal{K}, \, 0) \big]$$
and 
$$\phi_*(y_1) = \big[ (L^2(\RR^2) \otimes C(\Omega), \, i \otimes 1, \, 0) \big].$$
Since $\theta \circ \phi_\mathcal{K}$ acts as $i \otimes 1$ on $L^2(\RR^2) \otimes C(\Omega)$, we have $\phi_\mathcal{K}^*(y_2) = \phi_*(y_1)$ and thus, the first upper diagram in the proposition is commutative (the vertical arrows are isomorphisms).
\medskip
\item As above, we get the following commutative diagrams :
{\small $$\xymatrix{ K^{S^1}_0 \big( C_0(\RR^2) \rtimes \RR^2 \big) \ar[rrr]^{\big( \phi_{C_0(\RR^2)}^{\RR^2} \big)_*} \ar[rrrdd]_{\hspace{-2.5cm} \bigotimes \big( \phi_{C_0(\RR^2)}^{\RR^2} \big)^* \Big(J_{\RR^2}^{S^1} \big(\tau_\Omega([\partial_2]) \big) \Big)} &&& K_0^{S^1} \big ( C_0(\Omega \times \RR^2) \rtimes  \RR^2 \big) \ar[dd]^{\bigotimes J_{\RR^2}^{S^1} \big( \tau_\Omega([\partial_2]) \big)} \\
&&&\\
  &&& K^{S^1}_0 \big( C(\Omega) \rtimes \RR^2 \big)  }$$}

{\small $$\xymatrix{  K^{S^1}_0 \big( C_0(\RR^2) \rtimes \RR^2 \big) \ar[dd]_{\hspace{-1cm} \bigotimes J_{\RR^2}^{S^1}([\partial_2])}  \ar[rrrdd]^{\hspace{0.5cm} \bigotimes \big( \phi_{\CC}^{\RR^2} \big)_* \big( J_{\RR^2}^{S^1}([\partial_2]) \big)} &&&\\
&&&\\
  K^{S^1}_0 \big( C^*(\RR^2) \big)  \ar[rrr]_{\hspace{-0.3cm} \big( \phi_{\CC}^{\RR^2} \big)_*}  &&&  K^{S^1}_0 \big( C(\Omega) \rtimes \RR^2 \big)  }$$}

\noindent
But, ${\displaystyle \big( \phi_{C_0(\RR^2)}^{\RR^2} \big)^* \Big(J_{\RR^2}^{S^1} \big(\tau_\Omega([\partial_2]) \big) \Big) = J_{\RR^2}^{S^1} \Big( \phi'^* \big(\tau_\Omega([\partial_2]) \big) \Big)}$ where $\phi' = \phi_{C_0(\RR^2)}$ is $\RR^2$-equivariant.\\
Similarly, $\big( \phi_{\CC}^{\RR^2} \big)_* \big( J_{\RR^2}^{S^1}([\partial_2]) \big) = J_{\RR^2}^{S^1} \big( \phi_*([\partial_2]) \big)$.\\
We have $$\phi_*([\partial_2]) =\Big [ \big( H \otimes C(\Omega),M \otimes 1, F \otimes 1 \big) \Big] \in KK^{\RR^2 \rtimes S^1}_0 \big( C_0(\RR^2),C(\Omega) \big)$$ with the notation of the class $[\partial_2]$ and 
{\small $$\phi'^* \big(\tau_\Omega([\partial_2]) \big) = \Big[ \big ( H \otimes C(\Omega), (M \otimes i_{C(\Omega)}) \circ \phi' , F \otimes 1 \big) \Big] \in KK^G_0 \big( C_0(\RR^2),C(\Omega) \big)$$}
with $i_{C(\Omega)}:C(\Omega) \rightarrow \mathcal{L}(C(\Omega))$ the multiplication operator on $C(\Omega)$.\\
Since $(M \otimes i_{C(\Omega)}) \circ \phi'$ is the representation of $C_0(\RR^2)$ on $H \otimes C(\Omega)$ obtained by multiplication of $C_0(\RR^2)$ on $H$ and $M \otimes 1$ is the same representation, we thus have equality of these two classes and of the two diagonal homomorphisms in the above diagrams.\\
The bottom diagram of the proposition is thus commutative.
\begin{flushright}
$\square$
\end{flushright}
\end{itemize}

\noindent Let's consider the forgetful homomorphism $r^*:KK^{S^1}(A,B) \longrightarrow KK(A,B)$ for $A$ and $B$ two $C^*$-algebras endowed with $S^1$-actions.\\
This homomorphism commutes trivially with the homomorphism $\eta_1$ obtained by Morita equivalence (equivariant or not) :
$$\xymatrix{ K_0(\CC) \ar[d]_{\eta_1} & K_0^{S^1}(\CC) \ar[l]_{r^*} \ar[d]^{\eta_1} \\
K_0(\mathcal{K}) & K^{S^1}_0(\mathcal{K}) \ar[l]^{r^*} }$$
It also commutes trivially with the homomorphism $\Psi_*$ of proposition \refg{diagcomm}.\\
Finally, the definition of the partial descent homomorphism in \citeg{Chab} and the one of the usual descent homomorphism gives the equality $r^* \circ J^{S^1}_{\RR^2} = J_{\RR^2} \circ r'^*$ where  $r'^*:KK^{\RR^2 \rtimes S^1}(A,B) \longrightarrow KK^{\RR^2}(A,B)$ is the restriction homomorphism of $\RR^2 \rtimes S^1$ on $\RR^2$ (see \citeg{SkanKK}).\\
In conclusion, we obtain the following commutative diagram :
$$\xymatrix{ & K_0^{S^1}(\CC) \ar[ld]_{r^*} \ar[dd]^{\beta_1 \circ \delta_1}_\simeq \ar[r]^{ \hspace{-0.3cm}\phi_*} & K^{S^1}_0(C(\Omega)) \ar[dd]^{\beta \circ \delta}_\simeq \\
K_0(\CC) \ar[dd]_{\beta' \circ \delta'}^\simeq & &\\
& K^{S^1}_0(C^*(\RR^2)) \ar[ld]^{r^*} \ar[r]^{\hspace{-0.2cm} \phi^{\RR^2}_*} & K^{S^1}_0(C(\Omega) \rtimes \RR^2) \ar[d]^{r^*}\\
K_0(C^*(\RR^2)) \ar[rr]_{\phi^{\RR^2}_*} & & K_0(C(\Omega) \rtimes \RR^2) \\
K_0(C_0(\RR^2))  \ar[u]^{F_*}_\simeq&&
}$$
where $\beta'$ and $\delta'$ are obtained in a similar way as $\beta_1$ and $\delta_1$ but without the $S^1$-equivariance.\\
Thus, the projection generating the summand $\ZZ$ in $K_0^{S^1}(C(\Omega))$ comes from the constant projection equal to $1$ in $K^{S^1}_0(\CC)$ that is sent on the projection $1$ in $K_0(\CC)$.\\
Thus, $r^* \circ \beta \circ \delta \circ \phi_* ([1])= \phi^{\RR^2}_* \circ \beta' \circ \delta' \circ r^*([1]) = \phi^{\RR^2}_* \circ F_*(\pm \beta ott) $ where $\beta ott$ is the element $[Bott] - 1$ of $K_0(C_0(\RR^2))$ obtained from the Bott projection and $F$ is the Fourier transform giving the isomorphism $C^*(\RR^2) \simeq C_0(\RR^2)$ $\Big( r^*([1]) = [1]$ spans $K_0(\CC)$ and since $\beta' \circ \delta'$ is an isomorphism, the image of $[1]$ by this homomorphism has to be $\pm$ a generator of $K_0(C_0(\RR^2)) \Big)$.

\bigskip
\noindent
It then suffices to follow how $\tau^\mu_*$ is changed under all these homomorphisms.\\
To begin, the map $\tau^\mu_*$ on $K_0(C(\Omega) \rtimes \RR^2 \rtimes S^1))$ becomes $\tau'^\mu_*$ on $K^{S^1}_0(C(\Omega) \rtimes \RR^2)$ and on $K_0(C(\Omega)\rtimes \RR^2)$ where $\tau'^\mu$ is given on the dense subalgebra $C_c(\Omega \times \RR^2)$ by (see \citeg{Julg})
$$\tau'^\mu (f) = \int_\Omega f(\omega,0) d\mu(\omega).$$
Moreover, $\tau'^\mu_* \circ r^* = \tau'^\mu_*$.\\
$\tau'^\mu_*$ then reads on $K_0(C^*(\RR^2))$ as $\tau_*$ where $\tau$ is the trace on $C^*(\RR^2)$ defined by $\tau(f):=f(0)$ for any $f \in C_c(\RR^2)$.\\
We then have $\tau \circ F (f) = F(f)(0) = \int_{\RR^2} f(x) dx$ and thus, 
$$\tau_* \circ F_*(\beta ott) = \int_{\RR^2} \Big( Tr(Bott(x)) -1 \Big) \, dx =0.$$
Thereby, 
$$\tau'^\mu \Big( \beta \circ \delta \circ \phi_* ([1]) \Big ) = \tau'^\mu \Big( r^* \circ \beta \circ \delta \circ \phi_* ([1]) \Big)$$
and so 
$$\tau'^\mu \Big( \beta \circ \delta \circ \phi_* ([1]) \Big) = \tau'^\mu \Big( \phi^{\RR^2}_* \circ F_*(\pm \beta ott) \Big) = \tau_* \circ F_*(\pm \beta ott) =0.$$
We just proved :

\begin{theo} The summand $\ZZ.q_1$ in $K_0 \Big( C(\Omega) \rtimes \RR^2 \rtimes S^1 \Big)$ is traceless. \label{Z1}
\end{theo}

\bigskip
\noindent
\subsection{Study of the summand {\small $\check{H}^2_c \left(  (\Omega \setminus F)/S^1 \, ; \, \ZZ \right)$} in {\small $K_0 \big( C(\Omega) \rtimes \RR^2 \rtimes S^1 \big)$}}
\medskip
\noindent
We will now focus on the image of the summand $\check{H}^2_c \Big( (\Omega \setminus F)/S^1 \, ; \, \ZZ\Big)$ under $\tau^\mu_*$.

\bigskip
\noindent
For this, we construct a map $K_1 \big( C(\Omega) \big) \longrightarrow K_0 \big ( C(\Omega) \rtimes S^1 \big)$ which is onto on $\check{H}^2_c \Big( (\Omega \setminus F)/S^1 \, ; \, \ZZ\Big)$. This map, denoted $\beta_{S^1} \circ \delta_{S^1}$, is defined in a similar way as the one in proposition \refg{KK1}.\\
We consider the Dirac operator $\partial_1$ which is $S^1$-equivariant on the circle $S^1$. It defines a class $[\partial_1] \in KK_1^{S^1} \big ( C(S^1), \CC)$.\\
The desired map is given by the composition of the isomorphism $\delta_{S^1}$ (coming from the Morita equivalence of $C(\Omega)$ and $C(\Omega \times S^1) \rtimes S^1$) and the Kasparov product $\beta_{S^1}:=\bigotimes J_{S^1} \big( \tau_\Omega [\partial_1] \big)$.\\
The same reasoning as the one used at the beginning of lemma \refg{coho} allow us to state :

\begin{lem} \label{cohonulle}
$H^k(\Omega)=0$ for every $k \geqslant 4$.\\
Thus $K_0(C(\Omega)) \simeq \check{H}^0(\Omega ; \ZZ) \oplus \, \check{H}^2(\Omega ; \ZZ)$ and $K_1(C(\Omega)) \simeq \check{H}^1(\Omega ; \ZZ) \oplus \,\check{H}^3(\Omega ; \ZZ)$.
\end{lem}
\medskip
\noindent
We then prove that the $\check{H}^3 \big( \Omega \, ; \, \ZZ \big)$ summand of $K_1 \big( C(\Omega) \big)$ is isomorphic to the $\check{H}^2_c \Big( (\Omega \setminus F)/S^1 \, ; \, \ZZ\Big)$ summand of $K_0 \big ( C(\Omega) \rtimes S^1 \big)$ and then we show that this isomorphism can be read in $K$-theory as the map constructed above.

\bigskip
\noindent
First, we have the following commutative diagram :

{\small $$\xymatrix{ &  & K_1 \big( C_0(\Omega \setminus F) \big) \ar[r]^{i_*}  \ar[d]^{\beta'_{S^1} \circ \delta'_{S^1}} \ar[lld]_{\Psi} & K_1 \big( C(\Omega) \big) \ar[d]_{\beta_{S^1} \circ \delta_{S^1}} \\
K_0 \Big ( C_0 \big( ( \Omega \setminus F) / S^1  \big) \Big) &  & K_0 \big ( C_0(\Omega \setminus F) \rtimes S^1 \big)   \ar[ll]_{Morita \; eq.} \ar[r]_{i_*}  &  K_0 \big ( C(\Omega) \rtimes S^1 \big) }$$}

\noindent
where $\psi$ is the composition of $\beta'_{S^1} \circ \delta'_{S^1}$ and of the isomorphism induced by Morita equivalence of $C_0(\Omega \setminus F) \rtimes S^1$ and $C_0 \big( ( \Omega \setminus F) / S^1  \big)$.\\
Writing the long exact sequence of the relative cohomology groups of the pair $(\Omega, F)$, one can prove that $i_*:\check{H}^3(\Omega ; \ZZ) \rightarrow  \check{H}_c^3(\Omega \setminus F ; \ZZ)$ is an isomorphism since $\check{H}^2(F ; \ZZ)=\check{H}^3(F ; \ZZ)=0$.\\
It then suffices to show that the $\check{H}_c^3(\Omega \setminus F ; \ZZ)$ part of $K_1 \big( C_0(\Omega \setminus F) \big)$ is sent isomorphically on $\check{H}^2_c \Big( (\Omega \setminus F)/S^1 \, ; \, \ZZ \Big)$ and then to compare this isomorphism with $\beta'_{S^1} \circ \delta'_{S^1}$.\\
For this, we have the following proposition :

\begin{prop} (\citeg{Bred})   
The projection $\Omega \setminus F \longrightarrow (\Omega \setminus F)/S^1$ is a  $S^1$-principal bundle.
\end{prop} 

\noindent
Thanks to this proposition, to the Leray-Serre spectral sequence and to the resulting  Gysin sequence, we then get :

\begin{cor} \label{cohoreduc}
$$\begin{array}{ccl}
\check{H}^3_c \big( \Omega \setminus F \, ; \, \ZZ\big) & = & \check{H}^2_c \Big( \big( \Omega \setminus F \big) / S^1 \, ; \, \ZZ \Big) \otimes \check{H}^1(S^1 ; \ZZ) \\
& \simeq & \check{H}^2_c \Big( \big( \Omega \setminus F \big) / S^1 \, ; \, \ZZ\Big)
\end{array}$$
\end{cor}

\noindent
The isomorphism of this corollary is given by "integration along the fiber". It remains to see that this isomorphism is the same map as $\Psi$, under the Chern character.

\bigskip
\noindent
We want to prove that the following diagram is commutative :
$$\xymatrix{  \check{H}^3_c \big( \Omega \setminus F \, ; \, \ZZ \big) \ar@{^{(}->>}[r] & \check{H}^2_c \Big( \big( \Omega \setminus F \big) / S^1 \, ; \, \ZZ \Big) \\
K_1 \big( C_0(\Omega \setminus F) \big) \ar@{->>}[u]^{ch_3^{\ZZ}}  \ar[r]^{\hspace{-0.5cm}\psi}  & K_0 \Big ( C_0 \big( ( \Omega \setminus F) / S^1  \big) \Big) \ar@{^{(}->>}[u]_{ch^{\ZZ}}}$$

\noindent
But $\check{H}^3_c \big( \Omega \setminus F \, ; \, \ZZ \big) $ resp. $\check{H}^2_c \Big( \big( \Omega \setminus F \big) / S^1 \, ; \, \ZZ \Big)$ is spanned by $\check{H}^3_c(C \times D^2 \times S^1 ; \ZZ)$ resp. $\check{H}^2_c(C \times D^2 ; \ZZ)$ i.e by  $\check{H}^3_c(D^2 \times S^1 ; \ZZ)$ resp. $\check{H}^2_c(D^2 ; \ZZ)$ (where  $C$ is a Cantor set and $D^2$ is an open subset of $\RR^2$) since  $\Omega \setminus F$ resp. $\big( \Omega \setminus F \big) / S^1$ is covered by the closure of disjoint union of open sets of the form  $C \times D^2 \times S^1$ resp. $C \times D^2$ which intersections are of dimension $2$ resp. $1$.

\noindent
We then have the following diagram (the cohomology groups are with integer coefficients) :

\def \commutatif{\ar@{}[rd]|{ \circlearrowleft}}
{\scriptsize $$\hspace{0cm} \xymatrix{ \check{H}^3_c \big( \Omega \setminus F \big) \ar[rrr]^\simeq   & \commutatif & & \check{H}^2_c \Big( \big( \Omega \setminus F \big)/S^1 \Big) \\
\commutatif & \check{H}^3_c(D^2 \times S^1)  \ar[lu] \eq[r] & \check{H}^2_c(D^2) \ar[ru] \commutatif  & \\
& K_1 \big( C_0(D^2 \times S^1) \big) \ar@{->>}[u] \ar[r]_{\hspace{0.3cm}\psi'} \ar[ld] \commutatif & K_0 \big( C_0(D^2) \big) \eq[u] \ar[rd] &\\
K_1 \big(C_0(\Omega \setminus F) \big) \ar@{->>}[uuu] \ar[rrr]_{\psi} & & & K_0 \Big( C_0 \big( (\Omega \setminus F)/S^1 \big) \Big) \ar[uuu]_\simeq  }$$}

\noindent
where $\psi' :  K_1 \big( C_0(D^2 \times S^1) \big) \rightarrow  K_0 \big( C_0(D^2 \times S^1) \rtimes S^1 \big) \rightarrow K_0 \big( C_0(D^2) \big)$ is the composition of $\beta''_{S^1} \circ \delta''_{S^1}$ (obtained as above with $D^2 \times S^1$ instead of $\Omega$) and the isomorphism induced by Morita equivalence.\\
It is enough to show that the middle diagram is commutative to prove that the integration along the fiber at the level of cohomology groups is $\Psi$ at the level of $K$-groups.\\
Integration along the fiber sends a generator of $\check{H}^3_c(D^2 \times S^1 ; \ZZ)$ on a generator of $\check{H}^2_c(D^2 ; \ZZ )$, and so, we must show that $\psi'$ maps the Bott generator of $K_1 \big(C_0(D^2 \times S^1) \big)$ on the one of $K_0 \big(C_0(D^2) \big)$.\\
We can reduce the problem a little bit more thanks to the following diagram :
\begin{equation} \label{eqbott}
\xymatrix{K_1 \big( C(S^1) \big ) \ar[rr]^{\hspace{-0.5cm} \otimes_\CC \beta ott} \ar[d]_{\beta'' \circ \delta'' = \otimes_{C(S^1)} y'_1} & & K_1\big( C_0(D^2 \times S^1) \big) \ar[d]^{\beta' \circ \delta' = \otimes_{C_0(D^2 \times S^1)} y'_2} \\
K_0(\CC) \ar[rr]_{\hspace{-0.5cm}  \otimes_\CC Bott} & & K_0\big( C_0(D^2) \big)}
\end{equation}
the map $\otimes_\CC \beta ott$ ($ \beta ott \in KK_0(\CC,C_0(D^2))$) is given by Bott periodicity in $KK$-theory (\citeg{SkanKK} p.214).\\
We denote $\alpha \in KK_0(C_0(D^2),\CC)$ its inverse.\\
Let's $\varepsilon$ be the class, in $KK_0\big( C(S^1) , C(S^1 \times S^1) \rtimes S^1 \big)$, of the equivalence bimodule of the Morita equivalence of $C(S^1)$ with $C(S^1 \times S^1) \rtimes S^1$.\\
We then have  
$$y'_1 = \varepsilon \otimes_{C(S^1 \times S^1) \rtimes S^1} J_{S^1} \big( \tau_{C(S^1)} [\partial_1] \big)$$
and
$$y'_2=\tau_{C_0(D^2)} (\varepsilon) \otimes_{C_0(D^2) \otimes C(S^1 \times  S^1) \rtimes S^1} J_{S^1} \big( \tau_{C_0(D^2)\otimes C(S^1)} [\partial_1] \big).$$ 
Since $S^1$ is not acting on $D^2$ (the $S^1$-action on $\Omega \setminus F$ is given by action on $S^1$ in $C \times D^2 \times S^1$), we obtain 
$$ J_{S^1} \big( \tau_{C_0(D^2 \times S^1)} [\partial_1] \big) = \tau_{C_0(D^2)} \Big( J_{S^1} \big( \tau_{C(S^1)} [\partial_1] \big) \Big).$$ 
Thus (see \citeg{SkanKK}) :
$$y'_2=\tau_{C_0(D^2)} \Big( \varepsilon \otimes_{C(S^1 \times S^1) \rtimes S^1} J_{S^1} \big( \tau_{C(S^1)} [\partial_1] \big)\Big)= \tau_{C_0(D^2)} (y'_1)$$
The above diagram is thus commutative since, if $X \in K_1\big( C_0(D^2 \times S^1) \big)$, we have (thanks to commutativity of the Kasparov product over $\CC$) :
$$\begin{array}{rcl}
\big ( (X \otimes_{C_0(D^2)} \alpha ) \otimes_{C(S^1)} y'_1 \big) \otimes_{\CC} \beta ott &  =  & \beta ott \otimes_\CC \big ( (X \otimes_{C_0(D^2)}\alpha) \otimes_{C(S^1)} y'_1 \big) \\
& = & \big ( (\beta ott \otimes_\CC X) \otimes_{C_0(D^2)}\alpha \big ) \otimes_{C(S^1)} y'_1 \\
& = & \big( (X \otimes_\CC \beta ott) \otimes_{C_0(D^2)}\alpha \big ) \otimes_{C(S^1)} y'_1\\
& = & \big ( X \otimes_\CC (\beta ott \otimes_{C_0(D^2)} \alpha) \big ) \otimes_{C(S^1)} y'_1\\
& = & (X \otimes_\CC 1_\CC) \otimes_{C(S^1)} y'_1\\
& = & X \otimes_{C(S^1)} y'_1\\
& = & X \otimes_{C_0(D^2 \times S^1)} y'_2
\end{array}$$

\noindent
It then remains to prove that the generator $[u]$ of $K_1(C(S^1))$ is sent on the generator $[1]$ of $K_0(\CC)$.\\
But the map $\beta'' \circ \delta'' : K_1(C(S^1)) \longrightarrow K_0(\CC)$, in the diagram \textbf{(}\refg{eqbott}\textbf{)}, is given by $\otimes_{C(S^1)} [\partial_1]$ i.e the Kasparov product by the class of the Dirac operator on the circle $[\partial_1] \in KK_1(C(S^1),\CC)$ or, in other words, by the odd index of the Dirac operator twisted by a unitary in $K_1(C(S^1))$.\\
This map sends $[u]$ on $[1]$.\\
Thus, we just have proved :

\medskip
\noindent
\begin{theo} The map $K_1(C(\Omega)) \longrightarrow K_0 \big ( C(\Omega) \rtimes \RR^2 \rtimes S^1 \big)$ is surjective from $\check{H}^3(\Omega ; \ZZ)$ onto $\check{H}^2 \Big( (\Omega \setminus F)/S^1) \, ; \, \ZZ\Big)$.
\end{theo}

\bigskip
\noindent
\subsection{Index theorem to compute the trace}
\bigskip
\noindent
We now want to identify the map $\beta \circ \delta \circ \beta_{S^1} \circ \delta_{S^1}$ constructed from $K_1 \left( C(\Omega) \right)$ to $K_0 \left( C(\Omega)\rtimes \RR^2 \rtimes S^1 \right)$ with the Kasparov product by the class of the unbounded Kasparov cycle (see \citeg{BaajJulg} for definitions) 
$$\tilde{D_3}:=(C(\Omega) \rtimes \RR^2 \rtimes S^1 \oplus C(\Omega) \rtimes \RR^2 \rtimes S^1 ,M_{C(\Omega)},D_3)$$
where $D_3$ is the Dirac operator along the leaves of $\Omega$ and $M_{C(\Omega)}$ is the natural representation of $C(\Omega)$ on $C(\Omega) \rtimes \RR^2 \rtimes S^1 \oplus C(\Omega) \rtimes \RR^2 \rtimes S^1$.

\bigskip
\noindent
Using the measured index theorem for foliation due to Alain Connes \citeg{Con1} extended to foliated spaces in \citeg{MooSch}, we will have (see also  \citeg{DouHurKam}) :

\begin{theo} \hspace{0cm}\\ \label{indice}
$\forall \; b \in \text{Im}(\beta \circ \delta \circ \beta_{S^1} \circ \delta_{S^1})$, $\tau^\mu_*(b)= \Big \langle Ch_\tau([u]), [C_{\mu^t}] \Big \rangle$ where $C_{\mu^t}$ is the Ruelle-Sullivan current associated to $\mu^t$, $Ch_\tau$ is the tangential Chern character (see \citeg{MooSch}) and $[u] \in K_1(C(\Omega))$ is a lift of $b$ for $\beta \circ \delta \circ \beta_{S^1} \circ \delta_{S^1}$.
\end{theo} 

\bigskip
\noindent
$C_{\mu^t}$ defines a class in the tangential cohomology group $H_3^\tau(\Omega)$ and $Ch_\tau([u])$ is an element of the longitudinal cohomology group $H_\tau^3(\Omega)$ (see \citeg{MooSch}).\\
The representative of a class in $H^3_\tau(\Omega)$ locally looks like 
$$\sigma=a(x_1,x_2,\theta,\omega) dx_1 \wedge dx_2 \wedge d\theta.$$
The Ruelle-Sullivan current is then locally defined as the integral :
$$[C_{\mu^t}] \big( [\sigma] \big) = \int_{\Xi} \int_{\RR^2 \times S^1} a(x_1,x_2,\theta,\omega) dx_1  dx_2 d\theta d\mu^t(\omega).$$
To obtain theorem \ref{indice}, we proceed in several steps.\\
First, we identify $\beta_{S^1} \circ \delta_{S^1}$ with the Kasparov product by the class of the element in $KK_1(C(\Omega),C(\Omega)\rtimes S^1)$ induced by the Dirac operator $d_1$ of dimension $1$ along the leaves of $\Omega$ foliated by $S^1$.\\
Next, we identify $\beta \circ \delta$ with the Kasparov product by the class of the element in $KK_0(C(\Omega) \rtimes S^1,C(\Omega) \rtimes \RR^2 \rtimes S^1)$ induced by the tangential Dirac operator $D_2$ of dimension $2$ transverse to the inclusion of the foliation of $\Omega$ by $S^1$ in the foliation by $\RR^2 \rtimes  S^1$ (see \citeg{HilSkan}).\\
Finally, we prove that the composition of these two homomorphisms is given by $\otimes_{C(\Omega)} \big [\tilde{D_3} \big]$.

\bigskip
\noindent
For these computations, we will use the characterisation of the product of unbounded cycles made by Kucerovsky in \citeg{Kuc} (see also \citeg{BaajJulg} for the definitions).

\begin{defi}
Let $A$ and $B$ be two graded $C^*$-algebras.\\
An \textbf{unbounded Kasparov module} is a triple $(\mathcal{H},\phi,D)$ where $\HH$ is a graded Hilbert $B$-module, $\phi:A \rightarrow \mathcal{L}(\mathcal{H})$ is a $*$-homomorphism, and $D$ is a densely defined, selfadjoint and regular unbounded operator of order 1, such that :
\begin{enumerate}
	\item $\phi(a)(1+D^2)^{-1} \in \mathcal{K}(\HH)$,
	\item for all $a$ in some dense subalgebra of $A$, the domain of $D$ is stable under $\phi(a)$ and $[D,\phi(a)] = D\phi(a) - (-1)^{\partial a}\phi(a)D$ is defined on $\mathfrak{Dom} \, D$  and can be extended to an element of $\mathcal{L}(\HH)$.
\end{enumerate}
We denote $\Psi(A,B)$ the set of unbounded Kasparov modules.\\
If $\HH$ is ungraded, such a module is an \textbf{ungraded unbounded Kasparov module} and we denote $\psi^1(A,B)$ the set of such modules.
\end{defi}

\bigskip
\noindent
Julg and Baaj have proved the following result (see \citeg{BaajJulg} and \citeg{Kuc}):

\begin{prop}
If $(\HH,\phi,D) \in \Psi(A,B)$ then $(\HH,\phi,D(1+D^2)^{-1/2})$ defines a Kasparov cycle in $KK_0(A,B)$.\\
Similarly, if $(\HH,\phi,D) \in \Psi^1(A,B)$ then $(\HH,\phi,D(1+D^2)^{-1/2})$ defines a Kasparov cycle in $KK_1(A,B)$.
\end{prop}

\medskip
\noindent
There is then a map {\small $\beta^{BJ}:\Psi(A,B) \rightarrow KK_0(A,B)$} (resp. {\small $\Psi^1(A,B) \rightarrow KK_1(A,B)$}) defined by $\beta^{BJ}(\HH,\phi,D)= \big[ (\HH,\phi,D(1+D^2)^{-1/2}) \big]$.

\bigskip
\noindent
Kucerovsky then proved a criterion to compute the Kasparov product of two unbounded modules:
\begin{theo} \label{Kuce}
Fix $(E_1 \hat{\otimes}_{\phi_2} E_2, \phi_1 \hat{\otimes} 1, D) \in \Psi(A,C)$, $(E_1, \phi_1, D_1) \in \Psi(A,B)$ and $(E_2, \phi_2, D_2) \in \Psi(B,C)$ such that :
	\begin{enumerate}[(i)]
		\item for any $x$ in a dense subset of $\phi_1(A)E_1$,
				$$\left[ \left( 
					\begin{array}{cc}
						D & 0 \\
						0 & D_2
					\end{array} \right), \left(
					\begin{array}{cc}
						0 & T_x \\
						T_x^* & 0
					\end{array} \right) \right]$$
				is bounded on $\mathfrak{Dom} \, (D \oplus \, D_2)$ (where $\mathfrak{Dom} \, F$ is the domain of the unbounded operator $F$);
		\item $\mathfrak{Dom} \, D \subset \mathfrak{Dom} \, (D_1 \hat{\otimes} 1)$ (or vice versa);
		\item and $\big \langle (D_1 \hti 1 )x,Dx \big \rangle + \big \langle Dx, (D_1 \hti 1) x \big \rangle \geqslant \kappa \langle x , x \rangle$ for any $x$ in the domain;
	\end{enumerate}
where $x \in E_1$ is homogeneous and $T_x:E_2 \rightarrow E$ maps $e \mapsto x \hat{\otimes} e$.\\
Then $(E_1 \hat{\otimes}_{\phi_2} E_2, \phi_1 \hat{\otimes} 1, D)$ represents the Kasparov product of $(E_1, \phi_1, D_1)$ and $(E_2, \phi_2, D_2)$.
\end{theo}

\bigskip
\noindent
This theorem applies on unbounded Kasparov modules of $\Psi(A,B)$.\\
In the sequel, we use ungraded unbounded Kasparov module of $\Psi^1(A,B)$.\\
To apply the theorem of Kucerovsky, we thus need to use the isomorphism $KK_1(A,B) \simeq KK_0(A,B \otimes \CC_1)$ (where $\CC_1$ is the complex Clifford algebra of dimension $2$ spanned by an element $\alpha$ of degree $1$ satisfying $\alpha^*=\alpha$ and $\alpha^2=1$) given by : 
$$ \begin{array}{ccc}
KK_1(A,B) & \longrightarrow &  KK_0(A,B \otimes \CC_1)\\
\text{} \big [ (\HH,\phi,F) \big] & \longmapsto & \big[ (\HH \hat{\otimes} \CC_1 , \phi \hat{\otimes} Id, F \hat{\otimes} \alpha) \big]
\end{array}$$

\noindent
We also have : if $(\HH,\phi,D) \in \Psi^1(A,B)$ then {\small $\big( \HH \hat{\otimes} \CC_1,\phi \hat{\otimes} Id, D \hat{\otimes} \alpha \big) \in \Psi(A,B \otimes \CC_1)$}.

\bigskip
\noindent
To compute the Kasparov product of an element in $KK_1$ and an element in $KK_0$ or conversely, we use the two following lemmas (the proof can be found in \citeg{Hai-these} and is easily deduced from the theorem of Kucerovsky):
\medskip
\noindent
\begin{lem} \label{lempremierproduit}
Let $A,B,C$ be three ungraded $C^*$-algebras.\\
Let $(E_1 \hat{\otimes}_{\phi_2} E_2, \phi_1 \hat{\otimes} 1, F)$ and $(E_2, \phi_2, F_2)$ be two unbounded Kasparov modules in $\Psi^1(A,C)$ and $(E_1, \phi_1, 0) \in \Psi(A,B)$ such that :
	\begin{enumerate}[(a)]
		\item $E_1$ is trivially graded,
		\item for any $x$ in a dense subset $\mathcal{D}$ of $\phi_1(A)E_1$, 
				$$\left[ \left( 
					\begin{array}{cc}
						F & 0 \\
						0 & F_2
					\end{array} \right), \left(
					\begin{array}{cc}
						0 & T_x \\
						T_x^* & 0
					\end{array} \right) \right]$$
				is bounded on $\mathfrak{Dom} \, (F \oplus \, F_2)$, where $T_x:E_2 \rightarrow E_1 \hat{\otimes}_{\phi_2} E_2$ maps $e_2 \mapsto x \hat{\otimes} e_2$.
	\end{enumerate}
Then the module $(E_1 \hat{\otimes}_{\phi_2} E_2, \phi_1 \hat{\otimes} 1, F)$ represents the Kasparov product of the two modules $(E_1, \phi_1, 0)$ and $(E_2, \phi_2, F_2)$, i.e
$$\beta^{BJ}(E_1 \hat{\otimes}_{\phi_2} E_2, \phi_1 \hat{\otimes} 1, F) = \beta^{BJ}(E_1, \phi_1, 0) \otimes_B \beta^{BJ}(E_2, \phi_2, F_2).$$
\end{lem}

\bigskip
\noindent

\begin{lem} \label{lemdernierproduit}
Let $A,B,C$ be three ungraded $C^*$-algebras.\\
Let $(B \hat{\otimes}_{\phi_2} E_2, \phi_1 \hat{\otimes} 1, F)$ and $(B, \phi_1, F_1)$ be two unbounded Kasparov modules in $\Psi^1(A,B)$ and $(E_2, \phi_2, F_2) \in \Psi(B,C)$ such that $\phi_2$ is nondegenerated and :
	\begin{enumerate}[(a)]
	\item F is the closure of the sum $P_1 + P_2$ where $P_1$ is an unbounded selfadjoint operator of degree $0$ on $B \hti_{\phi_2} E_2$ and $P_2$ an unbounded selfadjoint operator of degree $1$ on $B \hti_{\phi_2} E_2$ with  $\mathfrak{Dom} \, F = \mathfrak{Dom} \, P_1 \cap \mathfrak{Dom} \, P_2$,
		\item for any $x$ in a dense subset $\mathcal{B}$ of $\phi_1(A)B$ (such that for all $b \in \mathcal{B}$ and $e_2 \in \mathfrak{Dom} \, F_2$, $b \hti e_2 \in \mathfrak{Dom} \, P_1 \cap \mathfrak{Dom} \, P_2$), 
				$$\left[ \left( 
					\begin{array}{cc}
						P_2 & 0 \\
						0 & F_2
					\end{array} \right), \left(
					\begin{array}{cc}
						0 & T_x \\
						T_x^* & 0
					\end{array} \right) \right]$$
				is bounded on $\mathfrak{Dom} \, (F \oplus \, F_2)$, where $T_x:E_2 \rightarrow B \hat{\otimes}_{\phi_2} E_2$ maps $e_2 \mapsto x \hat{\otimes} e_2$;
		\item For any $b \in \mathcal{B}$, 
$$ \begin{array}{ccc}
\mathfrak{Dom} \, F_2 & \longrightarrow & B \hti E_2 \\
e_2 & \longmapsto & P_1(b \hti e_2)
\end{array}$$
is bounded,
and
$$ \begin{array}{ccc}
\mathfrak{Dom} \, P_1 \cap \mathfrak{Dom} \, P_2  & \longrightarrow & E_2 \\
b' \hti e_2 & \longmapsto & T^*_{b} P_1 (b' \hti e_2)
\end{array}$$
can be extended in a bounded operator on $\mathfrak{Dom} \, F$.
		\item $\mathfrak{Dom} \, F \subset \mathfrak{Dom} \, (F_1 \hat{\otimes} 1)$ (or vice versa);
		\item $\langle (F_1 \hti 1 )x , P_2y \rangle + \langle P_2x , (F_1 \hti 1) y \rangle = 0$ for all $x, \, y$ in  $\mathfrak{Dom} \, F$.
		\item $(-1)^{\partial x} \left [ \langle (F_1 \hti 1 )x , P_1 x \rangle + \langle (P_1 \hti 1)x , (F_1 \hti 1) x \rangle \right ] \geqslant \kappa \langle x , x \rangle$ for all homogeneous $x$ in the domain $\mathfrak{Dom} \, F$.
	\end{enumerate}
Then $(B \hat{\otimes}_{\phi_2} E_2, \phi_1 \hat{\otimes} 1, F)$ represents the Kasparov product of $(B, \phi_1, F_1)$ and $(E_2, \phi_2, F_2)$, i.e
$$\beta^{BJ}(B \hat{\otimes}_{\phi_2} E_2, \phi_1 \hat{\otimes} 1, F) = \beta^{BJ}(B, \phi_1, F_1) \otimes \beta^{BJ}(E_2, \phi_2, F_2).$$
\end{lem}

\bigskip
\noindent
\subsubsection*{First Kasparov product}

\medskip
\noindent
We first prove that the homomorphism
$$K_1 \big( C(\Omega) \big) \rightarrow K_0 \big( C(\Omega) \rtimes S^1 \big)$$
is the Kasparov product by the class of the Dirac operator of dimension $1$ along the leaves of $\Omega$ foliated by $S^1$.\\
The decomposition of this map is given by :
$$\xymatrix{ K_1 \big( C(\Omega) \big) \ar[d]^{\otimes [(L^2(S^1)^{op} \otimes C(\Omega),M_{C(\Omega)},0)]} \\
K_1 \big( C(\Omega) \otimes C(S^1)\rtimes S^1 \big) \ar[d]^{\otimes [\Psi]} \\ 
K_1 \big( C(\Omega \times S^1 ) \rtimes S^1 \big) \ar[d]^{\otimes J_{S^1}(\tau_\Omega [\partial_1])} \\
K_0 \big( C(\Omega) \rtimes S^1 \big)
}$$
where $\Psi:  C(\Omega) \otimes C(S^1)\rtimes S^1 \rightarrow  C(\Omega \times S^1) \rtimes S^1$ is the densely defined $S^1$-equivariant map:$$\Psi(f)(\omega,\theta)=f(\omega.(0,\theta),\theta),$$
for $f \in C(\Omega \times S^1)$, $M_{C(\Omega)}$ is the multiplication by an element of $C(\Omega)$ and $L^2(S^1)^{op}$ is the space $L^2(S^1)$ endowed with the multiplication : $\lambda. \xi := \overline{\lambda} \xi$ for $\lambda \in \CC$ and $\xi \in L^2(S^1)$.\\
We can then simplify this map :

\medskip
\noindent
\begin{prop} \label{S1}
$$[\mathcal{E}] \otimes [\Psi] \otimes  J_{S^1}(\tau_\Omega [\partial_1]) = \big[ (C(\Omega)  \rtimes S^1,M_{C(\Omega)},d_1) \big],$$
where $\mathcal{E}=(L^2(S^1)^{op} \otimes C(\Omega),M_{C(\Omega)},0)$.
\end{prop}

\bigskip
\noindent
\textbf{Proof :}
\begin{enumerate}[(a)]
	\item In $KK_1 \big( C(\Omega) \otimes C(S^1)\rtimes S^1 , C(\Omega) \rtimes S^1 \big)$,
	$$[\Psi] \otimes  J_{S^1}(\tau_\Omega [\partial_1]) = \Psi^*( J_{S^1}(\tau_\Omega [\partial_1]))$$
	and thus:
	$$[\Psi] \otimes  J_{S^1}(\tau_\Omega [\partial_1])  =  \big[(L^2(S^1) \otimes C(\Omega) \rtimes S^1 , \pi_\alpha \circ \Psi,\partial_1 \otimes 1) \big],$$
	with $ J_{S^1}(\tau_\Omega [\partial_1]) = \big[(L^2(S^1) \otimes C(\Omega) \rtimes S^1 , \pi_\alpha ,\partial_1 \otimes 1) \big]$ (see \citeg{SkanKK}) where, for $\zeta \in C(S^1,C(\Omega))$, $a \in  C(S^1,C(S^1))$, $\xi \in L^2(S^1)$ and $f \in C(S^1 \times \Omega)$,
$$\pi_\alpha(\zeta \otimes a)(\xi \otimes f)(k,h,\omega):=\int_{S^1} a(g)(k)\xi(g^{-1}k) \zeta(g)(\omega)f(g^{-1}h,\omega.g) dg.$$

\medskip
\noindent
	Define
	$$U:L^2(S^1) \otimes C(S^1 \times \Omega) \rightarrow L^2(S^1) \otimes C(S^1 \times\Omega)$$
	by : $\forall \, f \in C(S^1 \times \Omega)$, $\zeta \in L^2(S^1)$
	$$U(\zeta 	\otimes f)(g,h,\omega):=\zeta(g)f(gh,\omega.g^{-1}).$$ 
This map extends to a unitary operator of $L^2(S^1) \otimes C(\Omega) \rtimes S^1$.

\bigskip
\noindent
	Moreover, the representation of $C(\Omega) \otimes C(S^1) \rtimes S^1$ on the module defining  $[\Psi] \otimes  J_{S^1}(\tau_\Omega [\partial_1])$ is given by the natural representation obtained from the action of $C(\Omega)$ on $C(\Omega) \rtimes S^1$ and of $C(S^1) \rtimes S^1$ by compact operator action on $L^2(S^1)$.\\
	The operator $\partial_1 \otimes 1$ is transformed in the operator $\partial_1 \otimes 1 + 1 \otimes d_1$.\\
Thus, 
$$ \begin{array}{c} 
[\Psi] \otimes  J_{S^1}(\tau_\Omega [\partial_1]) \\
\shortparallel \\
 \big[(L^2(S^1) \otimes C(\Omega) \rtimes S^1 , \pi_{C(S^1) \rtimes S^1} \otimes M_{C(\Omega)},			\partial_1 \otimes 1 + 1 \otimes d_1) \big] 
\end{array}.$$
	
	\item Since $L^2(S^1)^{op} \otimes_{C(S^1) \rtimes S^1} L^2(S^1) \simeq \mathbb{C}$, it remains to prove that 
$$(C(\Omega) \rtimes S^1, M_{C(\Omega)},d_1) \, , \, \mathcal{E}$$ 
$$\text{ and }$$ 
$$(L^2(S^1) \otimes C(\Omega) \rtimes S^1 , \pi_{C(S^1) \rtimes S^1} \otimes M_{C(\Omega)},\partial_1 \otimes 1 + 1 \otimes d_1)$$
satisfy the conditions of lemma \refg{lempremierproduit} to obtain the equality of the proposition.\\
We must verify that :
$$\left[ \left( 
					\begin{array}{cc}
						d_1 & 0 \\
						0 & \partial_1 \otimes 1 + 1 \otimes d_1
					\end{array} \right), \left(
					\begin{array}{cc}
						0 & T_{x \otimes f} \\
						T_{x \otimes f}^* & 0
					\end{array} \right) \right]$$
is a bounded operator on $\mathfrak{Dom} \, (d_1 \oplus \, (\partial_1 \otimes 1 + 1 \otimes d_1))$ for any $x \otimes f$ in a dense subset of $L^2(S^1) \otimes C(\Omega)$ (see \citeg{BaajJulg} ou \citeg{Vas} for the definition of the domain of the tensor product of an unbounded operator with the identity).\\ 
Thanks to the identification $L^2(S^1)^{op} \otimes_{C(S^1) \rtimes S^1} L^2(S^1) \simeq \mathbb{C}$ :\\
for all $x \otimes f \in \mathfrak{Dom} \, \partial_1 \otimes C_\tau^\infty(\Omega)$ (where $C_\tau^\infty(\Omega)$ is the algebra of tangentially smooth maps on $\Omega$, see \citeg{MooSch}),
$$\begin{array}{cccc} T_{x \otimes f}: & L^2(S^1) \otimes C(\Omega) \rtimes S^1 & \longrightarrow & C(\Omega) \rtimes S^1 \\
& \zeta \otimes b & \longmapsto & \langle x,\zeta\rangle_{L^2} M_{C(\Omega)}(f) b
\end{array}$$ and
$$\begin{array}{cccc} T^*_{x \otimes f}: & C(\Omega) \rtimes S^1 & \longrightarrow & L^2(S^1) \otimes C(\Omega) \rtimes S^1 \\
& b & \longmapsto & (\langle x,\zeta \rangle_{C(S^1) \rtimes S^1} . \, \zeta) \otimes M_{C(\Omega)}(\overline{f}) b \\
& & & \shortparallel \\
& & &  \overline{x} \otimes  M_{C(\Omega)}(\overline{f})b
\end{array}$$
where $\zeta \in L^2(S^1)$ is chosen such that $\langle \zeta,\zeta \rangle_{L^2}=1$ (for example the constant function equal to $1$).\\
One can easily obtain 
$$\hspace*{-1cm} \left[ \left( 
					\begin{array}{cc}
						d_1 & 0 \\
						0 & \partial_1 \otimes 1 + 1 \otimes d_1
					\end{array} \right), \left(
					\begin{array}{cc}
						0 & T_{x\otimes f} \\
						T_{x \otimes f}^* & 0
					\end{array} \right) \right]$$
$$= \left( 
					\begin{array}{cc}
						T_{x \otimes Grad_1(f)} - T_{\partial_1 x \otimes f} & 0\\
						0 & T^*_{\partial_1 x \otimes f} + T^*_{x \otimes Grad_1(f)} 
					\end{array} \right),$$
and thus this operator is bounded for any $x \otimes f$ in $\mathfrak{Dom} \, \partial_1 \otimes C_\tau^\infty(\Omega)$ which is dense in $L^2(S^1) \otimes C(\Omega)$ and this ends the proof.\\
$Grad_1$ is the gradient on $C_\tau^\infty(\Omega)$ along $S^1$ and is defined as follows :\\
Let $f$ be in $C_\tau^\infty(\Omega)$, the gradient of $f$ along $S^1$ is given by :
$$Grad_1(f)(\omega):=\lim_{\theta \rightarrow 0} \dfrac{f(\omega.\theta) - f(\omega)}{\theta}.$$
Let's remark that we have the following relation :\\
for any $f \in C_\tau^\infty(\Omega)$ and $\zeta \in\mathfrak{Dom} \, d_1$:
\begin{equation} \label{eq1}
d_1 \big( M_{C(\Omega)}(f) \zeta \big)=M_{C(\Omega)}(f) d_1 \zeta + M_{C(\Omega)} \big( Grad_1(f) \big) \zeta.
\end{equation}

\begin{flushright}
$\square$
\end{flushright}
\end{enumerate}	

\bigskip
\noindent
\subsubsection*{Second Kasparov product}
\medskip
\noindent
We now prove that the second map
$$K_0 \big( C(\Omega) \rtimes S^1 \big) \rightarrow K_0 \big( C(\Omega) \rtimes \RR^2 \rtimes S^1 \big)$$
is the Kasparov product by the class of the tangential Dirac operator $D_2$ of dimension $2$ transverse to the inclusion of the foliation by $S^1$ in the foliation by $\RR^2 \rtimes S^1$ of $\Omega$.

\bigskip
\noindent
Let's remind its definition :
$$\xymatrix{ K_0 \big( C(\Omega) \rtimes S^1\big) \ar[d]^{\otimes J_{S^1} [(L^2(\RR^2,\CC \oplus \CC)^{op} \otimes C(\Omega),M_{C(\Omega)},0)]} \\
K_0 \Big( \big( C(\Omega) \otimes C_0(\RR^2)\rtimes \RR^2 \big) \rtimes S^1 \Big) \ar[d]^{\otimes J_{S^1} \left [\Psi^{\RR^2} \right]} \\ 
K_0 \Big( \big( C_0(\Omega \times \RR^2 ) \rtimes \RR^2 \big) \rtimes S^1 \Big) \ar[d]^{\otimes J_{\RR^2 \rtimes S^1}(\tau_\Omega [\partial_2])} \\
K_0 \big( C(\Omega) \rtimes \RR^2 \rtimes S^1 \big)
}$$
where $\Psi^{\RR^2}:  C(\Omega) \otimes C(S^1)\rtimes S^1 \rightarrow  C(\Omega \times S^1) \rtimes S^1$ is the densely defined $\RR^2$-equivariant map:
$$\Psi^{\RR^2}(f)(\omega,x)=f(\omega-x,x),$$
for $f \in C(\Omega \times \RR^2)$.

\bigskip
\noindent
We define $\E{G}$ (where $G$ represents $\RR^2$ or $\RR^2 \rtimes S^1$) as the Hilbert $C^*$-module on $C(\Omega) \rtimes G$ obtained as the completion of $C_c(\Omega \times G; \CC \oplus \CC)$ for the following scalar product with values in $C_c(\Omega \times G)$
$$\langle f,f'\rangle  (\omega,g) := \int_G \langle f(\omega.h,h^{-1}) , f'(\omega.h,h^{-1}g)\rangle_\mathbb{H} dg$$
where $\langle  , \rangle_\mathbb{H}$ is given by :
$$\langle (\lambda , \nu) , (\lambda',\nu') \rangle_\mathbb{H}:=\lambda \overline{\lambda'} + \nu \overline{\nu'}.$$
This Hilbert module can be endowed with a $S^1$-action by rotation on the second summand of $\CC \oplus \CC$.

\bigskip
\noindent
As for the first Kasparov product, we can show :

\bigskip
\noindent
\begin{prop} \label{secondprodui}
$$ J_{S^1} [\mathcal{E}_2] \otimes J_{S^1} \big[ \Psi^{\RR^2} \big] \otimes  J_{\RR^2 \rtimes S^1}(\tau_\Omega [\partial_2]) = \big[ \E{\RR^2 \rtimes S^1},\pi_{C(\Omega) \rtimes S^1},D_2) \big],$$
where $\mathcal{E}_2=\big( L^2(\RR^2,\CC \oplus \CC)^{op} \otimes C(\Omega),M_{C(\Omega)},0 \big)$ and
$$\pi_{C(\Omega) \rtimes S^1}(f)(g)(\omega,x,\theta) = \int_{S^1} f(\omega,\alpha) \, \alpha \cdot g((\omega,0,\alpha)^{-1}(\omega,x,\theta)) d\alpha$$
for $f \in C(\Omega \times S^1)$, $g \in C_c(\Omega \times \RR^2 \times S^1,\CC \oplus \, \CC)$ and where $\alpha \cdot {}$ represents the action of $S^1$ on $\CC \oplus \CC$ given by trivial action on the first summand and the rotation of angle $\alpha$ on the second.
\end{prop}	

\newpage
\noindent
\textbf{Proof :}
\begin{enumerate}[(a)]
 \item We have
$$J_{S^1} [\mathcal{E}_2] \otimes J_{S^1} \big[ \Psi^{\RR^2} \big] \otimes  J_{\RR^2 \rtimes S^1}(\tau_\Omega [\partial_2]) =  J_{S^1} \Big( [\mathcal{E}_2] \otimes \big[ \Psi^{\RR^2} \big] \otimes  J^{S^1}_{\RR^2} (\tau_\Omega [\partial_2])\Big)$$ where $J^{S^1}_{\RR^2}$ is the partial descent homomorphism constructed in \citeg{Chab}.
\medskip
	\item Using the same proof as for proposition \refg{S1}, we can show that the homomorphism :
$$\xymatrix{ K_0^{S^1} \big( C(\Omega) \big) \ar[d]^{\otimes [\mathcal{E}_2]} \\
K_0^{S^1} \Big( C(\Omega) \otimes C_0(\RR^2)\rtimes \RR^2 \Big) \ar[d]^{\otimes \left [ \Psi^{\RR^2} \right ]} \\ 
K_0^{S^1} \Big( C_0(\Omega \times \RR^2 ) \rtimes \RR^2 \Big) \ar[d]^{\otimes J_{\RR^2}^{S^1}(\tau_\Omega [\partial_2])} \\
K_0^{S^1} \big( C(\Omega) \rtimes \RR^2 \big)
}$$
is the Kasparov product by $\Big[ \big( \E{\RR^2}, M_{C(\Omega)}, d_2 \big) \Big]$ where $d_2$ is the Dirac operator along the leaves of $\Omega$ foliated by $\RR^2$ that is $S^1$-invariant.

	\item Now, applying the descent homomorphism {\small $J_{S^1}$}, we obtain a class {\small $[(H',\pi',F')]$} given by :
			\begin{itemize}
				\item $H'=\E{\RR^2} \otimes_{C(\Omega) \rtimes \RR^2} C(\Omega) \rtimes \RR^2 \rtimes S^1$ where the left action of $C(\Omega) \rtimes \RR^2$ on $C(\Omega) \rtimes \RR^2 \rtimes S^1$ is given by :
$$b.\int_{S^1} a_h u_h dh := \int_{S^1} b*_{C(\Omega) \rtimes \RR^2} a_h u_h dh$$
for any $b \in C(\Omega) \rtimes \RR^2$ and $\int_{S^1} a_h u_h dh \in C(\Omega) \rtimes \RR^2 \rtimes S^1$ .
\medskip
				\item If we use the identification of $\E{\RR^2} \otimes_{ C(\Omega) \rtimes \RR^2} C(\Omega) \rtimes \RR^2 \rtimes S^1$ with $\E{\RR^2 \rtimes S^1}$, $F'$ is the tangential Dirac operator transverse to the inclusion of $S^1$ in $\RR^2 \rtimes S^1$ since $D_2(b * a_h) = (d_2 b) * a_h$ for any $h \in S^1$ and any $b \in C(\Omega) \rtimes \RR^2 \oplus \, C(\Omega) \rtimes \RR^2$.
\medskip
				\item $\pi'$ is the action of $C(\Omega) \rtimes S^1$ on $\E{\RR^2 \rtimes S^1}$ described in the proposition (see \citeg{SkanKK}).
			\end{itemize}

\begin{flushright}
$\square$
\end{flushright}
\end{enumerate}

\bigskip
\noindent
\subsubsection*{Kasparov product by the tangential Dirac operator along the leaves of $\Omega$}
\bigskip
\noindent
To end this section, it remains to show that the Kasparov product of 
$$\Big[ \big( C(\Omega) \rtimes S^1 , M_{C(\Omega)},d_1 \big) \Big]$$ 
with
$$\Big[ \big( C(\Omega) \rtimes \RR^2 \rtimes S^1 \oplus \, C(\Omega) \rtimes \RR^2 \rtimes S^1 , \pi_{C(\Omega) \rtimes S^1}, D_2 \big) \Big]$$ 
over $C(\Omega) \rtimes S^1$ is given by the class of the unbounded module 
$$\big( C(\Omega) \rtimes \RR^2 \rtimes S^1 \oplus \, C(\Omega) \rtimes \RR^2 \rtimes S^1, M_{C(\Omega)}, D_3 \big)$$ 
where $D_3$ is the tangential Dirac operator of dimension $3$ along the leaves of $\Omega$ foliated by $\RR^2 \rtimes S^1$.\\
For this, we apply lemma \refg{lemdernierproduit} to prove :

\begin{theo} \label{theoind}
$$\Big[ \big( C(\Omega) \rtimes S^1 , M_{C(\Omega)},d_1 \big) \Big] \otimes_{C(\Omega) \rtimes S^1} \Big[ \big( \mathcal{E}_{\RR^2 \rtimes S^1} , \pi_{C(\Omega) \rtimes S^1}, D_2 \big) \Big]$$ $$
= \Big[ \big( \mathcal{E}_{\RR^2 \rtimes S^1} , M_{C(\Omega)}, D_3 \big) \Big]$$
where $\mathcal{E}_{\RR^2 \rtimes S^1} = \E{\RR^2 \rtimes S^1}$.
\end{theo}

\noindent
	\textbf{Proof :}
		\begin{enumerate}[(a)]
			\item[] We prove that the three unbounded Kasparov modules
$$\Big[ \big( C(\Omega) \rtimes S^1 , M_{C(\Omega)}, d_1 - \tfrac{1}{2} Id \big) \Big] \, , \, \Big[ \big( \mathcal{E}_{\RR^2 \rtimes S^1} , \pi_{C(\Omega) \rtimes S^1}, D_2 \big) \Big]$$ 
and
$$ \Big[ \big( \mathcal{E}_{\RR^2 \rtimes S^1} , M_{C(\Omega)}, D_3  - \tfrac{1}{2} Id \big) \Big]$$
satisfy the hypotheses of lemma \refg{lemdernierproduit} with $\mathcal{B}=C(S^1,C_\tau^\infty(\Omega))$.\\
Since
$$ \Big[ \big( \mathcal{E}_{\RR^2 \rtimes S^1} , M_{C(\Omega)}, D_3 - \tfrac{1}{2}Id \big) \Big] = \Big[ \big( \mathcal{E}_{\RR^2 \rtimes S^1} , M_{C(\Omega)}, D_3 \big) \Big]$$
and
$$\Big[ \big( C(\Omega) \rtimes S^1 , M_{C(\Omega)},d_1 - \tfrac{1}{2} Id \big) \Big] = \Big[ \big( C(\Omega) \rtimes S^1 , M_{C(\Omega)},d_1 \big) \Big],$$
this would prove the theorem.
\medskip
			\item The operator $D_3$ is given, on the Hilbert field 
$$\Big( L^2(\{ \omega \} \times \RR^2 \times S^1, \CC \oplus  \CC) \Big)_{\omega \in \Omega}$$
associated to the Hilbert module $\E{\RR^2 \rtimes S^1}$ (see \citeg{Con2}), by the matrix 
$$d_3=\begin{pmatrix} \begin{matrix} 0 & -\frac{\partial}{\partial \theta} \\ 
							\frac{\partial}{\partial \theta} & 0 \end{matrix}  & \begin{matrix} -\frac{\partial}{\partial x_1} &  -\frac{\partial}{\partial x_2} \\ -\frac{\partial}{\partial x_2} & \frac{\partial}{\partial x_1}  \end{matrix}\\   
						 \begin{matrix} \frac{\partial}{\partial x_1} & \frac{\partial}{\partial x_2} \\
						\frac{\partial}{\partial x_2} & -\frac{\partial}{\partial x_1} \end{matrix} &   d''_1  
\end{pmatrix}$$
where, if $f \in C_c^\infty(\RR^2 \rtimes S^1,\CC)$ and $R_\theta$ is the rotation of angle $\theta$ on $\CC$,
$$d''_1(f)(x_1,x_2,\theta)=i R_{\theta} \frac{\partial g}{\partial \theta}(x_1,x_2,\theta)$$
and $g(x_1,x_2,\theta)=R_{-\theta}(f(x_1,x_2,\theta))$.\\
For any $\theta \in S^1$, define $U_\theta:\CC \oplus  \CC \longrightarrow \CC \oplus  \CC$ as the identity action on the first summand $\CC$ and the action by rotation of angle $-\theta$ on the second summand.\\
Thus,
$$d_3=\partial'_1 + \partial_2 \otimes 1 $$
where $\partial_2 \otimes 1$ is the tensor product of the Dirac operator on $\RR^2$ and the identity, and if $f \in C_c^\infty(\RR^2 \rtimes S^1,\CC \oplus  \CC)$, 
$$\partial'_1(f)(x_1,x_2,\theta)=iU^*_\theta\frac{\partial g}{\partial \theta}(x_1,x_2,\theta)$$
with $g(x_1,x_2,\theta)=U_\theta(f(x_1,x_2,\theta))$.\\
Thus, $D_3 = D'_1 + D_2$ where $D'_1$ (resp. $D_2$) is the operator given by $\partial'_1$ (resp. $\partial_2 \otimes 1$) on the Hilbert field associated to  $\E{\RR^2 \rtimes S^1}$.\\
Identifying $C(\Omega) \rtimes S^1 \otimes_{C(\Omega) \rtimes S^1} \mathcal{E}_{\RR^2 \rtimes S^1} $ with 
$$\pi_{C(\Omega) \rtimes S^1} (C(\Omega) \rtimes S^1) . \mathcal{E}_{\RR^2 \rtimes S^1} \simeq \mathcal{E}_{\RR^2 \rtimes S^1},$$
for  any $f \hti g \in \mathfrak{Dom}  \, (d_1 \otimes 1) \cap \mathfrak{Dom} \, D_2$,  
\hspace{-1cm}$$\begin{array}{ccc} 
					\left( D_3 - \tfrac{1}{2} Id \right) \left( f \hti g \right) & = &  \left( D_3 - \tfrac{1}{2} Id \right)  \left( \pi_{C(\Omega) \rtimes S^1}(f) g \right) \\
\end{array}$$
\hspace{1cm}$$\begin{array}{ccc} 
					& = &  \left( D'_1 - \tfrac{1}{2} Id \right)  \left( \pi_{C(\Omega) \rtimes S^1}(f) g \right) + D_2 \big( \pi_{C(\Omega) \rtimes S^1}(f) g \big) \\
					 & = & (-1)^{\partial g} \pi_{C(\Omega) \rtimes S^1}(d_1 f) g  - \tfrac{1}{2} (f \hti g)+ D_2 \big( \pi_{C(\Omega) \rtimes S^1}(f) g \big)\\
					& = & (-1)^{\partial g} (d_1 f) \hti g  - \tfrac{1}{2} (f \hti g) + D_2 \big( \pi_{C(\Omega) \rtimes S^1}(f) g \big).\\
				\end{array}$$
The sign is due to the fact that $\CC$ is acting on $\CC \oplus \CC$ by multiplication on the first summand and by conjugation on the second.

\bigskip
\noindent
We thus need to verify the hypotheses of the lemma for
$$P_1(f \hti g) = (-1)^{\partial f} d_1(f) \hti g - \tfrac{1}{2} f \hti g \text{ and } P_2 = D_2.$$

			\item If $x \in \mathcal{B}$, one can easily compute :
				$$\hspace*{-0.6cm} \left[ \left( 
					\begin{array}{cc}
						D_2 & 0 \\
						0 & D_2
					\end{array} \right), \left(
					\begin{array}{cc}
						0 & T_x \\
						T_x^* & 0
					\end{array} \right) \right]= \left( 
					\hspace{-0.2cm} \begin{array}{cc}
						T'_{Grad_2(x)} & 0\\
						0 & T'_{Grad_2(x^*)}  
					\end{array} \hspace{-0.2cm} \right)$$
where
\begin{enumerate}[(1)]
		\item 
$$D_2( \pi_{C(\Omega) \rtimes S^1}(f) g ) = \pi_{C(\Omega) \rtimes S^1}(f) D_2g + \pi_{C(\Omega) \rtimes S^1} \big( Grad_2(f) \big) \cdot g$$
since $d_2$ commutes with the $S^1$-action on $C(\Omega) \rtimes \RR^2 \oplus C(\Omega) \rtimes \RR^2$.\\
$Grad_2$ is the gradient on $C(S^1,C_\tau^\infty(\Omega))$ along $\RR^2$ and is defined as follows:
If $f \in C_\tau^\infty(\Omega)$ and $e_1,e_2$ are the two vectors of the canonical basis of $\RR^2$, then, if
$$\begin{array}{ccc}
\partial^1f(\omega):=\displaystyle{\lim_{y \rightarrow 0} \dfrac{f(\omega+ye_1) - f(\omega)}{y}}\\
\partial^2f(\omega):=\displaystyle{\lim_{y \rightarrow 0} \dfrac{f(\omega+ye_2) - f(\omega)}{y}}
\end{array},$$
the gradient along $\RR^2$ is given by :
$$\begin{array}{ccl}
Grad_2(f)(\omega) & := & \partial^1f(\omega) e_1 + \partial^2f(\omega) e_2 \\
\\
						& = &
\left( 
\begin{array}{cc}
0  & -\partial^1f(\omega)+i\partial^2f(\omega)\\
\partial^1f(\omega)+i\partial^2f(\omega) & 0
\end{array} \right)
\end{array}.$$
For any function $x$ in $C(S^1,C_\tau^\infty(\Omega))$, we still denote $Grad_2(x)$ the function of $C(S^1,C_\tau^\infty(\Omega)  \oplus C_\tau^\infty(\Omega))$ : $\theta \longmapsto Grad_2(x(\theta,.))$.\\
Moreover, we denoted $\pi_{C(\Omega) \rtimes S^1} \big( Grad_2(f) \big) \cdot g$ the product :
$$\left( \begin{array}{cc}
0  & \pi_{C(\Omega) \rtimes S^1}(-f_1+if_2) \\
\pi_{C(\Omega) \rtimes S^1}(f_1+if_2) & 0
\end{array}\right) \cdot g$$
if $Grad_2(f)= f_1 e_1 + f_2 e_2 $.

\medskip
\noindent
\item  for any $x \in \mathcal{B}$, $$\hspace{-1cm}\begin{array}{cccccc} T_x: &  \mathcal{E}_{\RR^2 \rtimes S^1} & \longrightarrow &  C(\Omega) \rtimes S^1 \otimes_{C(\Omega) \rtimes S^1}  \mathcal{E}_{\RR^2 \rtimes S^1} & \simeq  & \mathcal{E}_{\RR^2 \rtimes S^1}\\
& b & \longmapsto & x \otimes b & \simeq & \pi_{C(\Omega) \rtimes S^1} (x) b
\end{array}$$
\item for any $x \in \mathcal{B}$,
$$\begin{array}{cccc} T^*_x: &   C(\Omega) \rtimes S^1 \otimes_{C(\Omega) \rtimes S^1}  \mathcal{E}_{\RR^2 \rtimes S^1} & \longrightarrow &  \mathcal{E}_{\RR^2 \rtimes S^1} \\
& f \otimes  b & \longmapsto & \langle x,f\rangle .b
\end{array}$$
where $\langle x,f\rangle .b := \pi_{C(\Omega) \rtimes S^1} (\langle x,f\rangle )b$.

\medskip
\noindent
\item for any $x \in \mathcal{B}$, $T'_{Grad_2(x)} : \mathcal{E}_{\RR^2 \rtimes S^1} \longrightarrow \mathcal{E}_{\RR^2 \rtimes S^1}$ is given by multiplication by  
$$\left( \begin{array}{cc}
0  & \pi_{C(\Omega) \rtimes S^1}(-G_1+iG_2) \\
\pi_{C(\Omega) \rtimes S^1}(G_1+iG_2) & 0
\end{array}\right),$$
if $Grad_2(x)=G_1 e_1 + G_2 e_2$.
\end{enumerate}

\medskip
\noindent
\item  For any $b \in \mathcal{B}$ and $e_2 \in \mathfrak{Dom} \, D_2$ , $P_1(b \hti e_2) = (-1)^{\partial e_2} d_1(b) \hti e_2 - \tfrac{1}{2} b \hti e_2$ is bounded for fixed $b$ thus 
$$ \begin{array}{ccc}
\mathfrak{Dom} \, D_2 & \longrightarrow & B \hti E_2 \\
e_2 & \longmapsto & P_1(b \hti e_2) 
\end{array}$$
is bounded.

\bigskip
\noindent
Moreover, for any $b \in \mathcal{B}$ and $b' \hti e_2 \in \mathfrak{Dom} \, P_1 \cap \mathfrak{Dom} \, P_2$, 
$$T^*_{b} P_1 (b' \hti e_2) = (-1)^{\partial e_2} \langle b , d_1(b') \rangle . e_2 - \tfrac{1}{2} \langle b , b' \rangle . e_2.$$
Thereby, since $d_1$ is selfadjoint,  
$$T^*_{b} P_1 (b' \hti e_2) = \big[ T^*_{d_1(b)} - \tfrac{1}{2} T^*_b \big] (b' \hti e_2) $$
and
$$ \begin{array}{ccc}
\mathfrak{Dom} \, P_1 \cap \mathfrak{Dom} \, P_2  & \longrightarrow & E_2 \\
b' \hti e_2 & \longmapsto & T^*_{b} P_1 (b' \hti e_2)
\end{array}$$
extends to a bounded operator on $\mathfrak{Dom} \, \left( D_3 - \tfrac{1}{2} Id \right)$.

\item Using the ellipticity of $D_3$, we get
$$\mathfrak{Dom} \,  D_3 \subset \mathfrak{Dom} \, (d_1 \hti 1).$$
where $\mathfrak{Dom} \, \left( d_1 \hti 1 \right) := \text{Im}(1+d_1^2)^{-\frac{1}{2}} \hti 1$.\\
Thus, $\mathfrak{Dom} \, \left ( D_3 - \tfrac{1}{2} Id \right )  \subset \mathfrak{Dom} \, \left ( \left(d_1 - \tfrac{1}{2} Id \right ) \hti 1 \right)$.

\medskip
\noindent
\item Since $d_1 \hti 1 - \frac{1}{2} Id$ and $D_2$ are selfadjoint and 
$$\big( d_1 \hti 1 - \tfrac{1}{2} Id \big) D_2 + D_2 \big( d_1 \otimes 1 - \tfrac{1}{2} Id \big) = 0,$$
for any $x$ in the domain,
$$\big \langle \big( d_1 \hti 1 - \tfrac{1}{2} Id \big) x, D_2 x \big \rangle + \big \langle D_2 x , \big(d_1 \hti 1 - \tfrac{1}{2} Id \big) x \big \rangle = 0.$$

\item One can easily prove :
\begin{itemize}
		\item for any $x$ in the domain such that $x$ is of degree $0$, 
{\small $$\Big \langle \Big(\big( d_1 - \tfrac{1}{2} Id \big) \hti 1 \Big) x,  P_1 x \Big \rangle = \Big \langle \Big( \big( d_1 - \tfrac{1}{2} Id \big) \hti 1 \Big) x , \Big( \big( d_1 - \tfrac{1}{2} Id \big) \hti 1 \Big) x \Big \rangle \geqslant 0,$$}
		\item for any $x$ in the domain such that $x$ is of degree $1$,
{\small $$-\Big \langle \Big(\big( d_1 - \tfrac{1}{2} Id \big) \hti 1 \Big) x,  P_1 x \Big \rangle = \Big \langle \big( d_1 \hti 1 \big) x , \big( d_1 \hti 1 \big) x \Big \rangle - \tfrac{1}{4} \langle x , x \rangle \geqslant  - \tfrac{1}{4} \langle x , x \rangle ,$$}
\end{itemize}
Thus,
{\small $$(-1)^{\partial x} \left[ \Big \langle \Big(\big( d_1 - \tfrac{1}{2} Id \big) \hti 1 \Big) x,  P_1 x \Big \rangle + \Big \langle P_1 x , \Big(\big( d_1 - \tfrac{1}{2} Id \big) \hti 1 \Big) x \Big \rangle \right ] \geqslant - \tfrac{1}{4} \langle x , x \rangle,$$}
for any $x$ in the domain.\\
The last hypothesis of lemma \refg{lemdernierproduit} is thus satisfied and the theorem is thus proved.
\end{enumerate}

\begin{flushright}
$\square$
\end{flushright}

\subsection{Computation of the image under $\tau^\mu_*$ of the summands $\mathbb{Z}.q_i$}
\medskip
\noindent
In this section, we prove that $\tau_*^\mu \big(\beta \circ \delta (\tilde{q_i}) \big) \in \mu^t(C(\Xi,\ZZ))$ by an explicit computation, where $\tilde{q_i}$ is the lift of $q_i$ in $K_0^{S^1}(C(\Omega))$ built in Appendix \refg{appnoyaubord}.

\bigskip
\noindent
It suffices to make this computation for one of the $\tilde{q_i}$'s, for example $\tilde{q_2}$.\\
To compute this image, we use the index theorem for $\Omega$ foliated by $\RR^2$.\\
In fact, the following diagram is commutative  
$$\xymatrix{ K_0(C(\Omega))  \ar[d]^{\beta_2 \circ \delta_2}_\simeq && K^{S^1}_0(C(\Omega)) \ar[d]^{\beta \circ \delta}_\simeq \ar[ll]^{r^*}\\
K_0(C(\Omega) \rtimes \RR^2) \ar[dr]_{\tau_*'^\mu} & & K_0^{S^1}(C(\Omega) \rtimes \RR^2) \ar[ll]^{r^*} \ar[dl]^{\tau_*^\mu}\\
& \RR & 
}$$
where $\beta_2$ and $\delta_2$ are obtained similarly to $\beta$ and $\delta$ but forgetting the $S^1$-equivariance and $\tau_*'^\mu$ was defined in section \ref{copieZ}.\\
Thus, to compute $\tau_*^\mu(\beta \circ \delta (\tilde{q_2}))$, it suffices to compute $\tau_*'^\mu \big( \beta_2 \circ \delta_2 \circ r^*(\tilde{q_2}) \big)$.

\bigskip
\noindent
A nonequivariant version of proposition \refg{secondprodui} prove that $\beta_2 \circ \delta_2$ is given by the Kasparov product by the class of the unbounded Kasparov module given by the tangential Dirac operator along the leaves of $\Omega$ foliated by $\RR^2$ and the index theorem for foliated spaces gives :
$$\tau_*'^\mu\Big( \beta_2 \circ \delta_2 \circ r^* (\tilde{q_2}) \Big) = \left \langle Ch^2_\tau(r^*(\tilde{q_2})) \mid [C_{\nu_{Z}}] \right \rangle$$
where $Z$ is a transversal of the foliated space, $\nu_Z$ an invariant transverse measure, $[C_{\nu_{Z}}]$ its associated Ruelle-Sullivan current and $Ch^2_\tau(\tilde{q_2})$ is the Chern character component that is in $H^2_{\tau}(\Omega)$, the tangential cohomology of $\Omega$ foliated by $\RR^2$.

\bigskip
\noindent
Let's remind some notations from Appendix \refg{appnoyaubord}.\\
Let $\wo$ be the tiling in $F_1$ fixed by the rotation $R_\pi$ of angle $\pi$ around the origin and used  to build $\tilde{q_2}$.\\
$\Qi$ is the set of all the tilings with the same $1$-corona as $\wo$.\\
Fix $r_0<r<r_1$ with $r=(r_0+r_1)/2$ small enough and $r_0,r_1$ close enough to $r$.\\
Denote $\Omega_{r_0,r_1}:=\{e^{i\theta}w + v ; w \in \Qi, \theta \in \RR, v \in \RR^2, r_0 \leq ||v|| \leq r_1 \}$.\\
Since outside $\Omega_{r_0,r_1}$, the bundle defining $\tilde{q_2}$ is trivial  (see Appendix \ref{appnoyaubord}), it suffices to compute the Chern character on this corona.

\bigskip
\noindent
Let's recall that $X$ is formed by all the tilings with the same $1$-corona as $\wo$ but that are not fixed by $R_\pi$.\\
Define $X^0:=X \cup \{\wo\}$.\\
In the sequel, we will consider the transversal $Z=X \times S^1$, endowed with the invariant transverse measure $\mu^t_X \otimes d_{S^1}$ where $\mu^t_X$ is the measure induced on $X$ by $\mu$ and $d_{S^1}$ is the usual measure on $S^1$.

\bigskip
\noindent
We denote $l_k$ ($k=1,\ldots,4$) the four quadrants of the circle of radius $r$ in $\CC$ i.e $l_k=\left \{re^{ki\theta}, \theta \in \left [ 0;\frac{\pi}{2} \right]  \right \}$.\\
Define also 
$$\begin{array}{crcl}
& V_k & := & \{e^{i\theta} (w + v) ; w \in X^0, v \in l_k, \theta \in \RR \},\\
& Y_k & := & \{w + v ; w \in X^0, v \in l_k \}\\
\text{and} & \tilde{Y_k} & := & Y_k /S^1.
\end{array}$$

\noindent
Let's now fix a smooth function $\varepsilon:S^1 \rightarrow \CC$ satisfying $|\varepsilon(z)|=1$, $\varepsilon(1)=1$, $\varepsilon(i)=-1$,  $\varepsilon'(1)=0=\varepsilon'(i)$ and $\varepsilon(e^{i(\theta+\pi)})=\varepsilon(e^{i\theta})$.\\
We then define the following functions :\\
$f_1:\tilde{Y_1} \rightarrow \CC$ is given by $f_1(y)=\varepsilon (e^{i\theta})$, if $y \simeq (w,re^{i\theta})$ under the identifications $\tilde{Y_1} \simeq Y_1 \simeq X^0 \times l_1$.\\
$f_2:\tilde{Y_2} \rightarrow \CC$ is defined by $f_2(y)=-1$. \\
$f_3:\tilde{Y_3} \rightarrow \CC$ is given by $f_3(y)=-\varepsilon (e^{i\theta})$, if $y \simeq (w,re^{i\theta})$ under the identifications $\tilde{Y_3} \simeq Y_3 \simeq X^0 \times l_3$.\\
$f_4:\tilde{Y_4} \rightarrow \CC$ is defined by $f_4(y)=1$.\\
In the appendix, we didn't need a smooth function and we took $\varepsilon(z)=z^2$ but the construction of $\tilde{q_2}$ obtained in the appendix is still available with any map $\varepsilon$ satisfying the above conditions.

\bigskip
\noindent
Since $\{e^{i\theta}w_0 + v ;\theta \in \RR, v \in \RR^2, r_0 \leq ||v|| \leq r_1 \}$ is of $\mu$-measure zero (\citeg{RadSad}), it suffices to study $Ch^2_\tau(\tilde{q_2})$ on $\Omega_{r_0,r_1} \setminus  \{e^{i\theta}w_0 + v ;\theta \in \RR, v \in \RR^2, r_0 \leq ||v|| \leq r_1 \}$.\\
But this set admits a trivialization :
$$s:X \times ([0;\pi] \times S^1) \times S^1 \simeq \Omega_{r_0,r_1} \setminus  \{e^{i\theta}w_0 + v ;\theta \in \RR, v \in \RR^2, r_0 \leq ||v|| \leq r_1 \} $$
given by $(w,x,e^{i\theta},e^{i\alpha}) \mapsto e^{i\alpha} \big( w + (\frac{r_1-r_0}{\pi}x+r_0) e^{i\theta} \big)$.

\bigskip
\noindent
On $X \times ([0;\pi] \times S^1) \times S^1$ (that is a chart of $\Omega$ foliated by $\RR^2$), the projection associated to $\tilde{q_2}$ is given, modulo $s$, on $\CC \oplus \CC$ (with $S^1$-action by rotation on the first summand and trivial action on the second) by the matrix :
$$e(w,x,e^{i\theta},e^{i\alpha}) := \frac{1}{2} \begin{pmatrix} 1 + \cos(x) & \sin(x)e^{-i\alpha} \overline{f}\big([w+re^{i\theta}]\big) \\
\sin(x)e^{i\alpha} f\big([w+re^{i\theta}]\big) & 1 - \cos(x) \end{pmatrix}$$
where $f:\bigcup_{k=1}^4 \tilde{Y'_k} \longrightarrow S^1 \subset \CC$ is the continuous map given by $f(y)=f_k(y)$ if $y \in \tilde{Y'_k} \subset \tilde{Y_k}$ with $Y'_k:=\{w + v ; w \in X, v \in l_k \}$ and $\tilde{Y'_k} := Y'_k /S^1$.\\
$f([w+re^{i\theta}])$ doesn't depend on the chosen $w$ in $X$ but only on $e^{i \theta}$.

\bigskip
\noindent
Thus, $Ch^2_{\tau} (\tilde{q_2})= \dfrac{1}{2i\pi}Tr(edede)$ is the tangential form defined by
$$(w,x,z,e^{i\alpha}) \mapsto \frac{1}{4i\pi} \sin(x) \, \overline{f} \,  \frac{\partial f}{\partial z} \, dx \, dz.$$

\bigskip
\noindent
Thereby 
$$\left \langle Ch^2_\tau(\tilde{q_2}) \mid [C_{\nu_{Z}}] \right \rangle = \int_{X} \int_{S^1} \int_{[0;\pi] \times S^1} \frac{1}{4i\pi} \sin(x) \, \overline{f} \, \frac{\partial f}{\partial z} \, dx \, dz \, d_{S^1}(e^{i\alpha}) \, d\mu^t_X(\omega).$$
Moreover,
$$\frac{1}{4i\pi}\int_{[0;\pi] \times S^1}  \sin(x) \, \overline{f} \, \frac{\partial f}{\partial z} \, dx \, dz =  \frac{1}{4i\pi} \int_{[0;\pi]} \sin(x) dx \, \int_{S^1} \overline{f} \, \frac{\partial f}{\partial z} \, dz  = Ind_0(f)$$ 
where $Ind_0(f)$ is the winding number of $f$ around $0$ (we recall that $f\overline{f}=1$).\\
Since $Ind_0(f) \in \ZZ$, 
$$\left \langle Ch^2_\tau(\tilde{q_2}) \mid [C_{\nu_{Z}}] \right \rangle = l \mu^t_X(X), \text{ with } l \in \ZZ.$$
As $X_3:=X \cup \Xi$ is a transversal in $\Omega$ foliated by $\RR^2 \rtimes S^1$, $\mu$ induces an invariant transverse measure $\mu^t_3$ on $X_3$.\\
Then $\mu^t_X(X)=\mu^t_3(X)=\mu^t_3(X-v)=\mu^t_3(U\setminus \{\wo -v\})$, where $v \in \RR^2$ is such that $\wo-v \in \Xi$ and $U$ is the set formed by the tilings of $\Omega$ coinciding on the $1$-corona of $\wo$ translated by $-v$.\\
But $\mu^t_3(\{\wo - v\})=0$ thus $\mu^t_X(X)=\mu^t_3(U)=\mu^t(U)$.\\
Since $U$ is a clopen subset of $\Xi$, we deduce
$$\left \langle Ch^2_\tau(\tilde{q_2}) \mid [C_{\nu_{Z}}] \right \rangle = l \mu^t(U) \in \mu^t \big(C(\Xi,\ZZ) \big).$$
\begin{flushright}
$\square$
\end{flushright}

We thus obtained :

\begin{theo}
$$\tau_*^\mu \left( \bigoplus_{i=2}^7 \ZZ.\tilde{q}_i \right) \subset \mu^t \big( C(\Xi,\ZZ) \big).$$ \label{Z2}
\end{theo}

\section{Conclusion}

From theorem \refg{Z1}, \refg{indice} and \refg{Z2}, we thus have obtained :

\bigskip
\noindent
\textbf{Theorem :}  \textit{If $\TT$ is a pinwheel tiling, $\Omega = \Omega(\TT)$ its hull provided with an invariant ergodic probability measure $\mu$ and $\Xi$ its canonical transversal provided with the induced measure $\mu^t$, then
$$K_0(C(\Omega) \rtimes \RR^2 \rtimes S^1) \overset{\psi_2}{\simeq} \ZZ\oplus \ZZ^6 \oplus \, \check{H}^2_c \Big( (\Omega \setminus F)/S^1 \, ; \, \ZZ \Big).$$
And
\begin{itemize}
	\itemb $\tau^\mu_* \big( \psi_2(\ZZ) \big)=0$.
	\itemb $\tau^\mu_* \big( \psi_2(\ZZ^6) \big) \subset \mu^t(C(\Xi,\ZZ))$.
	\itemb $\forall b \in H$, $\exists [u] \in K_1 \big( C(\Omega) \big)$ such that :
$$\tau^\mu_*(b)=\tau^\mu_*([u] \otimes_{C(\Omega)} [D_3]) = \big \langle Ch_\tau ( [u]) , [C_{\mu^t}] \big \rangle. $$
\end{itemize}
}

To prove the gap-labeling conjecture for pinwheel tilings, it then remains to study the cohomological part of the index formula proved in the last point of this theorem.\\
This is done in \citeg{Hai-coho}, proving that 
$$\tau_*^\mu \Big( K_0 \big( C(\Omega) \rtimes \RR^2 \rtimes S^1 \big) \Big) = \mu^t \big( C(\Xi,\ZZ) \big) = \frac{1}{264} \ZZ \left[ \frac{1}{5} \right].$$

\newpage

\newpage
\appendix

\section{Proof of lemma \ref{noyaubord}} \label{appnoyaubord}
\medskip
\noindent
We prove in this appendix that the kernel of the connecting map in lemma \ref{noyaubord} is isomorphic to $\ZZ^7$ with explicit generators.

\bigskip
\noindent
	\subsection*{First inclusion}
We have already seen that $K^{S^1}_0(C(F)) \simeq \ZZ^{12}$ in section \ref{studyK} and we prove in this section that $\bigoplus \limits_{i=1}^7 \mathbb{Z}.q_i \subset \text{ Ker }\partial$ where
$$q_1=(1,1,1,1,1,1,0,0,0,0,0,0)$$
$$q_2=(0,1,1,1,1,1,1,0,0,0,0,0)$$
$$q_3=(1,0,1,1,1,1,0,1,0,0,0,0)$$
$$q_4=(1,1,0,1,1,1,0,0,1,0,0,0)$$
$$q_5=(1,1,1,0,1,1,0,0,0,1,0,0)$$
$$q_6=(1,1,1,1,0,1,0,0,0,0,1,0)$$
$$q_7=(1,1,1,1,1,0,0,0,0,0,0,1).$$
For this, we show the inclusion $\mathbb{Z}.q_1 \subset \text{ Ker }\partial$.\\
It suffices to show that $(1,\ldots,1,0,\ldots,0) \in \text{ Ker }\partial$.

\bigskip
\noindent
But $(1,\ldots,1,0,\ldots,0)$ is the projection in $K_0^{S^1} \big( C(F) \big)$ associated to the constant $S^1$-invariant projection $f$ equal to $1$ in $C(F)$.\\
The constant map $\Psi:=1$ in $C(\Omega)$ is then a $S^1$-invariant and selfadjoint map which lifts $f$ (i.e $\pi(\Psi)=f$ where $\pi : C(\Omega) \longrightarrow C(F)$).\\
Thus, the definition of $\partial$ gives :
$$\partial([f])=[\exp(-2i\pi \Psi)]=[1].$$
\medskip
Thereby $q_1 \in \text{ Ker }\partial$.

\bigskip
\noindent
The next step is to show that $\mathbb{Z}.q_2= \mathbb{Z}.(0,1,1,1,1,1,1,0,0,0,0,0) \subset \text{ Ker }\partial$.\\
The other inclusions are shown in a similar way.

\bigskip
\noindent
Let $\omega_0 \in F_1$ be a fixed tiling for the rotation $R_{\pi}$ of angle $\pi$ around the origin and let $\Qi$ be the set of the tilings having the same $1$-corona as $\wo$ (the $1$-corona of $\wo$ consists of all the tiles intersecting the tiles, in $\wo$, containing the origin). Set :
$$\Omega_r:=\left \{ e^{i \theta} \omega + v ; \theta \in \RR, w \in \Qi, v \in \RR^2 , ||v||=r \right \}$$ 
with $r$ small enough.\\
The $S^1$-action on $\Omega_r$ is free, thus we obtain a $S^1$-principal bundle $\Omega_r \rightarrow \Omega_r / S^1$.\\
It is clear, from the study of $\Omega/S^1$ made in \citeg{Hai-coho}, that $\Omega_r/S^1$ has no cohomology groups in degree greater than $2$ and thus the bundle $\Omega_r \rightarrow \Omega_r / S^1$ is trivializable. We then exhibit a trivialization of this bundle.

\bigskip
\noindent
We first exhibit a fundamental domain of $\Omega_r$ for its $S^1$-action.\\
The nonsymmetric $2$-coronas that can surround the $1$-corona around the origin in $\wo$ can be partitionned in two sets such that $R_{\pi}$ sends one on the other (a $2$-corona consists of all the tiles that intersect the $1$-corona).\\
Let $X_2$ be the union of the clopen sets of $\Qi$ associated to one of these sets.\\
Next, we consider the nonsymmetric $3$-coronas of the $1$-corona around the origin in $\wo$ such that the $2$-corona is symmetric.\\
We can partition as above and we obtain a clopen set $X_3$.\\
Then we consider nonsymmetric $4$-coronas of the $1$-corona around the origin in $\wo$ such that the $3$-corona is symmetric and so on.\\
We have the clopen sets $X_2, X_3, X_4, \ldots$ of $\Qi$ and if we denote $X = \bigcup \limits_{n \geqslant 2} X_n$, then $X$ is an open set of $\Qi$ and $X$, $\{\wo\}$, $-X$ is a partition of $\Qi$.\\
We denote $X^0:=X \cup \{\wo\}$.\\
In the sequel, we call $l_k$ ($k=1,\ldots,4$) the four quadrants of the circle of radius $r$ around the origin of $\CC$ i.e $l_k=\left \{re^{ki\theta}, \theta \in \left [ 0;\frac{\pi}{2} \right] \right \}$.\\
We also set 
$$V_k:=\{e^{i\theta} (w + v) ; w \in X^0, v \in l_k, \theta \in \RR \}$$
and $Y_k:=\{w + v ; w \in X^0, v \in l_k \}$.

\bigskip
\noindent
We first prove that $V_k \rightarrow V_k/S^1$ is a trivial principal bundle homeomorphic to $Y_k \times S^1$.
\medskip
\begin{itemize}
	\item[] Let $\Psi_k:Y_k \times S^1 \longrightarrow V_k$ be defined by $\psi_k(w+v,e^{i\theta}):=e^{i\theta}(w+v)$. We show that this is a homeomorphism.

\bigskip
\noindent
$\psi_k$ is easily seen to be surjective from the definition of $Y_k$ and $V_k$.

\bigskip
\noindent
If $\psi_k(w+v,e^{i\theta})=\psi_k(w'+v',e^{i\theta'})$ then 
$$e^{i(\theta-\theta')}(w+v)-v'=w'.$$
If $r$ is small enough, $e^{i(\theta-\theta')}(v)-v'$ can send $e^{i(\theta-\theta')}w$ on another tiling of $X^0$ only if it vanishes.\\
Thus, $e^{i(\theta-\theta')}(v)=v'$ and $e^{i(\theta-\theta')}(w)=w'$.\\
The only possibility is $\theta'=\theta$ modulo $\pi$.\\
But, if $\theta'=\theta+\pi+2l\pi$ ($l \in \ZZ$), then $v'=-v$, which is not possible since $v',v \in l_k$.\\
Thereby, $\theta'=\theta + 2l\pi$ and $(w+v,e^{i\theta})=(w'+v',e^{i\theta'})$.\\
And $\psi_k$ is injective.

\bigskip
\noindent
Since $\psi_k$ is a continuous bijection from the compact set $Y_k\times S^1$ on the compact set $V_k$, it is a homeomorphism.
\end{itemize}

\bigskip
\noindent
Thereby, $Y_k \simeq V_k / S^1$ and, denoting $\tilde{Y_k}:=Y_k/S^1$, $\tilde{Y_k}$ is homeomorphic to $Y_k$. We will call $s_k$ this homeomorphism.\\
As $Y_k \simeq X^0 \times l_k$, we also note that we can identify $\tilde{Y_k}$ with $X^0 \times l_k$.

\bigskip
\noindent
We next study the intersection $\tilde{Y_k} \cap \tilde{Y_{k'}}$ ($k \neq k'$) in $\Omega_r/S^1$.
\medskip
\begin{itemize}
	\item[] If $[w+v]=[w'+v']$ with $[w+v] \in \tilde{Y_k}$ and $[w'+v'] \in \tilde{Y_{k'}}$, there exists $\theta \in \RR$ such that
$$e^{i\theta}(w+v)=w'+v'.$$
As above, $v'$ must be equal to $e^{i\theta}v$ and thus $w'=e^{i\theta}w$.\\
There are then only two possibilities:
\medskip
\begin{itemize}
	\itemb if $w \neq w_0$, then $\theta=2l\pi$ ($l \in \ZZ$) and $v'=v$, $w'=w$.
	\itemb if $w=w_0$, then $\theta=l\pi$ ($l \in \ZZ$) and thus $w'=w_0$, $v'=(-1)^l v$.
\end{itemize}
Thereby, 
$$\tilde{Y_k} \cap \tilde{Y_{k'}} \simeq X^0 \times (l_k \cap l_{k'}) \bigcup \{w_0\} \times (-l_k \cap l_{k'}).$$
\end{itemize}

\bigskip
\noindent
We then define $4$ maps allowing us to glue the $(\tilde{Y_k} \times S^1)$'s to identify the union of these sets with $\Omega_r$.\\
We can identify each $Y_k$ with $X^0 \times l_k$ which allow us to set $p_k$ the projection on the second coordinate.\\
In the sequel, we consider the $l_k$'s as imbedded in $\CC$ which allow us to use the square of elements of $l_k$.\\
$f_1:\tilde{Y_1} \rightarrow \CC$ is then defined by $f_1(y)=\frac{p_1(s_1(y))^2}{r^2}$. Under the identifications $\tilde{Y_1} \simeq Y_1$ and $Y_1 \simeq X^0 \times l_1$, $f_1$ is obtained by : if $y \simeq (w,re^{i\theta})$, $f_1(y)=e^{2i\theta}$.\\
$f_2:\tilde{Y_2} \rightarrow \CC$ is defined by $f_2(y)=-1$. \\
$f_3:\tilde{Y_3} \rightarrow \CC$ is defined by $f_3(y)=-\frac{p_3(s_3(y))^2}{r^2}$. $f_3$ is obtained, under the identifications $\tilde{Y_3} \simeq Y_3 \simeq X^0 \times l_3$, in the following way : if $y \simeq (w,re^{i\theta})$, $f_3(y)=-e^{2i\theta}$.\\
$f_4:\tilde{Y_4} \rightarrow \CC$ is defined by $f_4(y)=1$.

\bigskip
We introduce these maps because the transition maps on the intersection of $\tilde{Y_k}$ with $\tilde{Y_{k'}}$ are given by $f_k \overline{f_{k'}}$.\\
Let's prove that we can then define a map $\phi: \left (\bigcup \limits_{k=1}^4 \tilde{Y_k} \right ) \times S^1 \longrightarrow \Omega_r$ by $\phi(y,e^{i\theta}):=f_k(y)e^{i\theta}s_k(y)$ and that this is a homeomorphism.

\bigskip
\noindent
\begin{itemize}
	\item[] If $y=[w+v] \in \tilde{Y_k} \cap \tilde{Y_{k'}}$, then there are two possibilities :
\begin{itemize}
	\itemb if $w \neq w_0$, then $k=k'\pm 1 \mod{4}$.\\
We can assume, even if it means interchanging the roles of $k$ and $k'$, that $k=k'+1 \mod{4}$. We then have $s_k([w+v])=w+ri^{k-1}$ and $s_{k'}([w+v])=w+ri^{k-1}$. One can easily prove that $f_k(y)=f_{k'}(y)$ and thus $f_k(y)e^{i\theta}s_k(y)=f_{k'}(y)e^{i\theta}s_{k'}(y)$.
	\itemb if $w=w_0$, then 
\begin{itemize}
	\item either $k=k' \pm 1 \mod{4}$ and we can suppose $k=k' +1 \mod{4}$.\\ 
Then $s_k([w+v])=w_0+ri^{k-1}$ or $s_k([w+v])=w_0+ri^k$. In the first case, we obtain $s_{k'}([w+v])=w_0+ri^{k-1}$ and in the second, we have $s_{k'}([w+v])=w_0+ri^{k-2}$.\\
In the first case, one can easily verify that $f_k(y)=f_{k'}(y)$ and, since $s_{k'}([w+v])=s_{k}([w+v])$, we thus obtain that 
$$f_k(y)e^{i\theta}s_k(y)=f_{k'}(y)e^{i\theta}s_{k'}(y).$$ 
In the second case, $f_k(y)=-f_{k'}(y)$ and we also obtain that $f_k(y)e^{i\theta}s_k(y)=f_{k'}(y)e^{i\theta}s_{k'}(y)$,  since $s_{k'}([w+v])=-s_{k}([w+v])$.
\medskip
	\item or $k=k' \pm 2 \mod{4}$ and if $s_k([w+v])=w_0+v$, then we have ${\displaystyle s_{k'}([w+v])=w_0-v}$.\\
Again, one can verify that  $f_k(y)e^{i\theta}s_k(y)=f_{k'}(y)e^{i\theta}s_{k'}(y)$ since, if $k=k' \pm 2 \mod{4}$, $f_k=-f_{k'}$.
\end{itemize}
\end{itemize}
Thus, $\phi$ is well defined.

\bigskip
\noindent
Now, assume $y \in \tilde{Y_k}$ and $y' \in \tilde{Y_{k'}}$ satisfy $\phi(y,e^{i\theta})=\phi(y',e^{i\theta'})$.\\
Note $s_k(y)=w+v$ and $s_{k'}(y')=w'+v'$.\\
One has $f_k(y) e^{i\theta} (w+v)  = f_{k'}(y') e^{i\theta'}( w' +v')$.\\
Thus, projecting in $\Omega_r/S^1$, $y=y'$ is an element of the intersection $\tilde{Y_k} \cap \tilde{Y_{k'}}$.\\
It is enough then to prove that $\theta= \theta' +2l\pi$, $l \in \ZZ$.\\
As above, we must have $v'=f_{k'}(y)^{-1} f_k(y) e^{i(\theta-\theta')} (v)$ and thereby 
$$w'=f_{k'}(y)^{-1} f_k(y) e^{i(\theta-\theta')} (w).$$
There are several cases :
\begin{itemize}
	\itemb either $f_{k'}(y)^{-1} f_k(y) e^{i(\theta-\theta')} =1$ and $w'=w$, $v'=v$. Then $k$ is equal to $k'$, $k'+1$ or $k'-1 \mod{4}$.\\
In the three cases, $f_k(y)=f_{k'}(y)$ thus $e^{i(\theta-\theta')}=1$. So $e^{i\theta}=e^{i\theta'}$.
	\itemb or $f_{k'}(y)^{-1} f_k(y) e^{i(\theta-\theta')} = - 1$ then $w=w_0=-w'$ and $v'=-v$. Since, for two opposite points, we have $f_k(y)=-f_{k'}(y')$, we also obtain that  $e^{i\theta}=e^{i\theta'}$.
\end{itemize}	
In every cases, we have $(y,e^{i\theta}) =(y',e^{i\theta'})$ thus $\phi$ is injective.

\bigskip
\noindent
One can easily see that $\phi$ is also surjective and continuous thus $\phi$ is a homeomorphism between $\left (\bigcup \limits_{k=1}^4 \tilde{Y_k} \right ) \times S^1$ and $\Omega_r$.
\end{itemize}

\bigskip
\noindent
Thereby, we have constructed a trivialization of the projection $\Omega_r \longrightarrow \Omega_r/S^1$.

\bigskip
\noindent
We can now construct the fiber bundle on $\Omega$ lifting
$$(0,1,1,1,1,1,1,0,0,0,0,0) \in K_0^{S^1}(C(F)).$$
To do this, let's take a tiling $w_0 \in F_1$ fixed by $R_{\pi}$ and any $r$ small enough.\\
The above construction gives $\Omega_r$, $\Omega_{\leq r}$ and $H=^{\;\;c}(\Omega_{< r})$ where
$$\Omega_{\leq r}:=\left \{ e^{i \theta} \omega + v ; \theta \in \RR, w \in \Qi, v \in \RR^2 , ||v|| \leq r \right \}$$ 
and
$$\Omega_{<r}:=\left \{ e^{i \theta} \omega + v ; \theta \in \RR, w \in \Qi, v \in \RR^2 , ||v||<r \right \}.$$
On $\Omega_{\leq r}$, we take the line bundle $\Omega_{s \leq r} \times \CC$ with the diagonal action of $S^1$ where the $S^1$-action on $\CC$ is given by rotation.\\
On $H$, we take the line bundle $H \times \CC$ with the diagonal action of $S^1$ but, this time, the $S^1$-action on $\CC$ is trivial.\\
The gluing map is then given on the intersection $\Omega_{\leq r} \cap H = \Omega_r \simeq \Omega_r/S^1  \times S^1$ by
$$f(w,e^{i\theta},z)=(w,e^{i\theta},e^{i\theta}z).$$
Thus, since $(0,1,1,1,1,1,1,0,0,0,0,0) \in K_0^{S^1}(C(F))$ lifts on an element of $K_0^{S^1}(C(\Omega))$, it is in the kernel of $\partial$.

\bigskip
\noindent
Thereby $\bigoplus \limits_{i=1}^7 \mathbb{Z}.q_i \subset \text{Ker }\partial$.

\bigskip
\noindent
\subsection*{Inclusion in the other direction}
We next prove that $\text{Ker }\partial \subset \bigoplus \limits_{i=1}^7 \mathbb{Z}.q_i$.\\
Let $(n_1, \ldots, n_6,n'_1,\ldots, n'_6)$ be in $\mathbb{Z}^{12}$.\\
We can assume that the $n'_i$'s are zero thanks to the last section, since, for example, as $(1,0,0,0,0,0,-1,0,0,0,0,0) \in \text{ Ker }\partial $, 
$$\partial(n_1, \ldots, n_6,n'_1,\ldots, n'_6) = \partial(n_1+n'_1, n_2,\ldots, n_6,0,n'_2,\ldots, n'_6).$$
Moreover, we can also assume that every $n_i+n'_i$ are in $\mathbb{N}$ thanks to the last section, since, if $n$ is the smallest integer among the $(n_i+n'_i)$'s, we have 
$$\partial(n_1+n'_1,\ldots,n_6+n'_6,0,\ldots,0)=\partial(n_1+n'_1-n,\ldots,n_6+n'_6-n,0,\ldots,0),$$
 with $n_i+n'_i-n \in \mathbb{N}$.\\
Let's fix $(n_1, \ldots, n_6,0,\ldots,0) \in \ZZ^{12}$ with $n_i \in \NN$.\\
$(n_1, \ldots, n_6,0,\ldots, 0)$ is associated with the $S^1$-invariant projection :
$$f(x) = \left \{
\begin{array}{cl}
P_{n_1} & \text{if } x \in F_1 \\
\vdots & \\
P_{n_6} & \text{if } x \in F_6 \\
\end{array} \right .$$
where
$$P_{n_i} := \left (
\begin{array}{ccc}
1 & &  \\
   & \ddots  &\\
 & & 1   \\
\end{array} \right ) \in \mathcal{M}_{n_i}(\mathbb{C}).
$$
\text{ }\\
We will exhibit a $S^1$-invariant projection in $\mathcal{M}_\infty(C(\Omega))$ which lifts this projection.\\
For this, consider 
$$\tilde{K_i}:=\Big \{ \omega \in \Omega \; ; \; d(\omega,F_i) \leqslant \dfrac{1}{n} \Big \}$$
 and
$$\tilde{V_i}:=\Big \{ \omega \in \Omega \; ; \; d(\omega,F_i) < \dfrac{1}{n-1} \Big \}$$
with $n$ big enough so that the $\tilde{V_i}$'s are disjoint.\\
Set $K_i:=p(\tilde{K_i})$ and $V_i:=p(\tilde{V_i})$ where $p : \Omega \longrightarrow \Omega / S^1$.\\
$K_i$ is then a compact subset of  the open set $V_i$ and it contains $[x_i]:=p(x_i)$ where $x_i$ is some element of $F_i$.\\
Thus, there exist maps $\phi_i \in C(\Omega / S^1)$, $i \in \{1, \ldots, 6 \}$, such that 
\begin{itemize}
\item $Supp(\phi_i) \subset V_i$;
\item $0 \leqslant \phi_i \leqslant 1$ and $\phi_i(x)=1$ if $x \in K_i$ ;
\item $\phi_1 + \ldots + \phi_6 = 1$ on $\bigcup K_i$.
\end{itemize}
\medskip
Set $\Psi:=\phi_1.P'_{n_1} + \ldots + \phi_6.P'_{n_6}$ where
$$P'_{n_i} := \left (
\begin{array}{cccccc}
1 & & & & & \\
   & \ddots & & & & \\
 & & 1 & & & \\
& & & 0 & & \\
& & & & \ddots & \\
& & & & & 0
\end{array} \right ) \in \mathcal{M}_{Max(n_i)}(\mathbb{C}).
$$
Thereby, $\Psi \in  \mathcal{M}_{Max(n_i)} \left( C \left ( \Omega / S^1 \right) \right)$ and $\Psi$ is selfadjoint.\\
We can then see $\Psi$ as a function of $C(\Omega)$, constant on the $S^1$-orbits and $\Psi$ can thus be seen as a $S^1$-invariant function.\\
One easily see that $\pi_*(\Psi)= f$ in $\mathcal{M}_\infty(C(F))$, thereby
$$\partial([f]) = [\exp(-2 i \pi \Psi)] \in K_1^{S^1} ( C_0 \big( \Omega \setminus F) \big).$$
By lemma \refg{coho}, we have an isomorphism between $K_1^{S^1}(C_0(\Omega \setminus F))$ and
$$\check{H}_c^1 \Big( (\Omega \setminus F)/S^1 , \mathbb{Z} \Big) = \left \{ [h] \; ; \; h:  \Big( (\Omega \setminus F)/S^1 \Big)^+ \longrightarrow S^1 \text{ continuous } \right \}$$ 
where $[h]$ is the class of continuous maps $\Big( (\Omega \setminus F)/S^1 \Big)^+ \longrightarrow S^1$ homotopic to $h$ and $\Big( (\Omega \setminus F)/S^1 \Big)^+$ is the Alexandroff compactification of $(\Omega \setminus F)/S^1$ (see \citeg{Hat} for the equality of the cohomology group and the set of classes of continuous maps from $\Big( (\Omega \setminus F)/S^1 \Big)^+$ to $S^1 \Big)$.\\
This isomorphism is given by $det_*$.\\
Thus : $\partial([f]) = [\exp(-2 i \pi (n_1.\phi_1 + \ldots + n_6. \phi_6))] \in \check{H}_c^1 \Big( (\Omega \setminus F)/S^1 \, ; \, \ZZ \Big)$.

\medskip
\noindent
Note that a continuous map on $\Big( (\Omega \setminus F)/S^1 \Big)^+$ can be seen as a continuous map on $\Omega / S^1$, constant on $F/S^1$.

\bigskip
\noindent
It remains to find the $(n_1, \ldots, n_6,0, \ldots,0)$'s for which there exists a continuous homotopy $H_t$ defined on $\Big( (\Omega \setminus F)/S^1 \Big)^+$ and with values in $S^1$ between $\exp(-2 i \pi (n_1.\phi_1 + \ldots + n_6. \phi_6))$ and $1$.\\
This is equivalent to know if there exists a continuous homotopy $H_t:\Omega/S^1 \rightarrow S^1$ between $exp(-2 i \pi (n_1.\phi_1 + \ldots + n_6. \phi_6))$ and $1$ which is constant on $F/S^1$ for any $t$.\\
Or, in other words, we want to know if we can find a continuous homotopy $h_t : \Omega/S^1 \longrightarrow \mathbb{R}$ between $h_0:=n_1.\phi_1 + \ldots + n_6. \phi_6$ and a continuous map $h_1:=g$ defined on $\Omega/S^1$ with values in $\mathbb{Z}$ and such that $h_t([x_i])-h_t([x_j]) \in \ZZ$ for any $t \in [0,1]$, $i, \,j \in \{1, \ldots, 6\}$.

\medskip
\noindent
Since $\Omega/S^1$ is connected, $h_t$ would be a homotopy between $n_1.\phi_1 + \ldots + n_6. \phi_6$ and a constant map $g$ on $\Omega/S^1$ equal to an integer and satisfying 
$$h_t([x_i])-h_t([x_j]) \in \ZZ \text{ pour tout } t, i,j.$$

\noindent
We prove that this is possible only if $n_i=n_j$ for all $i,j$.\\
Let $ev_i(h) : [0;1] \longrightarrow \mathbb{R}$ be the continuous map defined by $ev_i(h)(t):=h_t([x_i])$.\\
$ev_i(h)-ev_j(h)$ is then a continuous map with integer values and thus is constant equal to $n_i-n_j$.\\
Thereby, $\forall t \in [0;1]$, $h_t([x_i])-h_t([x_j]) = n_i-n_j$.\\
But, for $t=0$, $h_0=g=k \in \mathbb{Z}$ and so, if there exists such a homotopy $h_t$, then $n_i=n_j$ for all $i,j$.\\
Thus, 
$$(n_1, \ldots, n_6,0,\ldots,0) \in \text{Ker } \partial \, \Rightarrow  \, (n_1, \ldots, n_6,0,\ldots,0) \in \mathbb{Z}.(1, \ldots,1,0,\ldots,0).$$

\noindent
If now $(n_1, \ldots, n_6,n'_1,\ldots,n'_6) \in \text{Ker} \partial$, 
$$\begin{array}{rcl}
\partial(n_1, \ldots, n_6,n'_1,\ldots,n'_6) & = & \partial(n_1+n'_1, \ldots, n_6+n'_6,0,\ldots,0)\\
& = & \partial(n_1+n'_1-n, \ldots, n_6+n'_6-n,0,\ldots,0)
\end{array}$$
where $n= Min(n_i+n'_i)$.\\
From the above result, we have $n_i+n'_i-n=n_j+n'_j-n$, i.e $n_i+n'_i=n_j+n'_j$ for any $i,j$.\\
Let $k$ be this common integer, we have obtained that $n_i=-n'_i+k$, thus
$$\begin{array}{rl}
(n_1, \ldots, n_6,n'_1,\ldots,n'_6) & = \; (-n'_1+k, \ldots,-n'_6+k,n'_1,\ldots,n'_6)\\
& = \; k.q_1+n'_1(q_2-q_1)+\ldots + n'_6(q_7-q_1) \in \bigoplus \limits_{i=1}^7 \ZZ.q_i.
\end{array}$$
Thereby, we have proved that $\text{Ker } \partial \subset \bigoplus \limits_{i=1}^7 \ZZ.q_i $.

\bigskip
\noindent
This completes the proof of lemma \ref{noyaubord}.



\bigskip
\noindent

\newpage
\bibliographystyle{alpha}
\bibliography{Pinwheel}

\begin{thebibliography}{SBGC84}

\bibitem[BBG06]{BelBenGam}
J.~Bellissard, R.~Benedetti, and J.-M. Gambaudo.
\newblock \textit{Spaces of tilings, finite telescopic approximation and gap
  labelings}.
\newblock {\em Comm. Math. Phys.}, \textbf{261}:1--41, 2006.

\bibitem[Bel82]{Bel82}
J.~Bellissard.
\newblock \textit{Schr\"odinger's operator with an almost periodic potential :
  an overview}.
\newblock {\em Lecture Notes in Physics}, \textbf{153}, 1982.

\bibitem[Bel86]{Bel86}
J.~Bellissard.
\newblock \textit{K-theory of $C^*$-algebras in Solid State Physics}.
\newblock {\em Lecture Notes in Physics}, \textbf{257}:99--156, 1986.

\bibitem[Bel92]{Bel1}
J.~Bellissard.
\newblock \textit{Gap labelling theorems for Schr\"odinger operators}.
\newblock In {\em From number theory to physics (Les Houches, 1989)}, pages
  538--630, 1992.
\newblock Springer, Berlin.

\bibitem[BG03]{BenGam}
R.~Benedetti and J-M Gambaudo.
\newblock \textit{On the dynamics of $\mathbb{G}$-Solenoids. Applications to
  Delone sets}.
\newblock {\em Ergod. Th. \& Dynam. Sys.}, \textbf{29}:673--691, 2003.

\bibitem[BHZ00]{BelHerZar}
J.~Bellissard, D.J.L. Herrmann, and M.~Zarrouati.
\newblock \textit{Hulls of aperiodic solids and gap labeling theorems}.
\newblock {\em CRM Monogr. Ser.}, \textbf{13}:207--258, 2000.
\newblock A.M.S., Providence.

\bibitem[BJ83]{BaajJulg}
S.~Baaj and P.~Julg.
\newblock \textit{Th\'eorie bivariante de Kasparov et op\'erateurs non born\'es
  dans les $C\sp{*} $-modules hilbertiens}.
\newblock {\em C. R. Acad. Sci. Paris Sr. I Math.}, \textbf{296 no.
  21}:875--878, 1983.

\bibitem[BKL01]{BelKelLeg}
J.~Bellissard, J.~Kellendonk, and A.~Legrand.
\newblock \textit{Gap-labelling for three-dimensional aperiodic solids}.
\newblock {\em C.R.A.S, serie I}, \textbf{332}:521--525, 2001.

\bibitem[BOO02]{BenOyo}
M.-T. Benameur and H.~Oyono-Oyono.
\newblock \textit{Index theory for quasi-crystals. I. Computation of the
  gap-label group}.
\newblock \textbf{252}:137--170, 2002.
\newblock J. Funct. Anal.

\bibitem[Bre72]{Bred}
G.E. Bredon.
\newblock {\em Introduction to compact transformation groups}.
\newblock Pure and applied mathematics \textbf{46}, 1972.
\newblock Academic Press.

\bibitem[Cha99]{Chab}
J.~Chabert.
\newblock {\em Stabilit\'e de la conjecture de Baum-Connes pour certains
  produit semi-directs de groupes}.
\newblock PhD thesis, Univ. de la M\'edit\'erran\'ee Aix-Marseille II, 1999.

\bibitem[Con79]{Con1}
A.~Connes.
\newblock \textit{Sur la th\'eorie non commutative de l'int\'egration}.
\newblock {\em Lecture Notes in Math.}, \textbf{725}:19--143, 1979.
\newblock Springer, New York.

\bibitem[Con82]{Con2}
A.~Connes.
\newblock \textit{A survey of foliations and operator algebras}.
\newblock {\em Proc. Sympos. Pure Math.}, \textbf{38} part. \textbf{1}, 1982.
\newblock A.M.S. , Providence.

\bibitem[DHK91]{DouHurKam}
R.G. Douglas, S.~Hurder, and J.~Kaminker.
\newblock \textit{The Longitudinal Cocycle and the Index of Toeplitz
  Operators}.
\newblock {\em J. Funct. Anal.}, \textbf{101}:120--144, 1991.

\bibitem[Gre69]{Gre}
F.~P. Greenleaf.
\newblock \textit{Invariant means on topological groups and their
  applications}.
\newblock {\em Van Nostrand mathematical studies}, \textbf{16}, 1969.

\bibitem[Hat02]{Hat}
A.~Hatcher.
\newblock {\em \textit{Algebraic topology}}.
\newblock 1st ed., Cambridge University Press, 2002.

\bibitem[HRS05]{HolRadSad}
C.~Holton, C.~Radin, and L.~Sadun.
\newblock \textit{Conjugacies for Tiling Dynamical Systems}.
\newblock {\em Comm. Math. Phys.}, \textbf{254}:343--359, 2005.

\bibitem[HS87]{HilSkan}
M.~Hilsum and G.~Skandalis.
\newblock \textit{Morphismes $K$-orient\'es d'espaces de feuilles et
  fonctorialit\'e en th\'eorie de Kasparov (d'apr\`es une conjecture d'A.
  Connes)}.
\newblock {\em Ann. Sci. \'Ecole Norm. Sup. (4)}, \textbf{20 no. 3}:325--390,
  1987.

\bibitem[Jul81]{Julg}
P.~Julg.
\newblock \textit{$K$-th\'eorie \'equivariante et produits crois\'es}.
\newblock {\em Note C.R.A.S. Paris}, \textbf{292}:629--632, 1981.

\bibitem[Kas88]{Kaspa2}
G.G. Kasparov.
\newblock \textit{Equivariant $KK$-theory and the novikov conjecture}.
\newblock {\em Inv. Math.}, \textbf{91}:147--201, 1988.

\bibitem[Kas95]{Kaspa1}
G.G. Kasparov.
\newblock \textit{$K$-theory, group $C^*$-algebras and higher signatures}.
\newblock {\em London Math. Soc. Lecture Note Ser.}, \textbf{226}:101--146,
  1995.

\bibitem[Kel95]{Kelgap}
J.~Kellendonk.
\newblock \textit{Noncommutative geometry of tilings and gap labeling}.
\newblock {\em Rev. Math. Phys.}, \textbf{7}:1133--1180, 1995.

\bibitem[KP00]{KelPut}
J.~Kellendonk and I.F. Putnam.
\newblock {\em Tilings, $C^*$-algebras and K-theory}.
\newblock CRM monograph Series \textbf{13}, 177-206, 2000.
\newblock M.P. Baake and R.V. Moody Eds., A.M.S., Providence.

\bibitem[KP03]{KamPut}
J.~Kaminker and I.~Putnam.
\newblock \textit{A proof of the gap labeling conjecture}.
\newblock {\em Michigan Math. J.}, \textbf{51}:537--546, 2003.

\bibitem[Kuc97]{Kuc}
D.~Kucerovsky.
\newblock \textit{The $KK$-product of unbounded modules}.
\newblock {\em $K$-theory}, \textbf{11 no. 1}:17--34, 1997.

\bibitem[LP03]{LagPle}
J.C. Lagarias and P.A.B. Pleasant.
\newblock \textit{Repetitive Delone sets and quasicrystals}.
\newblock {\em Ergod. Th. \& Dynam. Sys.}, \textbf{23}:831--867, 2003.

\bibitem[Mat]{Mat}
M.~Matthey.
\newblock {\em K-theories, $C^*$-algebras and assembly maps}.
\newblock PhD thesis, Universit\'e de Neuch\^atel.

\bibitem[Mou]{Hai-coho}
H.~Moustafa.
\newblock \textit{PV cohomology of pinwheel tilings, their integer group of
  coinvariants and gap-labelling}.
\newblock \\http://arxiv.org/abs/0906.2107, to appear in Comm. Math. Phys.

\bibitem[Mou09]{Hai-these}
H.~Moustafa.
\newblock {\em Gap-labeling des pavages de type pinwheel}.
\newblock PhD thesis, Univ. Blaise Pascal, Clermont-Ferrand, 2009.
\newblock {\small
  \\http://math.univ-bpclermont.fr/$\sim$moustafa/These/these-Moustafa.pdf}.

\bibitem[MS06]{MooSch}
C.~C. Moore and C.~Schochet.
\newblock \textit{Global analysis on foliated spaces}.
\newblock {\em MSRI Publications}, \textbf{9}, 2006.

\bibitem[ORS02]{OrmRadSad}
N.~Ormes, C.~Radin, and L.~Sadun.
\newblock \textit{A homeomorphism invariant for substitution tiling spaces}.
\newblock {\em Geometriae Dedicata}, \textbf{90}:153--182, 2002.

\bibitem[PB09]{PeaBel}
J.~Pearson and J.~Bellissard.
\newblock \textit{Noncommutative Riemannian Geometry and Diffusion on
  Ultrametric Cantor Sets}.
\newblock {\em Journal of Noncommutative Geometry}, \textbf{3}:847--865, 2009.

\bibitem[Ped79]{Ped}
G.~K. Pederson.
\newblock \textit{$C^*$-algebras and their automorphism groups}.
\newblock {\em London Math. Society Monographs}, \textbf{14}, 1979.
\newblock Academic Press, London.

\bibitem[Pet05]{Pet}
S.~Petite.
\newblock {\em Pavages du demi-plan hyperbolique et laminations}.
\newblock PhD thesis, Univ. de Bourgogne I, 2005.

\bibitem[Rad94]{Rad}
C.~Radin.
\newblock \textit{The pinwheel tilings of the plane}.
\newblock {\em Ann. of Math.}, \textbf{139}:661--702, 1994.

\bibitem[Rad95]{Rad1}
C.~Radin.
\newblock \textit{Space tilings and substitutions}.
\newblock {\em Geom. Dedicata}, \textbf{55}:257--264, 1995.

\bibitem[Rie82]{Rief}
M.A. Rieffel.
\newblock \textit{Morita equivalence for operator algebras}.
\newblock {\em Proc. of Symposia in Pure Math.}, \textbf{38}:285--298, 1982.

\bibitem[RS98]{RadSad}
C.~Radin and L.~Sadun.
\newblock \textit{An algebraic invariant for substitution tiling systems}.
\newblock {\em Geom. Dedicata}, \textbf{73}:21--37, 1998.

\bibitem[Sad]{SadPriv}
L.~Sadun.
\newblock \textit{Private conversation in september 2007}.

\bibitem[SBGC84]{Sheetal}
D.~Shechtman, I.~Blech, D.~Gratias, and J.V. Cahn.
\newblock \textit{Metallic phase with long range orientational order and no
  translational symmetry}.
\newblock {\em Phys. Rev. Lett.}, \textbf{53}:1951--1953, 1984.

\bibitem[Seg68]{Gra}
G.~Segal.
\newblock \textit{Equivariant K-theory}.
\newblock {\em Inst. Hautes \'Etudes Sci. Publ. Math.}, \textbf{34}:129--151,
  1968.

\bibitem[Ska91]{SkanKK}
G.~Skandalis.
\newblock \textit{Kasparov's bivariant $K$-theory and applications}.
\newblock {\em Exposition. Math.}, \textbf{9}:193--250, 1991.

\bibitem[Spa66]{Span}
E.H. Spanier.
\newblock {\em Algebraic topology}.
\newblock McGraw-Hill series in higher mathematics, 1966.

\bibitem[Vas01]{Vas}
S.~Vassout.
\newblock {\em Feuilletages et r\'esidu non commutatif longitudinal}.
\newblock PhD thesis, Univ. Pierre et Marie Curie - Paris VI, 2001.

\bibitem[vE94]{vanElst}
A.~van Elst.
\newblock \textit{Gap labelling theorems for Schrod\"inger operators on the
  square and cubic lattices}.
\newblock {\em Rev. Math. Phys.}, \textbf{6}:319--342, 1994.

\bibitem[Ypm]{Ypma}
F.~Ypma.
\newblock \textit{Quasicrystals, $C^*$-algebras and $K$-theory}.
\newblock \\http://remote.science.uva.nl/~npl/fonger.pdf (2004).

\end{thebibliography}

\bigskip
\noindent
Univ. Blaise Pascal, Clermont-Ferrand, FRANCE\\
\textit{E-mail address :}  haija.moustafa@math.univ-bpclermont.fr\\
\textit{URL :} http://math.univ-bpclermont.fr/~moustafa

\end{document}